\documentclass[review,hidelinks,onefignum,onetabnum]{siamart220329}
\usepackage{stmaryrd}
\usepackage{amsmath, amssymb, amsfonts}
\usepackage{lipsum}
\usepackage{amsfonts}
\usepackage{amssymb}
\usepackage{graphicx}
\usepackage{epstopdf}
\usepackage{algorithmic}
\usepackage{multirow}
\ifpdf
  \DeclareGraphicsExtensions{.eps,.pdf,.png,.jpg}
\else
  \DeclareGraphicsExtensions{.eps}
\fi

\newsiamremark{remark}{Remark}
\newsiamremark{hypothesis}{Hypothesis}
\newsiamremark{example}{Example}
\crefname{hypothesis}{Hypothesis}{Hypotheses}
\newsiamthm{claim}{Claim}

\headers{high order PP SLSV for VP}{X.Y. Zhang, X.F. Cai, and W.X. Cao}

\title{The positivity-preserving high-order semi-Lagrangian spectral volume method  for Vlasov-Poisson equations \thanks{Submitted to the editors DATE.
}
}

\author{Xinyue Zhang\thanks{Laboratory of Mathematics and Complex Systems, School of Mathematical Sciences, Beijing Normal University, Beijing, China  
  (\email{zhang\_xinyue@mail.bnu.edu.cn}).}
  \and Xiaofeng Cai\thanks{Research Center for Mathematics, Advanced Institute of Natural Sciences, Beijing Normal University, Zhuhai 519087, China; Guangdong Provincial/Zhuhai Key Laboratory of Interdisciplinary Research and Application for Data Science, Beijing Normal-Hong Kong Baptist University, Zhuhai 519087, China.
  (\email{xfcai@bnu.edu.cn}).}
\and Waixiang Cao\thanks{Laboratory of Mathematics and Complex Systems, School of Mathematical Sciences, Beijing Normal University, Beijing, China  
 (\email{caowx@bnu.edu.cn}), corresponding author.}
}
\usepackage{amsopn}

\ifpdf
\hypersetup{
  pdftitle={The positivity-preserving high-order semi-Lagrangian spectral volume method  for Vlasov-Poisson equations},
  pdfauthor={X.Y. Zhang, X.F. Cai, and W.X. Cao}
}
\fi

\usepackage{lineno}
\nolinenumbers

\begin{document}

\maketitle

\begin{abstract}
  In this paper,  a novel high order semi-Lagrangian (SL) spectral volume (SV) method is proposed and studied for nonlinear Vlasov–Poisson (VP) simulations via operator splitting. The proposed  algorithm  combines both advantages of semi-Lagrangian and spectral volume approaches, exhibiting strong stability, robustness under large time steps, arbitrary high-order accuracy in space, local mass conservation, and positivity preservation.
Numerical study of  the SLSV method applied to
the  one-dimensional and two-dimensional transport equations, the Vlasov-Poisson system, the classical benchmark problems including Landau damping and two-stream instabilities is conducted,  
confirming the effectiveness,  accuracy, and robustness of  our algorithm 
in addressing complex nonlinear phenomena.  \end{abstract}



\section{Introduction}
The Vlasov equation, in its various formulations (e.g., Vlasov-Maxwell, Vlasov-Darwin, and Vlasov-Poisson), 
is extensively used in both astrophysical contexts \cite{ApplyVlasov1,ApplyVlasov2,ApplyVlasov3} and laboratory environments \cite{ApplyVlasov4,ApplyVlasov5,ApplyVlasov6,ApplyVlasov7,ApplyVlasov8}. 
The Vlasov-Poisson (VP) system, arising from plasma applications, is known as a fundamental model for collisionless plasmas with a negligible magnetic ﬁeld, and it serves as a fundamental framework for modeling the dynamics of collisionless plasma. 

The development of efficient numerical schemes for 
the Vlasov equation is extremely demanding yet challenging. 
One main difficulty in numerical simulation lies in 
its high dimensionality, as the  Vlasov equation  is typically formulated in a six-dimensional phase space (3 spatial and 3 velocity dimensions) plus time, results in substantial computational complexity and cost.  
On the other hand, the  Vlasov equation has many remark physical quantities such as mass, momentum, and energy convervations, as well as  the positivity of the probability density.  A good numerical scheme should be able to mimic as many of
these physical properties as possible, which  is non-trivial and consists of the other difficulty for numerically solving 
Vlasov equation. In addition, 
the linear dependence of the advection velocity on the velocity components in phase space imposes severe time-step restrictions, particularly for particles with moderate to high velocities, further exacerbating computational demands. 

Several numerical approaches have been developed to address the challenges of solving the Vlasov equation.  Particle-in-cell (PIC) methods have long been a standard tool for numerical simulation of the VP system \cite{PICforVPS}, which discretize the distribution function into macro-particles (Lagrangian representation) while solving electromagnetic fields on a mesh (Eulerian representation). This approach inherently ensures positivity and mass conservation, avoids small time-step restrictions by evolving particles in a Lagrangian framework, and leverages standard mesh-based methods for field equations. However, it suffers from the sampling noise of order $\mathcal{O} (1/\sqrt{N} )$, which prevents accurate description of physics of interest in many cases. These limitations highlight the trade-offs in achieving computational efficiency and accuracy. We refer to 
% the classic textbook 
\cite{PICforVPS} for a more detailed review of PIC methods.
    The semi-Lagrangian (SL) method operates by advancing the solution along characteristic lines, similar to PIC methods, while interprets the solution on a fixed grid, akin to the Eulerian approach. Initially, the probability density function (PDF) is defined on a grid. The PDF is then evolved forward in time using Lagrangian dynamics, and the updated PDF is projected back onto the original mesh. This approach combines the advantages of particle methods, such as the absence of small time-step restrictions, with the structured grid framework, enabling high-order accuracy and computational efficiency. 
  SL methods are unconditionally stable at the algorithmic level. This allows larger time steps, enhancing computational efficiency. As proven in \cite{YangCaiQiuSLDG2020}, the semi-Lagrangian discontinuous Galerkin (SLDG) method achieves unconditional stability, optimal convergence, and superconvergence in long-time simulations.  Due to these advantageous properties, SL methods have garnered substantial interest in the plasma simulation community, as evidenced by recent advancements and applications in this field \cite{SLforVPS1,SLforVPS2}. The integration of the SL method with the discontinuous Galerkin (DG) approach has also been extensively studied (see \cite{Cai2017SLDG,Cai2018SLDG,Cai2022SLDG,2011splitSLDG,ECSLDG2023,ifeSLDG}).

The main purpose of the current study is to propose a novel high-order, grid-based semi-Lagrangian spectral volume (SLSV) method method for solving the  Vlasov-Poisson equation. We concentrate on the 1 spatial and 1 velocity dimensions plus time equation.
Specifically, 
our discretization framework is based on high-order spectral volume (SV) scheme coupled with an operator splitting semi-Lagrangian 
method. 
The SV method was initially formulated and subsequently developed for hyperbolic equations by Wang and his collaborators \cite{SV2,SV3,SV1}, and is viewed as a generalization of the classical Godunov finite volume method \cite{Riemann-exact,Riemann-appro}. The SV method exhibits several advantageous properties, including high-order accuracy, compact stencils, and geometric flexibility (applicable to unstructured grids). Notably, as the SV method preserves conservation laws on finer meshes, it may achieve higher resolution for discontinuities compared to other high-order methods (see \cite{SunWang2004}). 
Over the past decades, the SV method has been extensively applied to solving various partial differential equations (PDEs), such as the shallow water wave equation \cite{SVShallow}, the Navier-Stokes equations \cite{SVNS1,SVNS2}, and electromagnetic field equations \cite{SV1}, among others. A rigorous mathematical analysis of the SV method, including $L^2$ stability, accuracy, error estimates, and superconvergence, was studied in \cite{SV-Cao} within the framework of the Petrov-Galerkin method. It was also demonstrated in \cite{SV-Cao} that a specific class of SV schemes is equivalent to upwind DG schemes when applied to linear constant hyperbolic equations. 
In this work, as an extension of the SV method to the Vlasov-Poisson system, we apply the SV method to plasma-related problems. 
We demonstrate that this methodology offers a promising approach, delivering highly accurate numerical results, preserves
 local mass, energy, and positivity while maintaining relatively low computational cost.

The structure of this paper is organized as follows. In Section 2, we provide 
 a concise introduction to the spatial discretization framework of the SV method and propose 
 a SLSV scheme  for one-dimensional linear transport equation. 
 Section 3 elaborates on the proposed SLSV method applied  to two-dimensional transport equations with 
 the operator splitting.   In Section 4, we adopt the SLSV method for VP system and a detailed algorithm is provided. 
 Several numerical experiments are presented in Section 5 to demonstrate the accuracy, 
efficiency, and conservation properties of the methods. Finally, in  Section 6 we give some concluding remarks.

\section{The SLSV method for one-dimensional advection equations}\label{Sec:1D}

In this section, we will present the semi-Lagrangian spectral volume (SLSV) method for the one-dimensional transport equation and introduce the theoretical foundations that underpin the SLSV scheme.

Consider the following one-dimensional transport equation 
\begin{equation}\label{transporteq}
	u_t+(a(x,t)u)_x=0,\quad x\in \Omega:=[x_a,x_b], \ u(x,0)=u_0(x),
\end{equation}
 with periodic boundary condition.  Here $u_0$ and  the velocity field
 $a(x,t)$  are smooth.

\subsection{The standard SV method}

Let the computational domain $\Omega $ is divided into $N$ nonoverlapping elements, i.e., 
$x_a=x_{\frac 1 2}<x_{\frac 3 2} < \cdots <x_{N+\frac 1 2}=x_b$. 
Denote the $i$th element, called an SV(spectral volume), its center and size as $I_i=[x_{i-\frac 1 2},x_{i+\frac 1 2}]$, $x_i=(x_{i+\frac 1 2}+x_{i-\frac 1 2} )/ 2$ and $h_i = x_{i+\frac 1 2}-x_{i-\frac 1 2}$, respectively. 

Denote by $T_i,~i=1,2,\cdots, N$  the the linear transformation  mapping  the reference volume $[-1,1]$ to each SV $I_i$. That is, 
    \begin{equation}\label{affine transformation}
        x = T_is = \frac {h_i}{2} s + x_i,\quad s \in [-1,1],\ x\in I_i.
    \end{equation}
For any positive integer $l$,  let $ \mathcal{Z}_l = {1,2,\cdots,l},~ \mathcal{Z}_l^0={0,1,2,\cdots,l}$, and  $-1=s_0<s_1<\cdots<s_k<s_{k+1}=1$ be $k+2$ points in the reference volume $[-1,1]$.  Define 
$$
x_{i,p} = T_is_p, \quad i\in \mathcal{Z}_N , ~p \in \mathcal{Z}_{k+1}^0. 
$$
Here $x_{i,0}=x_{i-\frac 12}, x_{i,k+1}=x_{i+\frac 12}$.
% , and $h_i=x_{i+\frac 12}-x_{i-\frac 12}$ 
These $k$+$2$ points partition the SV $I_i$ into $k+1$ subintervals $I_{i,p} = [x_{i,p},x_{i,p+1}], i\in \mathcal{Z}_N, p\in\mathcal{Z}_{k+1}$, called CV(control volume).
 All the CVs constitute the dual partition of $\Omega$.  
 Define the (discontinuous) finite element space  by 
    \begin{equation*}\label{trailespace}
        \mathcal{U}_{h}=\mathcal{U}_{h}^{k}=\{v \in L^{2}(\Omega):\left.v\right|_{I_i} \in  \mathcal{P}^k,  i \in \mathcal{Z}_N\},
    \end{equation*}
where $\mathcal{P}^k$  denotes the finite element space  of polynomials  with degree not greater than $k$.

The  idea of the SV method is based on the local conservation on each CV.  To be more precise, 
The standard SV method  for \eqref{transporteq}  reads as: find a $u_h(\cdot, t) \in \mathcal{U}_h$ such that 
\[
  \int_{I_{i,p}}\partial_t u_h dx+\widehat {au_h}(x_{i,p+1})-\widehat {au_h}(x_{i,p})=0, 
\]
where $\widehat {au_h}$ denotes the numerical flux of $au_h$, which  should be carefully 
designed to ensure the stability of SV schemes.  
Well-known monotone fluxes include the Lax-Friedrichs flux
upwind flux, Lax-Friedrichs flux, the Godunov flux and HLLC flux, and so on  (see, e.g.,\cite{Riemann-appro,LeVeque}). 
Noticing that different choices of $\{s_i\}_{i=1}^{k}$ leads to different SV scheme, which may have effect on the stability of the SV scheme. Throughout of this paper,  we choose Gauss points (which are shown to be most robust and stable numerically)
 to construct CVs.

\subsection{The semi-Lagrangian SV method}

We note that \eqref{transporteq}  admots 
 the characteristic equation
\[
  \frac{dx(t)}{dt} = a(x(t), t).
\]
 Consider the adjoint problem of \eqref{transporteq} 
 \[
    \varphi_t+a(x,t)\varphi_x=0,\ \ t\in [t^n,t^{n+1}],\  \varphi(x,t^{n+1})=\varphi_0(x). 
 \]
  Then the solution stays constant
   along a 
 characteristic trajectory determined by \eqref{transporteq}, and thus (\cite{limiter2})
 \[
     \frac{d}{dt} \int_{\tilde I_i(t)} u(x,t)\varphi(x,t) dx=0, 
 \] 
  where $\tilde I_i(t)$ is a dynamic interval bounded by the characteristics from cell boundaries of $I_i$ at $t = t^{n+1}$. 
  Then the weak formulation of semi-Lagrangian reads 
  \begin{equation}\label{weak-form}
     \int_{ I_i} u(x,t^{n+1})\varphi_0 dx= \int_{I_i^*} u(x,t^{n})\varphi(x,t) dx, 
 \end{equation}
   where $I_i^*=[x_{i-\frac 12}^*,x_{i+\frac 12}^*]$ with  $x_{i\pm \frac 12}^*$
being the feet of trajectory emanating from $(x_{i\pm \frac 12},t^{n+1})$ at time level $t^n$.

\begin{figure}[!h]
    \centering
    \includegraphics[scale=0.15]{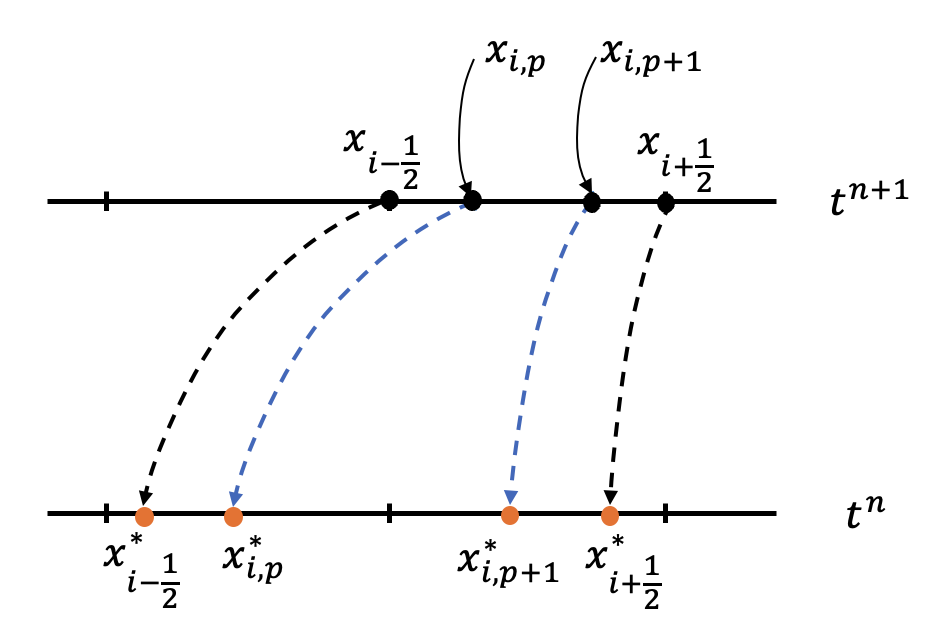}
    \includegraphics[scale=0.15]{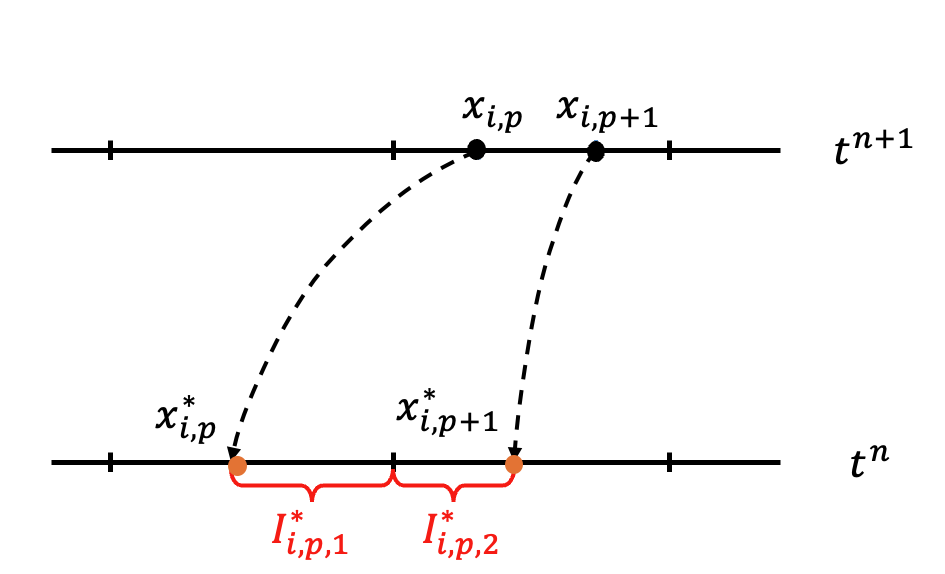}
\caption{\label{Fig:upstream} Schematic illustration for one-dimensional SLSV schemes.}
\end{figure}

The SLSV method is based on characteristic Galerkin weak formulation in \eqref{weak-form} while incorporates
sub-interval-level conservation property of SV method. To be more precise,  
 The SLSV method for  \eqref{transporteq}  reads as: find a $u_h\in \mathcal{U}_{h}$ such that 
 \begin{equation}\label{RHSSV}
	\int_{I_{i,p}} u_h\left(x, t^{n+1}\right) d x =\int_{I_{i,p}^*} u_h\left(x, t^n\right) d x, \quad u \in \mathcal{U}_{h},
\end{equation}
where ${I}^*_{i,p}(t)$ is a dynamic moving cell, emanating from a CV $I_{i,p}$ at $t^{n+1}$ (see Figure \ref{Fig:upstream} (left)).
We note that to update the numerical solution $u_h$ at time $t^{n+1}$ , we need to evaluate the integral on
the right-hand side (RHS) of \eqref{RHSSV}, i.e.,  the information over a CV from the solution at $t^n$ is required.

Now we divide procedure of SLSV algorithm into the following steps:

\begin{itemize}
\item[i.] Choose $k+2$ sub-nodes in each SV $I_i, i \in \mathcal{Z}_N$ and then 
construct CVs $I_{i,p}, p \in \mathcal{Z}_{k+1}^0$. 

\item[ii.] Trace characteristics backward in time to $t^n$ for each subnode $x_{i,p},~i\in \mathcal{Z}_N,~p\in \mathcal{Z}_{k+1}^0$ and then obtain the new coordinate $x_{i,p}^*$, see Figure \ref{Fig:upstream} (left), by numerically solving the following final-value problem (trajectory equation):
    $$
    \frac{d x(t)}{d t}=a(x(t), t) \text { with } x\left(t^{n+1}\right)=x_{i,p},\quad i\in \mathcal{Z}_N,~p\in \mathcal{Z}_{k+1}^0,
    $$
with a high order numerical integrator such as a fifth order Runge-Kutta method. 

\item[iii.] Detect intervals/sub-intervals within the upstream CV, $I_{i,p}^{\star}=\displaystyle\bigcup_l I_{i,p,l}^{\star}$, which are the intersections between $I_{i,p}^{\star}$ and the grid elements ($l$ is the index for sub-interval). For instance, see Figure \ref{Fig:upstream}(right), there are two sub-intervals: $I_{i,p, 1}^{\star}=\left[x_{i,p}^{\star}, x_{i-\frac{1}{2}}\right]$ and $I_{i,p, 2}^{\star}=\left[x_{i-\frac{1}{2}}, x_{i,p+1}^{\star}\right]$. 

\item[iv.]  Update the numerical solution at $t^{n+1}$ by using \eqref{RHSSV}  with 
  the right hand side (RHS) of \eqref{RHSSV}  approximated by
\begin{equation}\label{RHSSSV}
	\int_{I_{i,p}^*} u_h\left(x, t^n\right) d x = \sum_{l } \int_{I_{i,p,l}^*} u_h\left(x, t^n\right) d x.
\end{equation}
\end{itemize}

We note that  the integration in \eqref{RHSSSV}  can be calculated exactly 
since  $u_h$ is a polynomial of degree $k$.
 
 \section{ SLSV with operator splitting for 2D transport equation}\label{Sec:2d transport}

   In this section, we discuss the  SLSV  method applied to the two-dimensional  transport equation
    via operator splitting.  Consider the following 2D transport equation
    \begin{equation}{\label{2D transport}}
        u_t+(a(x, y, t) u)_x+(b(x, y, t) u)_y=0,\quad (x, y) \in \Omega,
    \end{equation}
with a proper initial condition $u(x, y, 0)=u_0(x, y)$. Here, $(a(x, y, t), b(x, y, t))$ is a prescribed velocity field depending on time and space. Boundary conditions are periodic for simplicity. 

Let $\Omega= [x_a,x_b]\times[y_a,y_b]$,   partitioned into $N_x$ and $N_y$  nonoverlapping elements in $x$- direction and $y$- direction, respectively. 
That is,  
    $$x_a=x_{\frac 1 2}< \cdots<x_{i + \frac 1 2} < \cdots <x_{N_x+\frac 1 2}=x_b,$$ 
    $$y_a=y_{\frac 1 2}< \cdots<y_{j + \frac 1 2} < \cdots <y_{N_y+\frac 1 2}=y_b.$$
Denoted the element $A_{i j}=\tau_i^x \times \tau_j^y = \left[x_{i-\frac{1}{2}}, x_{i+\frac{1}{2}}\right] \times\left[y_{j-\frac{1}{2}}, y_{j+\frac{1}{2}}\right]$ as a SV, and its center and size in $x$- direction, $y$- direction as
$x_i=\frac{x_{i+\frac 1 2}+x_{i-\frac 1 2} } {2},y_j=\frac{y_{j+\frac 1 2}+y_{j-\frac 1 2} } {2},$
$h^x_i = x_{i+\frac 1 2}-x_{i-\frac 1 2},h^y_j = y_{j+\frac 1 2}-y_{j-\frac 1 2}, i \in \mathcal{Z}_{N_x}, j \in \mathcal{Z}_{N_y},$ 
respectively. 
 
 We next present the SLSV scheme for the two-dimensional  transport equation \eqref{2D transport}, inspired by the 
 idea of Strang operator splitting proposed by Cheng and Knorr in \cite{Split}.  In the two-dimensional splitting setting, 
 we adopt the piecewise ${\mathcal Q}^k$ tensor-product polynomial space, 
 and each SV $A_{i j}$ is partitioned into $(k+1)^2$ CVs by using the $k$ Gauss points along the $x$- and $y$- directions, 
 respectively. 
  
  Now we are ready to present a general framework of the splitting SLSV scheme. 
 First,  we split  \eqref{2D transport} into two 1D advection problems in both $x$- and $y$-directions:
    \begin{align}
        \label{2d transport split1} u_t+(a(x, y, t) u)_x=0,\\ 
        \label{2d transport split2} u_t+(b(x, y, t) u)_y=0.
    \end{align}
The primary advantage of this splitting approach lies in its ability to facilitate the direct application of the one-dimensional SLSV scheme to two-dimensional problems, thereby ensuring high computational efficiency. 
 Then the split equations \eqref{2d transport split1}  and \eqref{2d transport split2}  are solved based on   the 1D SLSV formulation 
 via Strang operator splitting at the temporal instance $t^n$ over a time step $\Delta t$ with $\Delta t=t^{n+1}-t^n$. The specific steps of the process are described in Algorithm \ref{alg:2d transport} provided below.

\begin{algorithm}
\caption{SLSV scheme via Strang splitting for 2D passive-transport equation.}
\label{alg:2d transport}
\begin{algorithmic}[1]
\STATE{\(\frac{1}{2} \Delta t\) step on \(u_t+(a(x, y, t) u)_x=0\).}
\STATE{\(\Delta t\) step on \(u_t+(b(x, y, t) u)_y=0\). }
\STATE{\(\frac{1}{2} \Delta t\) step on \(u_t+(a(x, y, t) u)_x=0\).}
\end{algorithmic}
\end{algorithm}

Next,   we present a detailed procedure of the SLSV scheme  applied to the quasi-1D problem \eqref{2d transport split1} in Step 1 of Algorithm \ref{alg:2d transport}. 
To effectively reduce the original formulation to a purely one-dimensional problem, the equation \eqref{2d transport split1} is evaluated along the Gaussian lines aligned with the $y$-direction within each Cartesian grid cell, i.e., 
    \begin{equation*}\label{gaussline}
        y = y_{j,m},\quad y_{j,m}:=T_j \eta_m, \quad  m\in \mathcal{Z}_{k+1},
    \end{equation*}
where $\eta_m$ are the roots of the $(k+1)$th degree Legendre polynomial in $[-1,1]$,  and 
$T_j$ denotes the $y$- direction affine transformation from $[-1,1]$ to $[y_{j-\frac 12}, y_{j+\frac 12}]$, i.e.,  
\begin{equation}\label{y-affine transformation}
    y = T_js = \frac {h_j^y}{2} s + y_j,\quad s \in [-1,1].
\end{equation}
Along each Gauss line, we adopt the 1D SLSV scheme to solve the following  constant coefficient advection:
    \begin{equation}\label{quasi-1D equation}
        u_t+(a(x, y_{j,m}, t) u)_x=0.
    \end{equation}

   Further details pertaining to the solution of the quasi-1D problem \eqref{2d transport split1} in the Step 1 of Algorithm \ref{alg:2d transport} will be demonstrated in Algorithm \ref{alg:quasi-1d}.    The Step 2 and Step 3 of
     Algorithm \ref{alg:2d transport}  follow the same procedure. 
     We note  that in Step 2, the positions of $x$ and $y$ are interchanged. This allows us to obtain the numerical solution of $u$ from $t^n$ to $t^{n+1}$.

\begin{algorithm}
    \caption{The proposed SLSV algorithm for solving quasi-1D advection equations of the form $u_t+(a(x, y, t) u)_x=0$ with a time step $\frac {\Delta t} {2}$.}
    \label{alg:quasi-1d}
    \begin{algorithmic}[1]
    \STATE {Start with the current solution:
    		$$\left.u_h\left( x,y,t^n\right)\right|_{A_{i j}}:=\displaystyle\sum_{\ell=1}^{(k+1)^2} U_{i j}^{(\ell)} \varphi^{(\ell)}(\xi, \eta).$$
    		The canonical variables $(\xi, \eta)$ are linearly related to the variables $(x, y)$ via affine transformation $T_i$ in \eqref{affine transformation} 			and $T_j$ in \eqref{y-affine transformation}. The term $ U_{i j}^{(\ell)}$ represents the coefficient corresponding to each basis function $				\varphi^{(\ell)}(\xi, \eta)$ in the SV $A_{i,j}$,  at $t^n$.}
     \STATE {In each spectral volume $A_{i,j}$, calculate the values of $u_h$
    		along the Gaussian line $y=y_{j,m}$ at $t^n$, that is,  
    		\begin{equation*}\label{1d quasi assumption}
        			\left.u_h\left( x,y_{j,m},t^n\right)\right|_{A_{i j}}=\sum_{\ell=1}^{k+1} U_{i j m}^{(\ell)} \Phi^{(\ell)}\left(\xi\right).
   		 \end{equation*}
    		Here  $ {U}_{i j m}^{(\ell)}$ represents the coefficient corresponding to each basis function $\Phi^{(\ell)}\left(\xi\right)$ 
    		along the $x$ direct. }
     \STATE {Solve the equation \eqref{quasi-1D equation} by using the 1D SLSV scheme \eqref{weak-form} with the given values of 							$u_h\left( x,y_{j,m},t^n\right)$ at time $t^n$, and   derive the numerical  solution for each Gaussian line $y=y_{j,m}$ at $t^{n}+\frac {\Delta t}{2}$, 		see Figure \ref{Fig:oper_split}(b), which is
    			\begin{equation}\label{1d quasi solusion}
        				\left.u_h^{new}\left( x,y_{j,m}\right)\right|_{A_{i j}}=\sum_{\ell=1}^{k+1} \widehat{U}_{i j m}^{(\ell)} \Phi^{(\ell)}\left(\xi\right).
   			 \end{equation}
    		Here the term $ \widehat{U}_{i j m}^{(\ell)}$ denotes the coefficient associated with the basis function $\Phi^{(\ell)}\left(\xi\right)$.}
     \STATE  {Obtain  the values of $u_h^{new}$ at the  $(k+1)^2$ Gauss points from \eqref{1d quasi solusion}, and then 
    		reconstruct $u_h^{new}$ at time $t^{n}+\frac{\Delta t}{2}$ by using the Gauss Lagrange interpolation. }
    \end{algorithmic}
\end{algorithm}

    \begin{figure}[!h]
        \centering
        \includegraphics[scale=0.5]{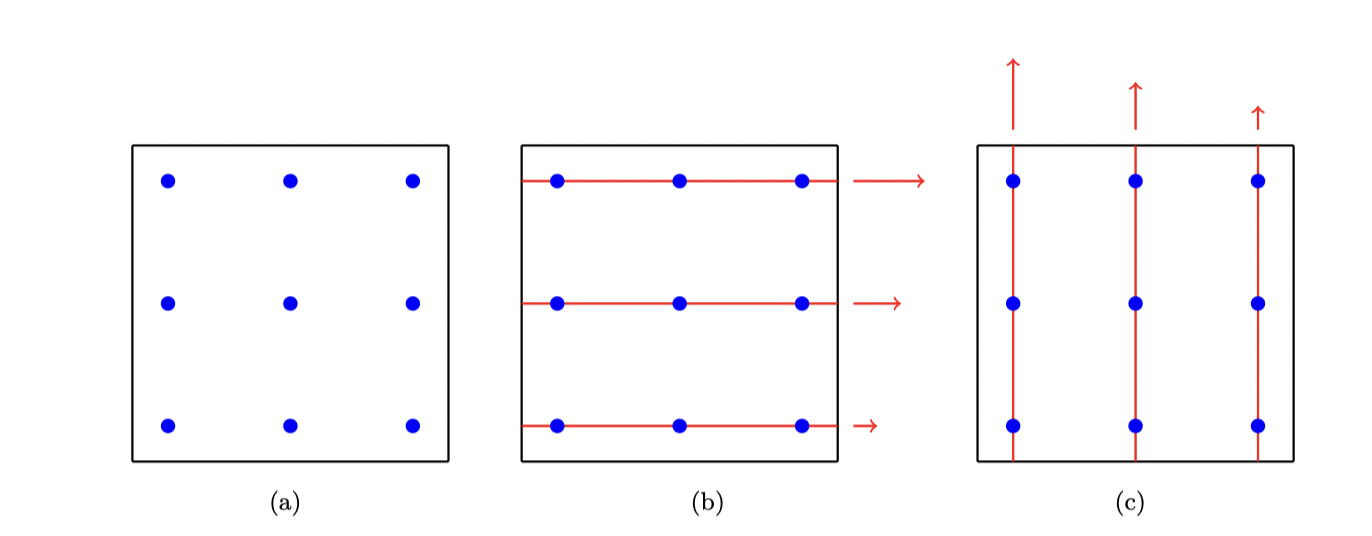}
    \caption{\label{Fig:oper_split} Schematic illustration of the 2D SLSV scheme via Strang splitting with $k = 2$.}
    \end{figure}

     \begin{remark}
        The splitting 2D SLSV formulation inherits many desired properties from the base 1D formulation, such as the high-order accuracy in space, the unconditional stability, and the mass conservation. Meanwhile, a splitting error in time is introduced. For the Strang operator splitting, the error is of second order. Higher order splitting methods can be constructed in the spirit of composition methods \cite{HighOrder-split1,HighOrder-split2}, while the number of intermediate stages accumulates substantially; hence, the computational cost, increases exponentially with the order of the splitting method. For example, a fourth-order splitting SLDG method is developed in \cite{HighOrder-split2} for solving the Vlasov-Poisson system.
     \end{remark}

\section{SLSV method for Vlasov-Poisson system}\label{Sec:2d VPS}

In this section, we discuss the application of the SLSV method for 
 the Vlasov-Poisson  system 
\begin{equation}\label{VPS}
    \frac{\partial f}{\partial t}+\mathbf{v} \cdot \nabla_{\mathbf{x}} f+\mathbf{E}(\mathbf{x}, t) \cdot \nabla_{\mathbf{v}} f=0,
\end{equation}
and
\begin{equation}\label{VPS Electric}
    \mathbf{E}(\mathbf{x}, t)=-\nabla_{\mathbf{x}} \phi(\mathbf{x}, t), \quad-\Delta_{\mathbf{x}} \phi(\mathbf{x}, t)=\rho(\mathbf{x}, t),
\end{equation}
where $\mathbf{x}$ and $\mathbf{v}$ are coordinates in phase space $(\mathbf{x}, \mathbf{v}) \in \mathbb{R}^3 \times \mathbb{R}^3, \mathbf{E}$ is the electric field, $\phi$ is the self-consistent electrostatic potential and $f(\mathbf{x}, \mathbf{v},t)$ is probability distribution function (PDF) which describes the probability of finding a particle with velocity $\mathbf{v}$ at position $\mathbf{x}$ at time $t$. The probability distribution function couples to the long range fields via the charge density, $\rho(t, \mathbf{x})=\int_{\mathbb{R}^3} f(\mathbf{x}, \mathbf{v},t) d \mathbf{v}-1$, where we take the limit of uniformly distributed infinitely massive ions in the background. Equations \eqref{VPS} and \eqref{VPS Electric} have been nondimensionalized so that all physical constants are one.

  For simplicity, 
we consider the VP system \eqref{VPS} with the one-dimensional (1D) physical space and 1D velocity space on the two dimensional domain $\Omega=\Omega_x \times \Omega_v$:
    \begin{align}
        \label{1+1 VPS f} & f_{t}+v f_{x}+E(x,t) f_{v}=0, \\
        \label{1+1 VPS electric} & E_x=\rho(x,t)-\rho_0.
    \end{align}
where $\Omega_x$ is a bounded domain with periodic boundary conditions. As for $v$- direction, set  $\Omega_v= [ -v_{\max }, v_{\max }] $ with $v_{\max }$ chosen large enough so that a zero boundary condition is reasonably imposed.
The total and background densities are
    \begin{equation*}\label{back densities}
        \rho(x, t):=\int_{-\infty}^{\infty} f(x, v, t) d v, \quad \text { and } \rho_0:=\frac{1}{2 v_{\max }} \int_{-v_{\max }}^{v_{\max }} \rho(x, t) d x. 
    \end{equation*}
Note that $\rho_0$ is a constant in time, due to conservation of mass on the periodic domain  $\left[-v_{\max }, v_{\max }\right]$.

 Although the PDF is not itself a physical observable, its moments represent various physically observable quantities:

$$
\begin{aligned}
& \rho(x,t):= \displaystyle\int f(x,v,t) d v, \quad \text { (mass density)}, \\
& \rho u(x,t):= \displaystyle\int v f(x,v,t) d v, \quad \text { (momentum density), } \\
& \mathcal{E}(x,t):= \displaystyle \frac{1}{2} \int v^2 f(x,v,t) d v, \quad \text { (energy density). }
\end{aligned}
$$

The Vlasov-Poisson system contains an infinite number of quantities that are conserved in time, 
which  can be used to measure the performance of  a numerical algorithm. 
We list four conservation quantities that will be used in our numerical experiment to 
test the efficiency of our scheme.
\[
\begin{aligned}
& \|f\|_{L_1}:=\iint |f(x,v,t)| d v d x, \\
& \|f\|_{L_2}:=\left(\iint f^2(x,v,t) d v d x\right)^{\frac{1}{2}}, \\
& \text { Total energy }:=\iint v^2 f(x,v,t) d v d x+ \int \mathbf{E}^2 d x,\\
& \text { Entropy }:=-\iint \ f(x,v,t) \log (f(x,v,t)) d v d x.
\end{aligned}
\]

Below we describe in detail how SLSV can be applied to the Vlasov-Poisson system.
Based on the SLSV method via Strang operator splitting applied to the two-dimensional transport equation, we split the equation \eqref{1+1 VPS f} into two 1D advection problems based on the quadrature nodes in both $x$- and $v$- directions:
    \begin{align*}
        & f_{t}+v f_{x}=0, \\
        & f_{t}+E(x,t) f_{v}=0.
    \end{align*}
Built upon the 1D SLSV scheme, evolve the split equations via operator splitting at time $t^n$ over a time step $\Delta t$ : 
    \begin{algorithm}
        \caption{SLSV scheme via operator splitting for Vlasov-Poisson system.}
        \label{alg:2d VPS}
        \begin{algorithmic}[1]
        \STATE {\(\frac{1}{2} \Delta t\) step on \(f_t + \mathbf{v} \cdot f_x = 0\).}
        \STATE {Solve \(\nabla^2 \phi = \rho^{n+\frac 1 2} - \rho_0\), and compute \(\mathbf{E}^{n+\frac 1 2} = \nabla \phi\).}
        \STATE {\(\Delta t\) step on \(f_t + \mathbf{E}^{n+\frac 1 2} \cdot f_x = 0\).}
        \STATE {\(\frac{1}{2} \Delta t\) step on \(f_t + \mathbf{v} \cdot f_x = 0\).}
        \end{algorithmic}
    \end{algorithm}

The SLSV scheme for splitting equation in Algorithm \ref{alg:quasi-1d} can then be inserted into Steps $1,3$, and $4$ of Algorithm \ref{alg:2d VPS} for solving the $1+1$ dimensional Vlasov-Poisson system. We note that for Steps $1$ and $4$, $a(x,y,t)=v$ in \eqref{2d transport split2}; while for Step $3$, $b(x,y,t)=E^{n+\frac{1}{2}}(x)$. 
For  Step $2$, any numerical scheme capable of solving the Poisson equation with high precision can be employed. Given that the test function space is a discontinuous finite element space, the local discontinuous Galerkin (LDG) method is selected in this study due to its superior compatibility and adaptability.

We end with this section the discuss of positivity-preserving property.  Note that our algorithm satisfies the local conservation 
property and thus cell average is preserved at all time, which indicates that  the cell average  at $t^{n+1}$ remains positive, provided that the numerical solution $f^n$ at $t^n$ is positive. Furthermore, by 
defining 
\[
        \tilde{f}(x, v)=\theta\left(f\left(x, v, t^n\right)-\bar{f}\right)+\bar{f}, \theta=\min \left\{\left|\frac{\bar{f}}{m^{\prime}-\bar{f}}\right|, 1\right\}, 
\]
where $\bar{f}$ is the cell average of the numerical solution and $m^{\prime}$ is the minimum value of $f\left(x, v, t^n\right)$ over $A_{ij}$, 
  we obtain a point-wise positive function at each element $A_{i,j}$.   We would like to point out that 
  the new defined function, referred as positivity-preserving (PP) limiter (see, e.g.,\cite{PPlimiter}) is of high order approximation 
  and do not have effect on the error accuracy. We also refer to  \cite{limiter1,limiter2} for the detailed discuss of 
  the positivity-preserving limiter.

\section{Numerical experiments}

In this section, a sequence of numerical experiments are conducted to test the performance, accuracy, conservation 
property of the proposed SLSV schemes. 
The  one-dimensional linear advection equation is first studied in Section \ref{Sec:num-1d}
 which is solved by using the  SLSV scheme with large time step. 
Numerical  results show that the amplification time step does not affect the optimal convergence order, indicating the high performance and efficiency of our algorithm. 
In Section \ref{Sec:num-2d}, we test the order of accuracy of the SLSV method with dimensional splitting for two-dimensional linear passive transport problems. The performance on shape preservation of the SLSV scheme is also demonstrated. 
 Section \ref{Sec:VPS} is dedicated to the study of SLSV for the Vlasov-Poisson  system, including 
 the forced Vlasov-Poisson equation with known exact solutions, weak and strong Landau damping, and two-stream instability.

In all of our numerical tests,  uniform meshes are obtained by equally dividing the 
computational domain into $N$ ( for 1D) or $N\times N$ (for 2D) spectral volumes. 
 We take the time step size
 as  $\Delta t=\frac{CFL}{\frac{v_{\max }}{\vartriangle x }+\frac{\max (\left\lvert E\right\rvert )}{\vartriangle v}},$ 
where CFL is specified for different runs.   

\subsection{1D Linear Passive‑Transport Problems}\label{Sec:num-1d}
    \begin{example}(\textbf{1D Linear equation})
    \end{example}
    We  consider the following 1D transport equation
        \begin{equation*}\label{1Dlinear}
            u_t+u_x=0,~(x,t)\in[0,2\pi]\times (0,20],
        \end{equation*}
    with smooth initial data $u(x)=\sin(x)$ and periodic boundary condition. The exact solution is 
    $u(x,t)=\sin(x-t)$. 
    % The meshes are  obtained by equally dividing $[0,2\pi]$ into $N$ elements. 
    
     Listed in Table \ref{Table:SV1Dlinear} are
    the $L^2$ and $L^{\infty}$ errors and the associated convergence rates calculated from $ {\mathcal P}^k, k=1,2,3$, SLSV  method  with $\Delta t = 0.4h$ and $\Delta t = 2.4h$. Optimal convergence rates 
      are observed. In addition, the errors  for $\Delta t = 2.4h$ are either
    comparable or slightly smaller than those  for $\Delta t = 0.4h$, indicating that 
    the large CFL does not  affect the approximation accuracy. 

        \begin{table*}[!ht]
            \centering
            \footnotesize
            \caption{\label{Table:SV1Dlinear}1D linear equation: errors and convergence orders of $L^2$ and $L^{\infty}$  at $T$=20.} 
            \vspace{0.4em} \centering%  xiaowu
            \begin{tabular}{c|c|cccc|cccc}%{\textwidth}
            %\toprule[0.001pt]
            \hline
            % \multirow{2}{*}{k}   &\multirow{2}{*}{N}   & \multicolumn{4}{l}{$P^k$ SLSV with $CFL=0.4$}  & \multicolumn{4}{l}{$P^k$ SLSV with $CFL=2.4$} \\
            % \cmidrule(r){3-6} \cmidrule(r){7-10} 
            \multirow{2}{*}{k}   &\multirow{2}{*}{N}   & \multicolumn{4}{c|}{$CFL=0.4$}  & \multicolumn{4}{c}{$CFL=2.4$} \\
            \cline{3-10}
            & &$L^2$ &Order&$L^{\infty}$& Order&$L^2$ &Order&$L^{\infty}$& Order\\
            %\midrule[0.001pt]
            \hline
            ~ &20 &7.72E-03 &- &1.37E-02&- &3.77E-03 &-&1.21E-02&-\\
            ~ &40 &1.96E-03 &1.98 &3.78E-03&1.86 &9.30E-04 &2.02    &3.01E-03 &2.00 \\
            1 &80 &4.92E-04 &1.99 &9.76E-04&1.95 &2.38E-04 &1.96    &7.77E-04 &1.96 \\
            ~ &160 &1.21E-04 &2.02 &2.36E-04&2.05 &4.67E-05 &2.35    &1.37E-04 &2.50 \\
            ~ &320 &3.08E-05 &1.98 &6.23E-05&1.92 &1.16E-05 &2.01    &3.36E-05 &2.03 \\
            \hline
            ~ &20 &3.18E-04 &- &4.70E-04&- &9.52E-05 &-&1.94E-04&-\\
            ~ &40 &4.00E-05 &2.99 &5.82E-05&3.01 &1.23E-05 &2.95    &2.52E-05 &2.94 \\
            2 &80 &4.98E-06 &3.00 &7.26E-06&3.00 &1.49E-06 &3.05    &2.81E-06 &3.16 \\
            ~ &160 &6.23E-07 &3.00 &9.33E-07&2.96 &1.73E-07 &3.10    &4.76E-07 &2.56 \\
            ~ &320 &7.80E-08 &3.00 &1.13E-07&3.05 &2.13E-08 &3.02    &5.41E-08 &3.14 \\
            \hline
            ~ &20	&4.44E-06	&-	&7.76E-06	&- &1.99E-06&-	&4.08E-06	&- \\
            ~ &40	&2.81E-07	&3.98	&4.15E-07	&4.22 &1.21E-07	&4.05 &2.22E-07	&4.20\\
            3 &80	&1.78E-08	 &3.98 	&2.75E-08	&3.92 &7.99E-09	&3.91 &1.66E-08	&3.74\\
            ~ &160	&1.09E-09	 &4.02 	&1.80E-09	&3.94 &3.85E-10	&4.38 &8.96E-10	&4.21\\
            ~ &320	&7.01E-11	 &3.96 	&1.09E-10	&4.05 &2.46E-11	&3.97 &6.10E-11	&3.88\\
            \hline
            \end{tabular}
            %\vspace{\baselineskip}
        \end{table*}

\subsection{2D Linear Passive‑Transport Problems}\label{Sec:num-2d}
    \begin{example}(2D linear equation) 
    \end{example}
    Consider
        \begin{equation*}\label{2d linear equation}
            \left\{\begin{array}{l}
            u_t+u_x+u_y=0, (x,y) \in[0,2 \pi] \times [0,2 \pi], \\
            u(x, y, 0)=\sin (x+y),
            \end{array}\right.
        \end{equation*}
    with the periodic boundary conditions in both $x$ and $y$ directions. The exact solution is
        $$
        u(x, y, t)=\sin (x+y-2t).
        $$
 Table \ref{Tab: 2D linear} presents
  the $L^2$ and $L^{\infty}$ errors and the associated orders of accuracy when the ${\mathcal Q}^k,k=1,2,$ SLSV methods are 
  used with $\Delta t=0.5 h$ and $\Delta t=2.5 h$. Again 
  expected orders of accuracy and comparable errors
  are observed for both the two different CFLs.

        \begin{table*}[!ht]
            \centering
            \footnotesize
            \caption{\label{Tab: 2D linear}2D linear equation: errors and convergence orders of $L^2$ and $L^{\infty}$  at $T=\pi$.} 
            \vspace{0.4em} \centering%  xiaowu
            \begin{tabular}{c|c|cccc|cccc}%{\textwidth}
            %\toprule[0.001pt]
            \hline
            % \multirow{2}{*}{k}   &\multirow{2}{*}{N}   & \multicolumn{4}{l}{$P^k$ SLSV with $CFL=0.4$}  & \multicolumn{4}{l}{$P^k$ SLSV with $CFL=2.4$} \\
            % \cmidrule(r){3-6} \cmidrule(r){7-10} 
            \multirow{2}{*}{$k$}   &\multirow{2}{*}{$N^2$}   & \multicolumn{4}{c|}{ $CFL=0.5$}  & \multicolumn{4}{c}{$CFL=2.5$} \\
            \cline{3-10}
            & &$L^2$ &Order&$L^{\infty}$& Order&$L^2$ &Order&$L^{\infty}$& Order\\
            %\midrule[0.001pt]
            \hline
            ~ &$20^2$ &     1.78E-02 & - &     3.24E-02 & -  &     1.33E-02 &     - &     1.71E-02 &     -   \\
            ~ &$40^2$ &     4.53E-03 &     1.97 &     8.33E-03 &     1.96 & 2.91E-03 &     2.19 &     3.86E-03 &     2.15 \\
            1 &$80^2$ &     1.14E-03 &     1.99 &     2.11E-03 &     1.98 &     6.66E-04 &     2.13 &     9.41E-04 &     2.04 \\
            ~ &$160^2 $&     2.86E-04 &     2.00 &     5.30E-04 &     1.99  &    1.54E-04 &     2.12 &     2.17E-04 &     2.12 \\
            ~ &$320^2$ &     7.16E-05 &     2.00 &     1.33E-04 &     2.00  &     4.98E-05 &     1.99 &     1.01E-04 &     1.98 \\
            \hline
            ~ &$20^2 $&     2.25E-04 &  -  &     8.13E-04 & -  &     4.49E-04 &     - &     5.77E-04 &     - \\
            ~ &$40^2$ &     2.80E-05 &     3.01 &     1.06E-04 &     2.94  &     5.47E-05 &     3.04 &     7.18E-05 &     3.01 \\
            2 &$80^2$ &     3.49E-06 &     3.00 &     1.35E-05 &     2.98 &     6.85E-06 &     3.00 &     9.19E-06 &     2.97 \\
            ~ &$160^2$ &     4.36E-07 &     3.00 &     1.70E-06 &     2.99  &     8.24E-07 &     3.05 &     1.11E-06 &     3.05 \\
            ~ &$320^2$ &     5.44E-08 &     3.00 &     2.13E-07 &     3.00  &     7.03E-08 &     2.98 &     1.49E-07 &     2.95 \\
             \hline
            \end{tabular}
            %\vspace{\baselineskip}
        \end{table*}

    \begin{example}(Rigid body rotation)
    \end{example}
    Consider
        \begin{equation}\label{2D rigid}
             u_t-(y u)_x+(x u)_y=0, (x,y) \in[-2 \pi, 2 \pi]^2,
        \end{equation}
     with periodic boundary conditions.  Two
       initial condition are  tested, i.e.,  the circular symmetry one 
       $u(x, y, 0)=\exp \left(-x^2-y^2\right)$ and the asymmetric
       one $u(x, y, 0)=\exp \left(-x^2-10 y^2\right)$.      % The uniform mesh with $N_x=N_y=N$ SVs are used. 
   % We  adopt the ${\mathcal Q}^k$ SLSV methods with $k=1,2$ to solve this problem. 
    Table \ref{Tab: rigid1} and    Table \ref{Tab: rigid2} list 
  the $L^2$ and $L^{\infty}$ errors and the associated orders of accuracy
    calculated from   the ${\mathcal Q}^k$ SLSV methods with $k=1,2$ for the  symmetry one  and  asymmetric one, 
   respectively. 
     In our numerical experiment, we take 
    $\Delta t=0.5 h$ and $\Delta t=2.5 h$ and  $T=2 \pi$.  As expected, the $(k+1)$th order of accuracy for the ${\mathcal Q}^k$ splitting SLSV  scheme is observed in both the $L^2$ norm and $L^{\infty}$ norm.
        \begin{table*}[!ht]
            \centering
            \footnotesize
            \caption{\label{Tab: rigid1}Rigid body rotation: errors and convergence orders of $L^2$ and $L^{\infty}$  with  the 
            initial condition $u(x,y,0)=e^{-(x^2+y^2)}$ at $T=2\pi$.} 
            \vspace{0.4em} \centering%  xiaowu
            \begin{tabular}{c|c|cccc|cccc}%{\textwidth}
            %\toprule[0.001pt]
            \hline
            \multirow{2}{*}{$k$}   &\multirow{2}{*}{$N^2$}   & \multicolumn{4}{c|}{ $CFL=0.5$}  & \multicolumn{4}{c}{$CFL=2.5$} \\
            \cline{3-10}
            & &$L^2$ &Order&$L^{\infty}$& Order&$L^2$ &Order&$L^{\infty}$& Order\\
            %\midrule[0.001pt]
            \hline
            ~ &$20^2$ &     1.17E-02 &- &     1.62E-01 &- &     1.02E-02 &- &     1.37E-01 &-\\
            ~ &$40^2$ &     2.49E-03 &     2.23 &     5.98E-02 &     1.44 &     2.22E-03 &     2.20 &     5.56E-02 &     1.30 \\
            1 &$80^2$ &     5.28E-04 &     2.24 &     1.56E-02 &     1.94 &     4.77E-04 &     2.22 &     1.51E-02 &     1.88 \\
            ~ &$160^2 $&     1.22E-04 &     2.12 &     3.93E-03 &     1.99 &     1.12E-04 &     2.09 &     3.86E-03 &     1.97 \\
            ~ &$320^2$ &     2.96E-05 &     2.04 &     9.83E-04 &     2.00 &     2.76E-05 &     2.02 &     9.75E-04 &     1.99 \\
            \hline
            ~ &$20^2 $&     1.46E-03 &- &     2.28E-02 &- & 1.44E-03 &- &     1.80E-02 &-\\
            ~ &$40^2$ &     1.48E-04 &     3.30 &     3.53E-03 &     2.69 &   1.70E-04 &     3.08 &     2.93E-03 &     2.62\\
            2 &$80^2$ &     1.63E-05 &     3.18 &     3.63E-04 &     3.28 &   1.95E-05 &     3.12 &     4.06E-04 &     2.85\\
            ~ &$160^2$ &     1.99E-06 &     3.03 &     4.53E-05 &     3.00 &  2.41E-06 &     3.02 &     5.16E-05 &     2.98\\
            ~ &$320^2$ &     2.48E-07 &     3.00 &     5.60E-06 &     3.02 &    3.02E-07 &     3.00 &     6.51E-06 &     2.99\\
             \hline
            \end{tabular}
            %\vspace{\baselineskip}
        \end{table*}

        \begin{table*}[!ht]
            \centering
            \footnotesize
            \caption{\label{Tab: rigid2}Rigid body rotation: errors and convergence orders of $L^2$ and $L^{\infty}$  with the initial condition $u(x,y,0)=e^{-(x^2+10y^2)}$ at $T=2\pi$.} 
            \vspace{0.4em} \centering%  xiaowu
            \begin{tabular}{c|c|cccc|cccc}%{\textwidth}
            %\toprule[0.001pt]
            \hline
            \multirow{2}{*}{$k$}   &\multirow{2}{*}{$N^2$}   & \multicolumn{4}{c|}{ $CFL=0.5$}  & \multicolumn{4}{c}{$CFL=2.5$} \\
            \cline{3-10}
            & &$L^2$ &Order&$L^{\infty}$& Order&$L^2$ &Order&$L^{\infty}$& Order\\
            %\midrule[0.001pt]
            \hline
            ~ &$20^2$ &    2.59E-02 &- &     4.28E-01 &- &    2.44E-02 &- &     4.07E-01 &- \\
            ~ &$40^2$ &    1.82E-02 &     0.51 &     2.94E-01 &     0.54 &    1.66E-02 &     0.56 &     2.76E-01 &     0.56 \\
            1 &$80^2$ &    8.98E-03 &     1.02 &     1.40E-01 &     1.07 &     7.69E-03 &     1.11 &     1.25E-01 &     1.15 \\
            ~ &$160^2 $&    3.24E-03 &     1.47 &     5.41E-02 &     1.37 &     2.59E-03 &     1.57 &     4.55E-02 &     1.46 \\
            ~ &$320^2$&    9.12E-04 &     1.83 &     1.54E-02 &     1.81 &     7.06E-04 &     1.88 &     1.25E-02 &     1.86 \\
            \hline
            ~ &$20^2 $&     1.44E-02 &- &     2.71E-01 &- &     1.44E-02 &- &     2.71E-01 &- \\
            ~ &$40^2$ &     3.94E-03 &     1.87 &     6.51E-02 &     2.06 &     3.94E-03 &     1.87 &     6.51E-02 &     2.06 \\
            2 &$80^2$ &     5.01E-04 &     2.97 &     9.48E-03 &     2.78 &     5.01E-04 &     2.97 &     9.48E-03 &     2.78 \\
            ~ &$160^2$ &     4.38E-05 &     3.52 &     1.07E-03 &     3.15 &     4.38E-05 &     3.52 &     1.07E-03 &     3.15 \\
            ~ &$320^2$ &     4.26E-06 &     3.36 &     1.33E-04 &     3.00 &     4.26E-06 &     3.36 &     1.33E-04 &     3.00 \\
             \hline
            \end{tabular}
            %\vspace{\baselineskip}
        \end{table*}

    To test the splitting error in time,  we plot in Table \ref{Tab:rigid2 temporal order} the 
    $L^2$ and $L^{\infty}$ errors calculated from  ${\mathcal Q}^2$ SLSV-split scheme  at $T=20\pi$ (10 periods of rotation)
    with the initial condition  chosen as $u(x,y,0)=e^{-(x^2+10y^2)}$ and the time step size $\tau=\text{ CFL} h$, where CFL varies while 
    the mesh size $h$ is fixed (a fixed mesh with $160 \times 160$ cells). 
     We observe a  second-order splitting error and 
    the time error  increases with respect to the CFL number.

        \begin{table*}[!ht]
            \centering
            \footnotesize
            \caption{\label{Tab:rigid2 temporal order}Rigid body rotation: temporal errors and convergence orders for ${\mathcal Q}^2$ SLSV method with  the initial condition $u(x,y,0)=e^{-(x^2+10y^2)}$ at $T=20\pi$. } 
            \vspace{0.4em} \centering%  xiaowu
            \begin{tabular}{c|cccc}%{\textwidth}
            %\toprule[0.001pt]
            \hline
            CFL    &$L^2$ &Order&$L^{\infty}$& Order\\
            %\midrule[0.001pt]
            \hline
            5 &     7.99E-04 &-      &     1.60E-02 &-     \\
            10 &     1.02E-03 &     0.37 &     1.69E-02 &     0.08 \\
            15 &     1.89E-03 &     1.51 &     2.69E-02 &     1.15 \\
            20 &     3.27E-03 &     1.91 &     4.39E-02 &     1.71 \\
            25 &     5.10E-03 &     1.99 &     6.76E-02 &     1.93 \\
            30 &     7.33E-03 &     1.99 &     9.65E-02 &     1.96 \\
            35 &     9.92E-03 &     1.97 &     1.31E-01 &     1.97 \\
            \hline
            \end{tabular}
            %\vspace{\baselineskip}
        \end{table*}

\begin{figure}[!h]
            \centering
            \includegraphics[width=0.5\textwidth]{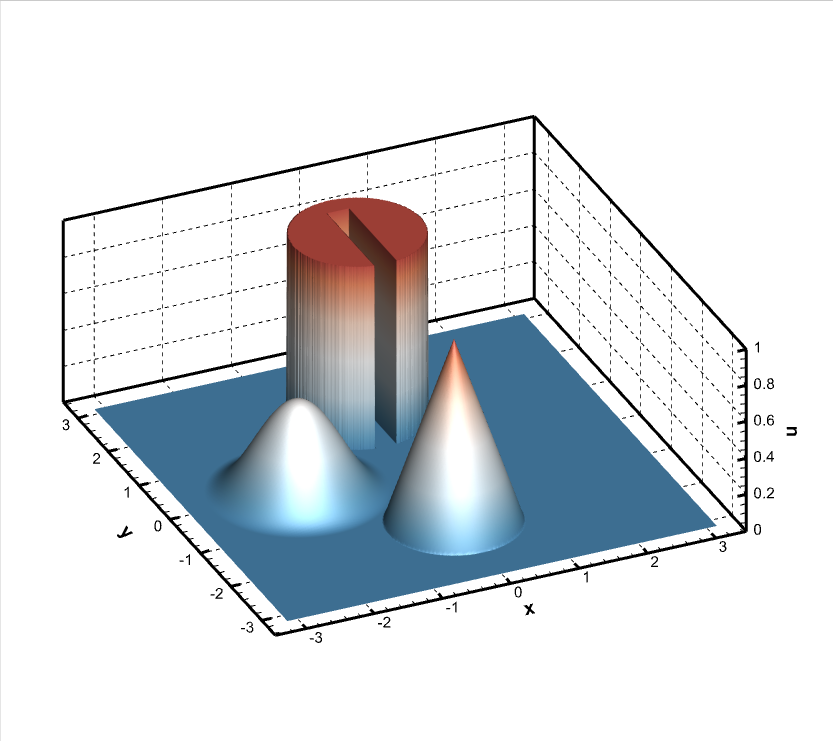}
            \caption{\label{Fig:con initial}  Plots of the initial profle. The mesh of $400 \times 400$ is used.}
        \end{figure}

 To test the shape preservation of the SLSV method for the rigid body rotation \eqref{2D rigid}, we consider 
  the computational domain  $[-\pi, \pi]$ with the same initial condition used in  \cite{coneinitial}, which consists of a slotted disk, a cone as well as a smooth hump (see  Figure \ref{Fig:con initial}). 
  We note that  
  the oscillations will be introduced for high-order schemes for solving this problem as the solution is discontinuous. 
        
\begin{figure}[!h]
            \centering
            \includegraphics[width=0.37\textwidth]{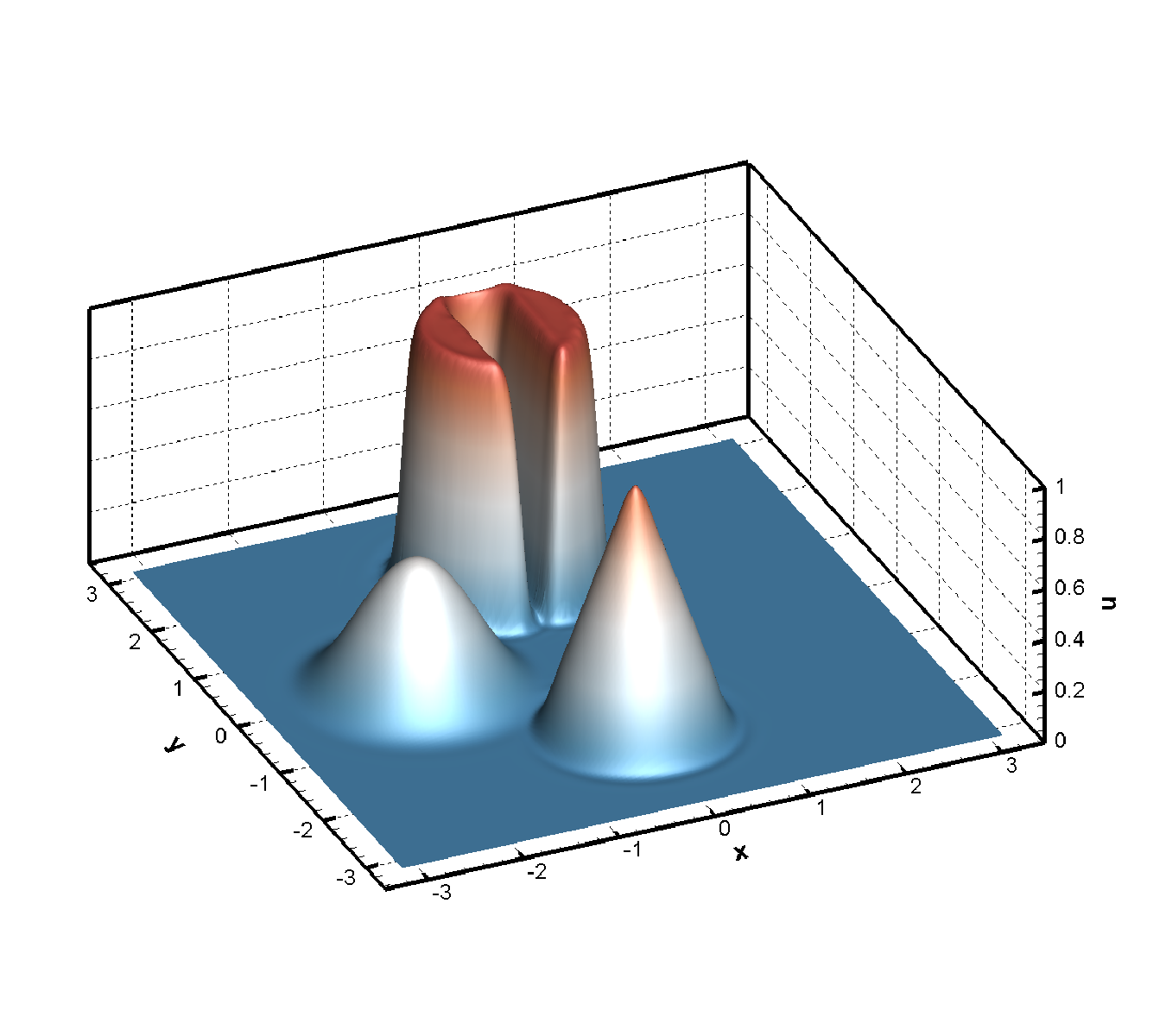}
            \includegraphics[width=0.37\textwidth]{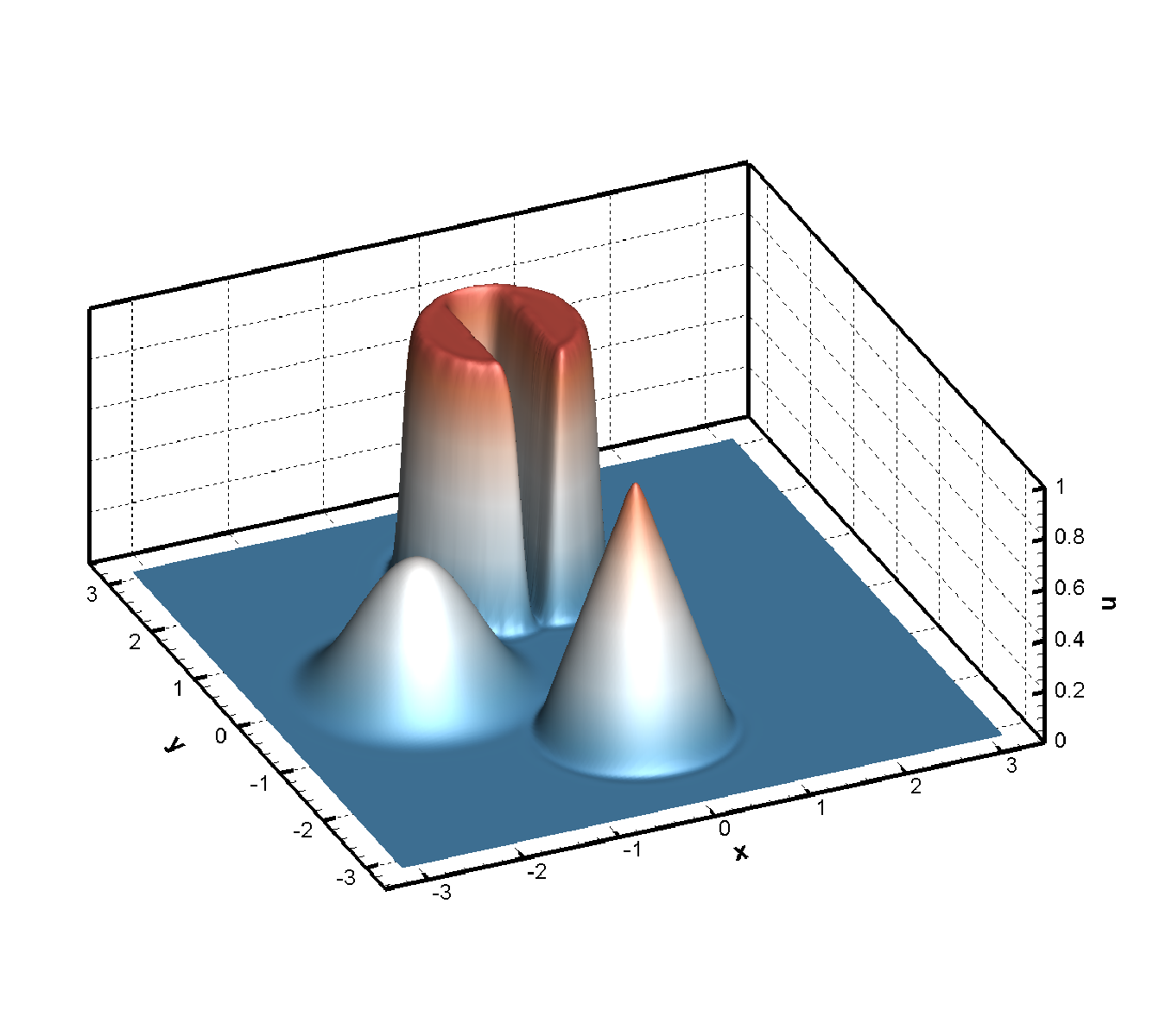}
            \caption{\label{Fig:rigid cone} Rigid body rotation: Plots of the numerical solutions of SLSV schemes with initial data Figure \ref{Fig:con initial}  at a mesh of $160 \times 160$ and $T = 12\pi$.  ${\mathcal Q}^2$ SLSV-split + WENO limiter is used with           
             % Upper left: mesh of $100 \times 100$ with CFL = 2.2. Upper right: mesh of $100 \times 100$ with CFL = 10.2. Bottom left: mesh of $160 \times 160$
             CFL = 2.2 (left) and  CFL = 10.2 (right).}
        \end{figure}

\begin{figure}[!h]
            \centering
            \includegraphics[width=0.37\textwidth]{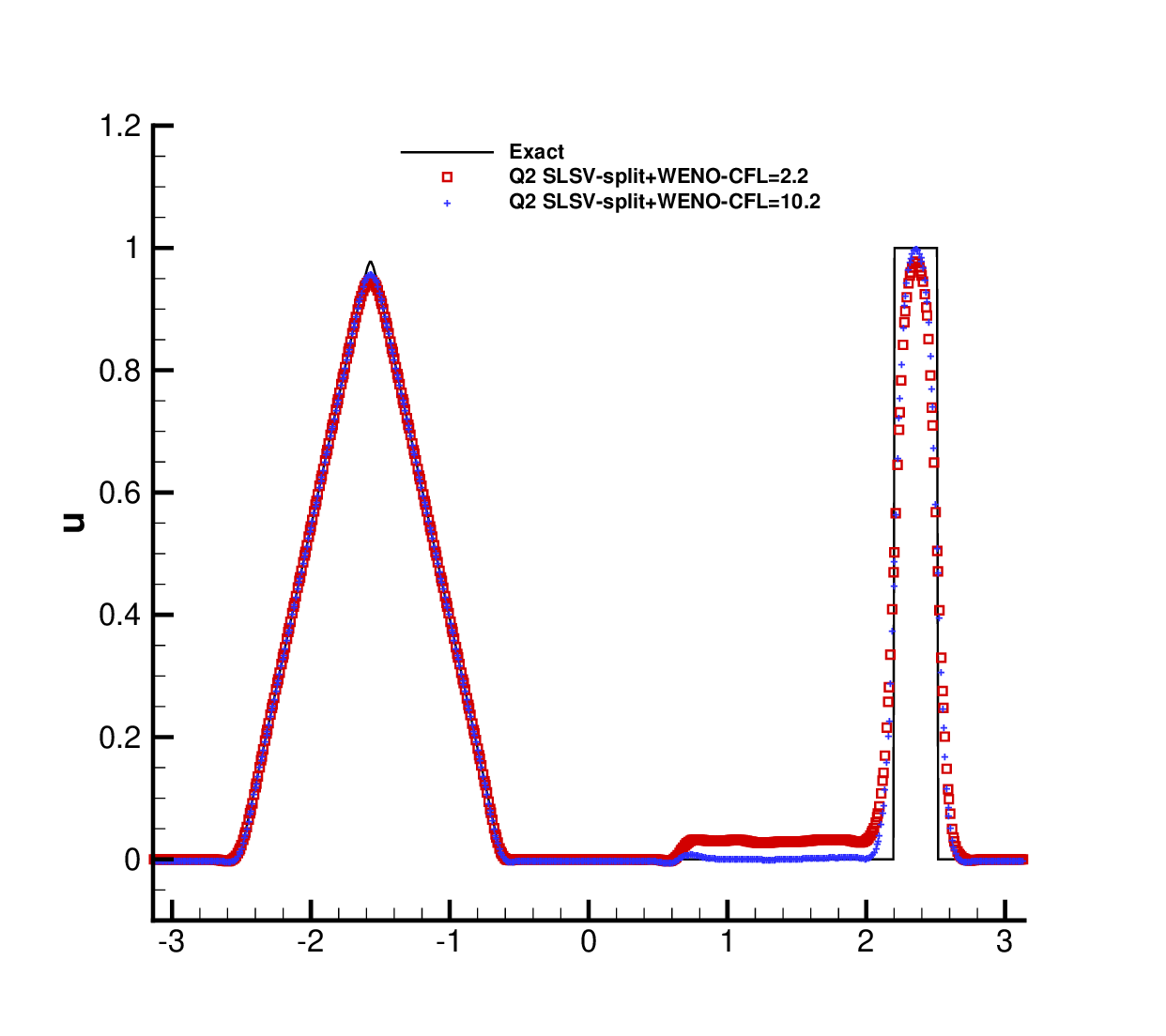}
            \includegraphics[width=0.37\textwidth]{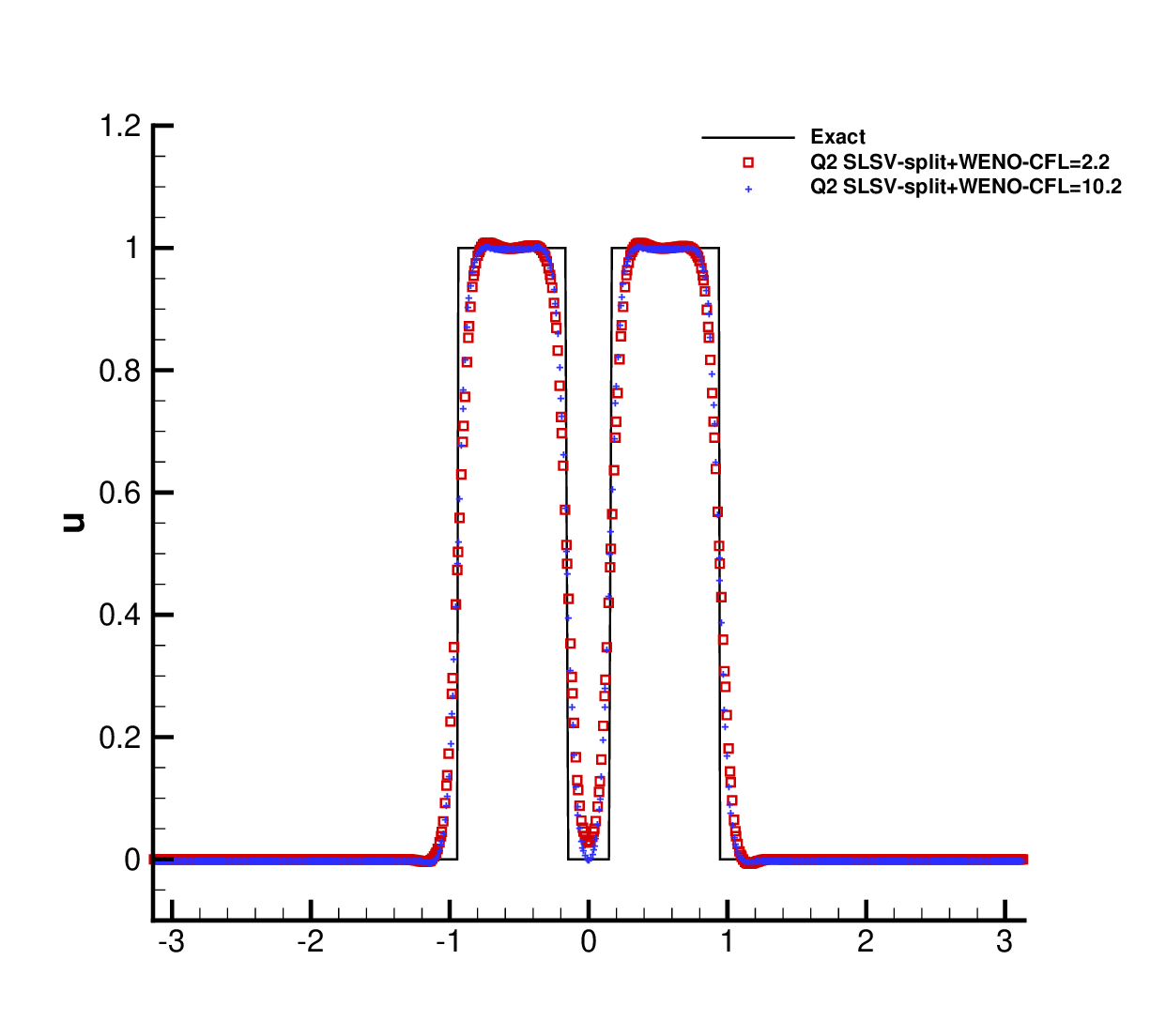}
            \caption{\label{Fig:rigid cone cut} Rigid body rotation: Plots of 1D cuts of the numerical solution for SLSV schemes with initial data Figure \ref{Fig:con initial} at a mesh of $160\times160$ and $T = 12\pi$.  Left: numerical solution at $x = 0 + \pi/100$. Right: numerical solution at $y = \pi/2 + \pi/100$. 
            %Upper: mesh of $100\times100$. Bottom: mesh of $160\times160$.
         }
        \end{figure}

  In our numerical experiment,  we adopt a simple WENO limiter in \cite{Wenolimiter} to suppress the oscillations. 
 Plotted in Figure \ref{Fig:rigid cone} is the numerical solution at time $T=12\pi$, i.e., 
 after six full revolutions. To demonstrate the shape preservation property, 
  we plot the one-dimensional cuts of the numerical and exact solutions in Figure \ref{Fig:rigid cone cut}. 
   As we may observe,   the WENO limiter effectively suppresses the unphysical oscillations of the numerical solution, 
  and the solution of the SLSV scheme with a larger time step ($\mathrm{CFL}=10.2$) is less dissipative than that of the SLSV scheme with a smaller time step ($\mathrm{CFL}=2.2$).

   \begin{example}(Swirling deformation flow)\label{Swirling deformation flow accurate}
    \end{example} 
    We consider the following Swirling deformation flow  with periodic boundary
        \begin{equation*}\label{2D swirling}
            u_t-\left(\cos ^2\left(\frac{x}{2}\right) \sin (y) g(t) u\right)_x+\left(\sin (x) \cos ^2\left(\frac{y}{2}\right) g(t) u\right)_y=0,(x, y) \in[-\pi, \pi]^2,
        \end{equation*}
    where $g(t)=\cos \left(\frac{\pi t}{T}\right) \pi, T=1.5$. The initial condition is set to be the following smooth cosine bells (with $C^5$ smoothness):
        \begin{equation*}\label{2D swirling cosine bell}
             u(x, y, 0)= \begin{cases}r_0^b \cos ^6\left(\frac{r^b}{2 r_0^b} \pi\right), & \text { if } r^b<r_0^b, \\ 0, & \text { otherwise, }\end{cases}
        \end{equation*}
    where $r_0^b=0.3 \pi$, and $r^b=\sqrt{\left(x-x_0^b\right)^2+\left(y-y_0^b\right)^2}$ denotes the distance between $(x, y)$ and the center of the cosine bell $\left(x_0^b, y_0^b\right)=(0.3 \pi, 0)$. The swirling deformation flow introduced in \cite{coneinitial} is more challenging than the rigid body rotation due to the space- and time-dependent flow field. In particular, along the direction of the flow, the initial function becomes largely deformed at $t=T/2$, then goes back to its initial shape at $t=T$ as the flow reverses. 
    If this problem is solved up to $T$,  the procedure is referred as a one full evolution.  

    We present in  Table \ref{Tab:swirling}  the convergence  for ${\mathcal Q}^1$ and ${\mathcal Q}^2$ SLSV scheme in terms of the $L^2$ and $L^{\infty}$ errors. A bit order increase is observed for the ${\mathcal Q}^1$ SLSV scheme on both $L^2$ and $L^{\infty}$ errors. In addition, the third order convergence for the ${\mathcal Q}^2$ SLSV scheme measured by both errors are observed. 

        \begin{table*}[!ht]
            \centering
            \footnotesize
            \caption{\label{Tab:swirling} Swirling deformation flow: errors and convergence orders of $L^2$ and $L^{\infty}$ at $T=1.5$.} 
            \vspace{0.4em} \centering%  xiaowu
            \begin{tabular}{c|c|cccc|cccc}%{\textwidth}
            %\toprule[0.001pt]
            \hline
            \multirow{2}{*}{$k$}   &\multirow{2}{*}{$N^2$}   & \multicolumn{4}{c|}{ $CFL=0.5$}  & \multicolumn{4}{c}{$CFL=2.5$} \\
            \cline{3-10}
            & &$L^2$ &Order&$L^{\infty}$& Order&$L^2$ &Order&$L^{\infty}$& Order\\
            %\midrule[0.001pt]
            \hline
            ~ &$20^2$ &     2.09E-02 &- &     3.82E-01 &-    &     9.88E-03 &     -&     1.87E-01 &     -\\
            ~ &$40^2$ &     7.02E-03 &     1.57 &     1.45E-01 &     1.40&     2.15E-03 &     2.20 &     5.14E-02 &     1.79 \\
            1 &$80^2$ &     1.47E-03 &     2.26 &     3.31E-02 &     2.13&4.03E-04 &     2.41 &     1.28E-02 &     2.08 \\
            ~ &$160^2 $&     2.50E-04 &     2.55 &     6.02E-03 &     2.46 &     8.47E-05 &     2.25 &     2.75E-03 &    2.22  \\
            ~ &$320^2$&     4.68E-05 &     2.42 &     1.10E-03 &     2.45   &     3.28E-05 &     2.17 &     8.11E-04 &     2.28 \\
            \hline
            ~ &$20^2$ &     4.63E-03 &- &     8.27E-02 &- &     1.13E-03 &     - &     1.97E-02 &     -\\
            ~ &$40^2$ &     5.83E-04 &     2.99 &     1.05E-02 &     2.97 &     1.26E-04 &     3.17 &     2.85E-03 &     2.79 \\
            2 &$80^2$ &     6.33E-05 &     3.20 &     1.23E-03 &     3.09&     1.49E-05 &     3.08 &     4.15E-04 &     2.78\\
            ~ &$160^2$ &     7.51E-06 &     3.08 &     1.54E-04 &     3.01 &     1.86E-06 &     3.01 &     4.91E-05 &     3.08 \\
            ~ &$320^2$ &     9.26E-07 &     3.02 &     1.94E-05 &     2.98 &     1.19E-06 &     2.99 &     2.37E-05 &     2.95\\
             \hline
            \end{tabular}
            %\vspace{\baselineskip}
        \end{table*}

\begin{figure}[!h]
            \centering
            \includegraphics[width=0.37\textwidth]{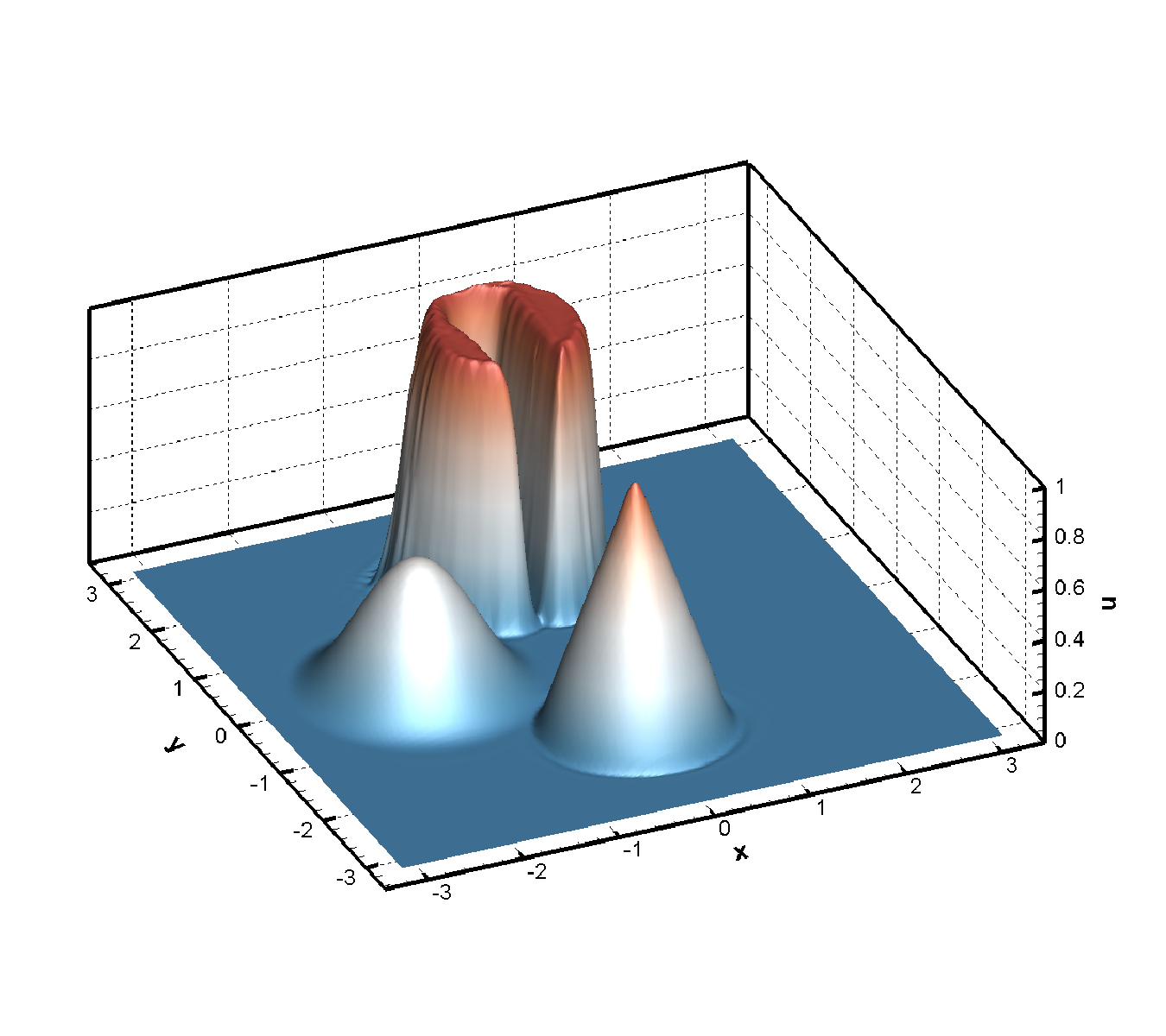}
            \includegraphics[width=0.37\textwidth]{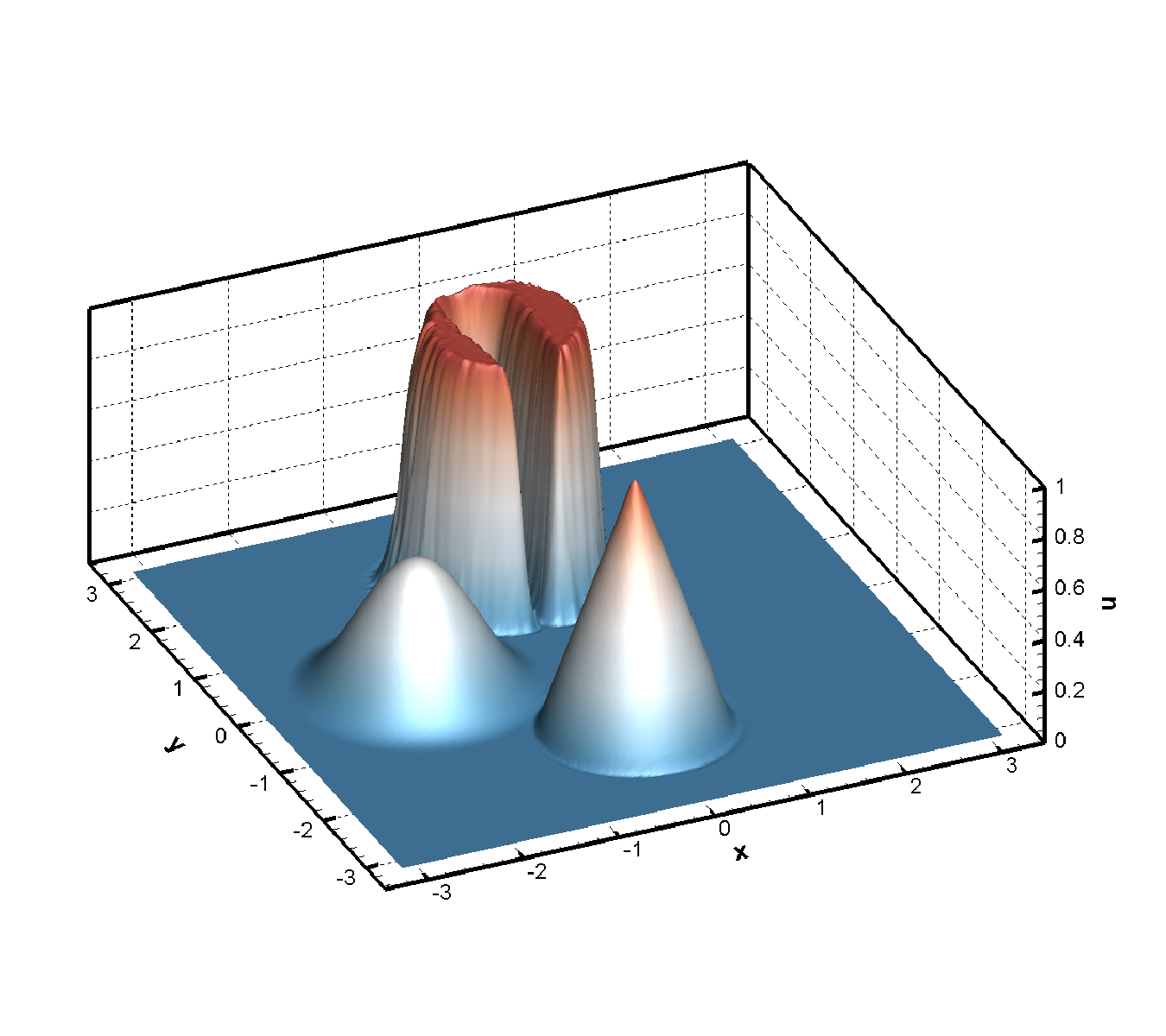}
            \caption{\label{Fig:Swirling-field-1.5} Swirling deformation flow: Plots of the numerical solutions of SLSV schemes  with initial data Figure \ref{Fig:con initial} up to $t = 1.5$. The ${\mathcal Q}^2$ SLSV-split+WENO is used with a mesh of $160\times160$. 
            Left:  CFL = 2.2. Right:  CFL = 10.2. % Upper: mesh of $100\times100$. Bottom: mesh of $160\times160$.
            }
        \end{figure}

\begin{figure}[!h]
            \centering
            \includegraphics[width=0.37\textwidth]{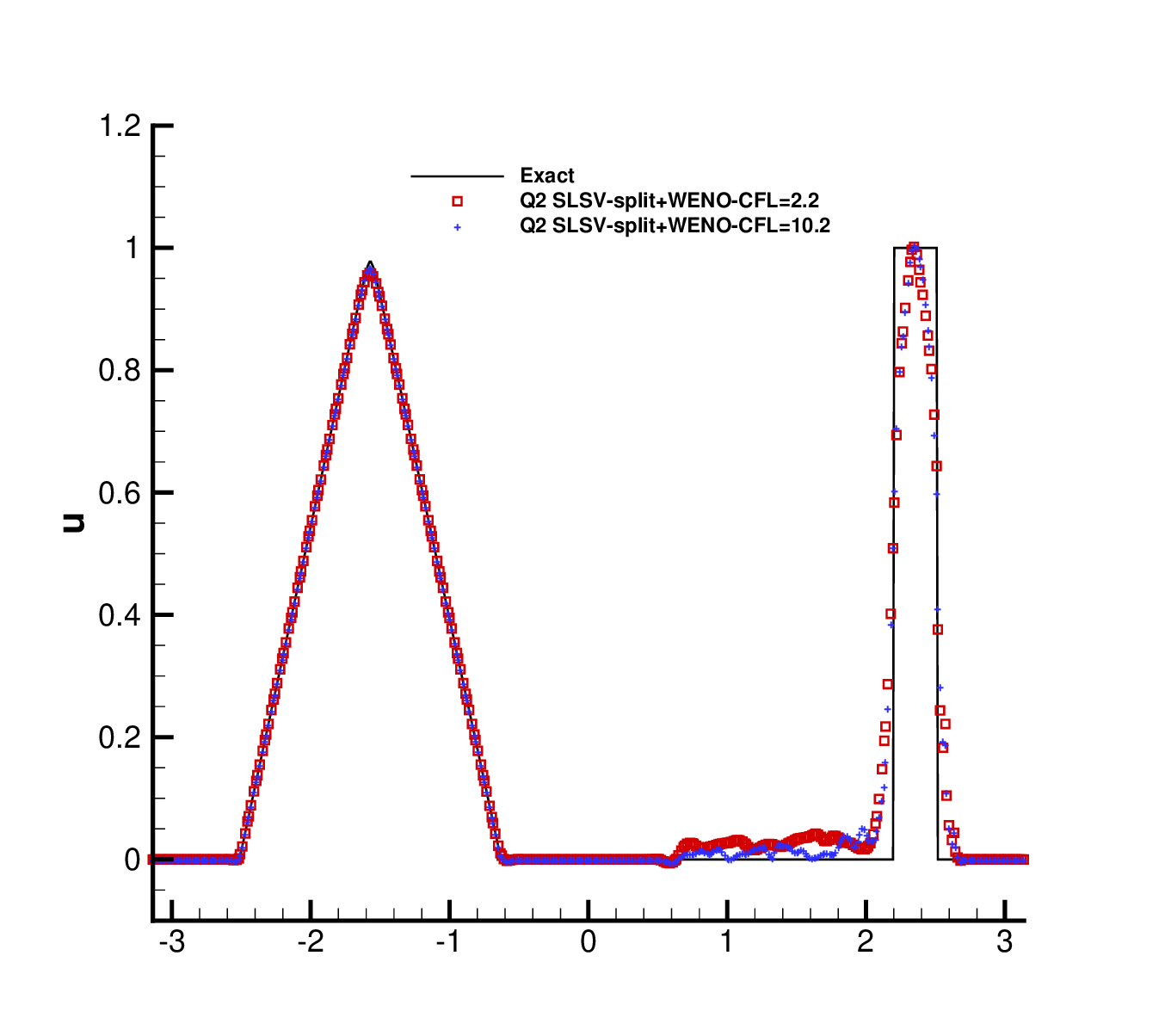}
            \includegraphics[width=0.37\textwidth]{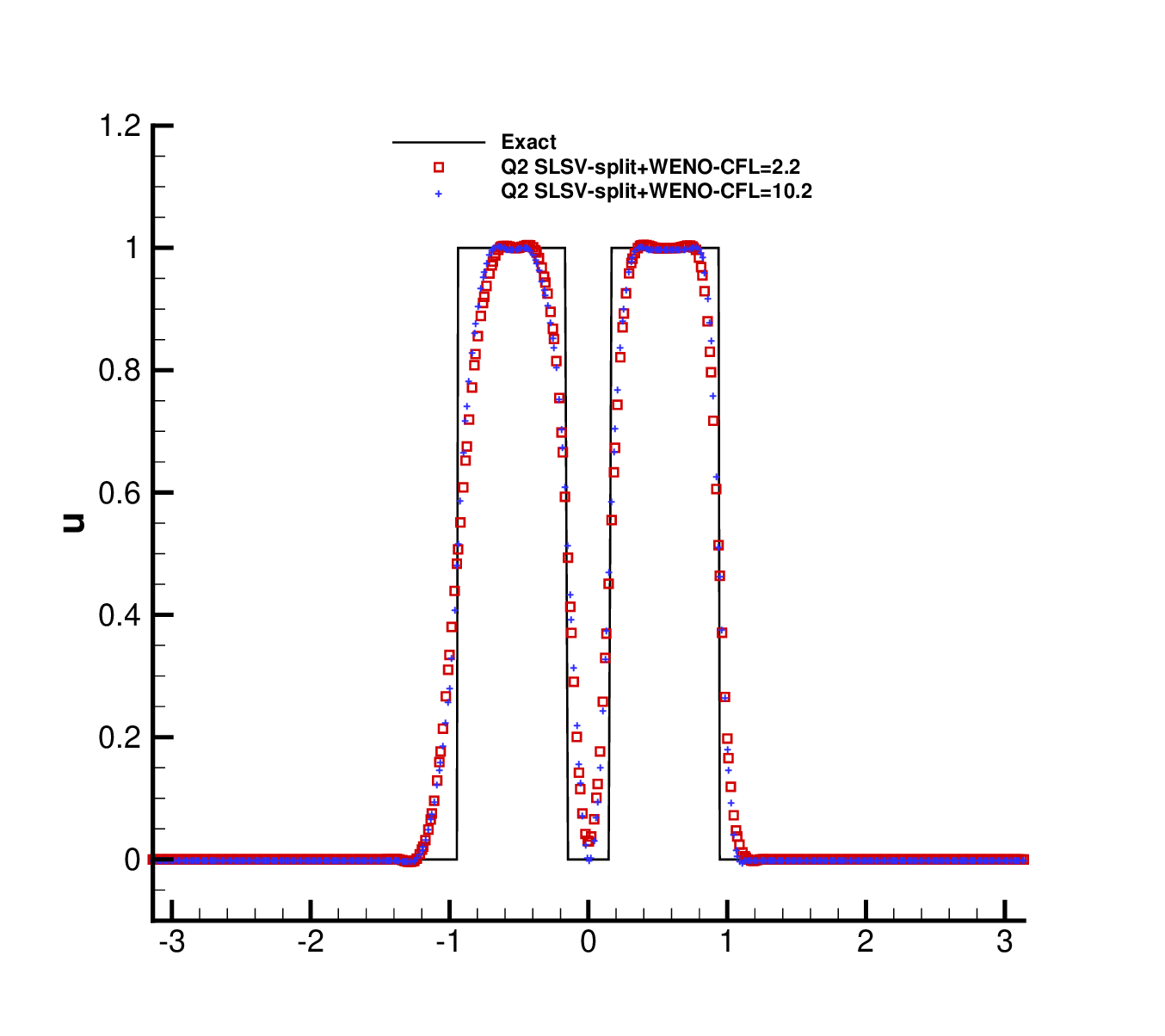}
            \caption{\label{Fig:Swirling-field-1.5-cut} Swirling deformation flow: Plots of 1D cuts of the numerical solution for SLSV schemes with initial data Figure \ref{Fig:con initial} at a mesh of $160\times160$ and $t = 1.5$.  Left: numerical solution at $x = 0 + \pi/100$. Right: numerical solution at $y = \pi/2 + \pi/100$. }
        \end{figure}

\begin{figure}[!h]
            \centering
            \includegraphics[width=0.37\textwidth]{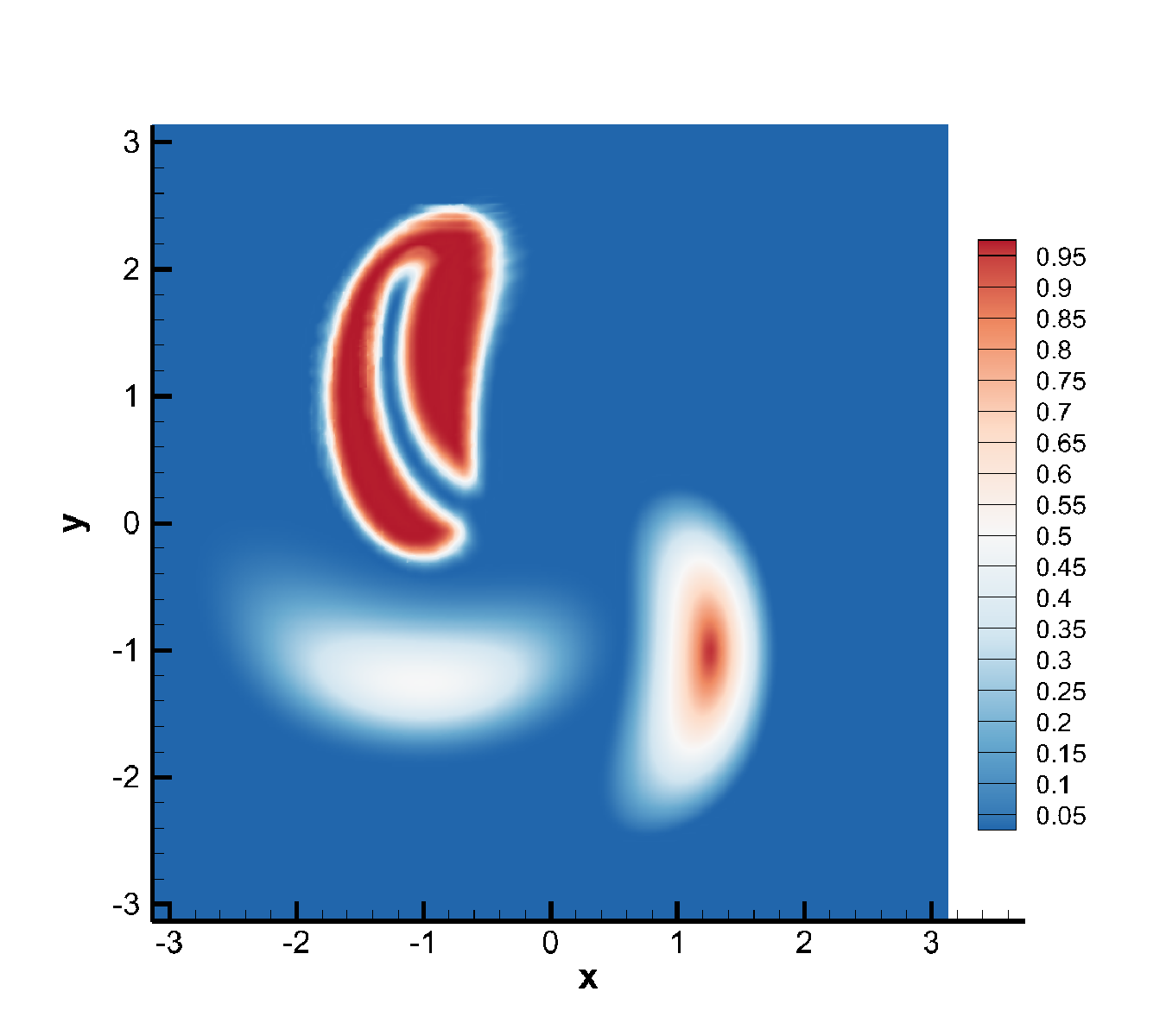}
            \includegraphics[width=0.37\textwidth]{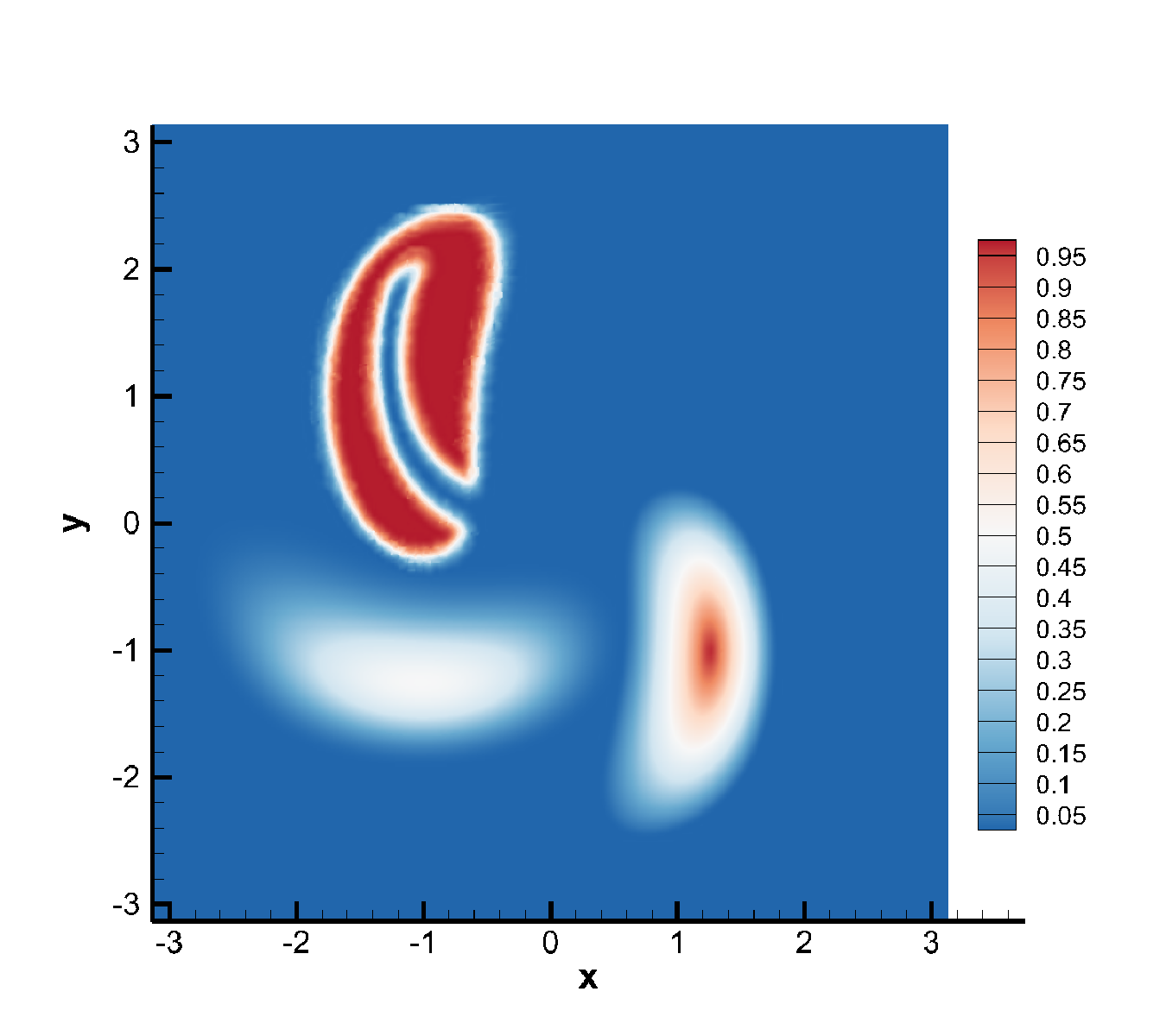}
            \caption{\label{Fig:Swirling-field-0.75} Swirling deformation flow:  Plots of the numerical solutions of SLSV schemes  with initial data Figure \ref{Fig:con initial} up to $t = 0.75$. The ${\mathcal Q}^2$ SLSV-split+WENO is used using a mesh of $160\times160$. Left:  CFL = 2.2. Right:  CFL = 10.2. % Upper: mesh of $100\times100$. Bottom: mesh of $160\times160$.
           }
        \end{figure}

Plotted in Figure \ref{Fig:Swirling-field-1.5} are the 
numerical solutions calculated from the SLSV method  with the initial condition in Figure \ref{Fig:con initial}.  
To better compare performance of different configurations, we show 1D cuts of the numerical solutions along with the exact solution in Figure \ref{Fig:Swirling-field-1.5-cut}. We observe 
that some mild spurious oscillations appear without using limiters, which can be removed by further coupling the WENO limiter. In addition, positivity of the numerical solution is guaranteed when the BP filter is applied. 
Figure \ref{Fig:Swirling-field-0.75} shows the contour plots of the numerical solution for SLSV methods at the final integration time 0.75, when the solution is quite deformed. 

\subsection{Vlasov‑Poisson System}\label{Sec:VPS}

In this section, we demonstrate the following different aspects via extensive numerical tests of the proposed algorithm on a set of benchmark VP examples.

\begin{figure}[!h]
            \centering
            \includegraphics[width=0.37\textwidth]{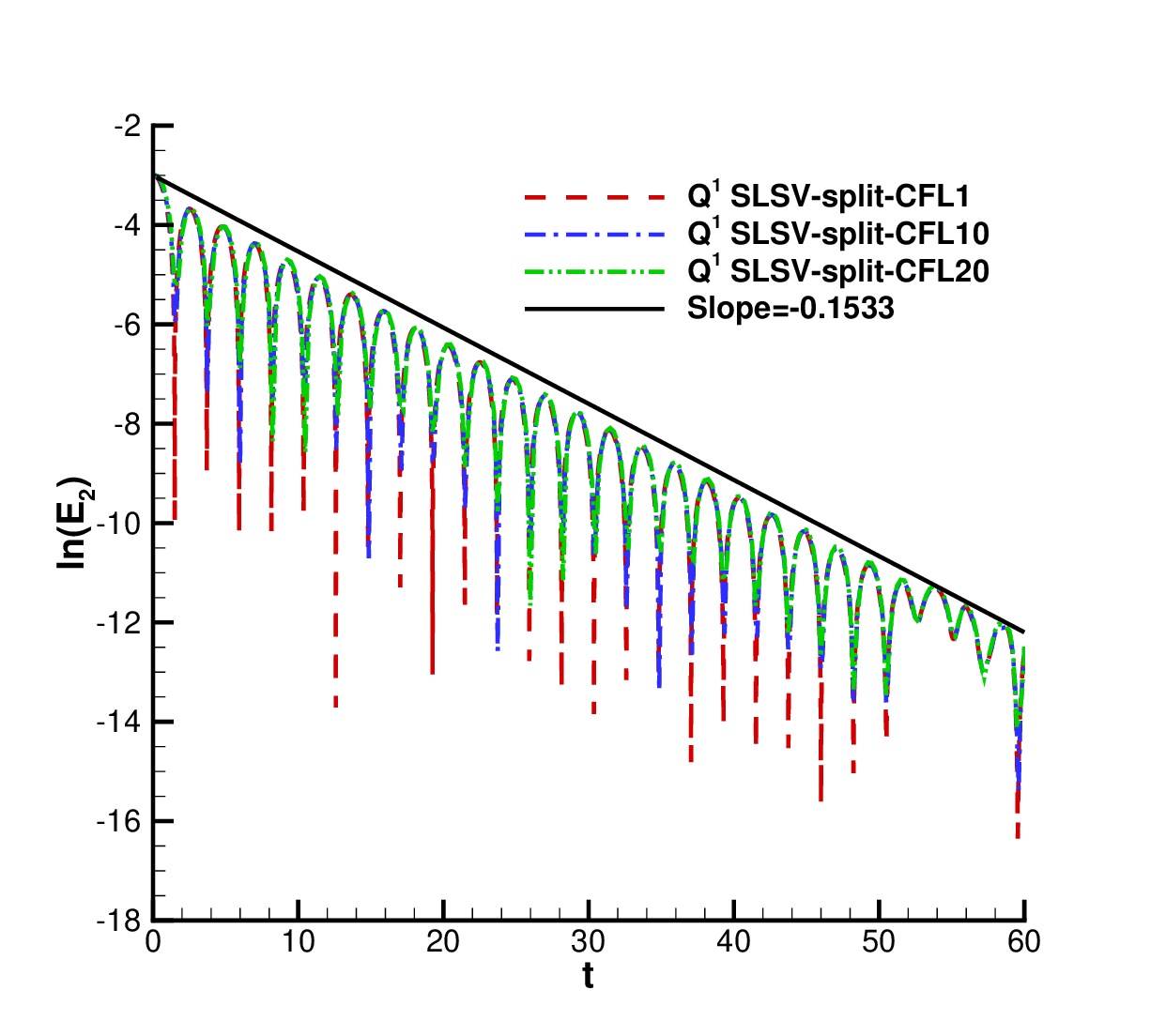}
            \includegraphics[width=0.37\textwidth]{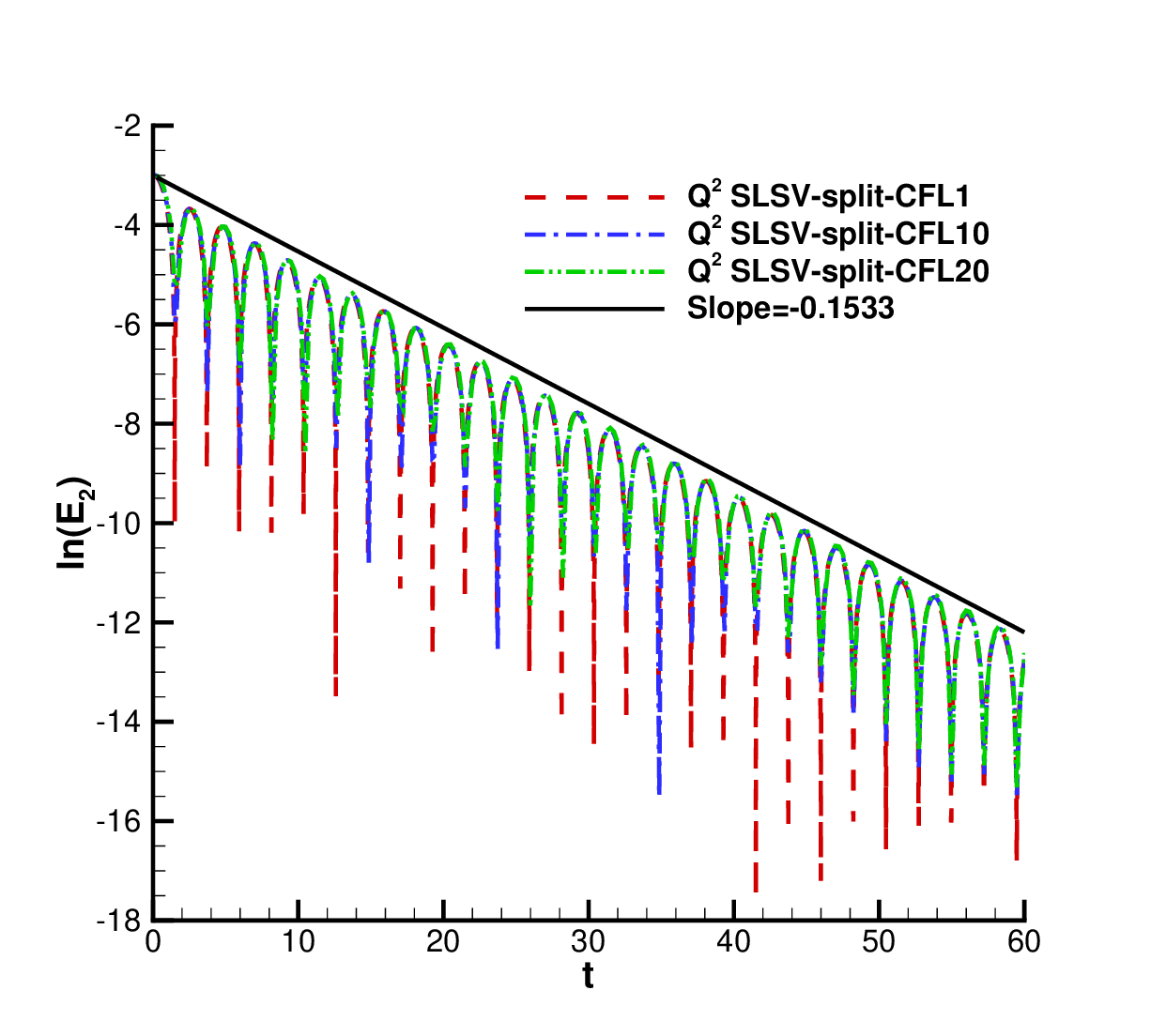}
            \caption{\label{Fig:weak_E2} Weak Landau damping: The SLSV schemes are equipped with the PP limiter. Time evolution of the electric field in $L^2$, using a mesh of $160 \times 160$ elements.}
        \end{figure}

\begin{itemize}
    \item[1.] Using a high spatial order. For weak Landau damping, we benchmark the numerical damping rate of the electrostatic field $E$ against the theoretical value from the linear theory. In particular, we test  ${\mathcal Q}^1$ and ${\mathcal Q}^2$ SLSV using CFL numbers as large as 20. In Figure \ref{Fig:weak_E2}, we show
     the advantage of using the third order SLSV scheme (compared with a second order one) for its superior performance in capturing the correct damping rate with CFL equals to 20.
    \item[2.] Spatial order of convergence. We test the spatial order of convergence in Table \ref{Tab: strong landau damping} for strong Landau damping, and in Table \ref{Tab: Two stream1 accruate} for two stream instabilities, respectively. The proposed SLSV method has a spatially optimal convergence order.
    \item[3.] Preservation of mass and other physical norms. 
     Tracking relative deviations of these quantities numerically provides a good measurement of the quality of numerical schemes. Our proposed SLSV scheme is locally and globally mass conservative. We will show comparable (sometimes superior) performance in preserving the physical norms for the proposed SLSV scheme with large CFLs. 
\end{itemize}

   \begin{example}(Weak Landau damping.)
    \end{example} 
    Consider weak Landau damping for the VP system. The initial condition is set to be the following perturbed equilibrium
        \begin{equation}\label{weak landou damping initial}
            f(x, v, t=0)=\frac{1}{\sqrt{2 \pi}}(1+\alpha \cos (k x)) \exp \left(-\frac{v^2},{2}\right)
        \end{equation}
    with $\alpha=0.01$ and $k=0.5$. Our computational domain is $[0,4 \pi] \times\left[-v_{\max }, v_{\max }\right]$. We truncate the velocity domain at $v_{\max }=2 \pi$. We also refer to \cite{Weaklandau1,Weaklandau2,Weaklandau3,Weaklandau4,VPS1,VPS2,VPS3} for the 
 numerically study of the weak Landau damping problem. 
            
\begin{figure}[!h]
            \centering
            \includegraphics[width=0.37\textwidth]{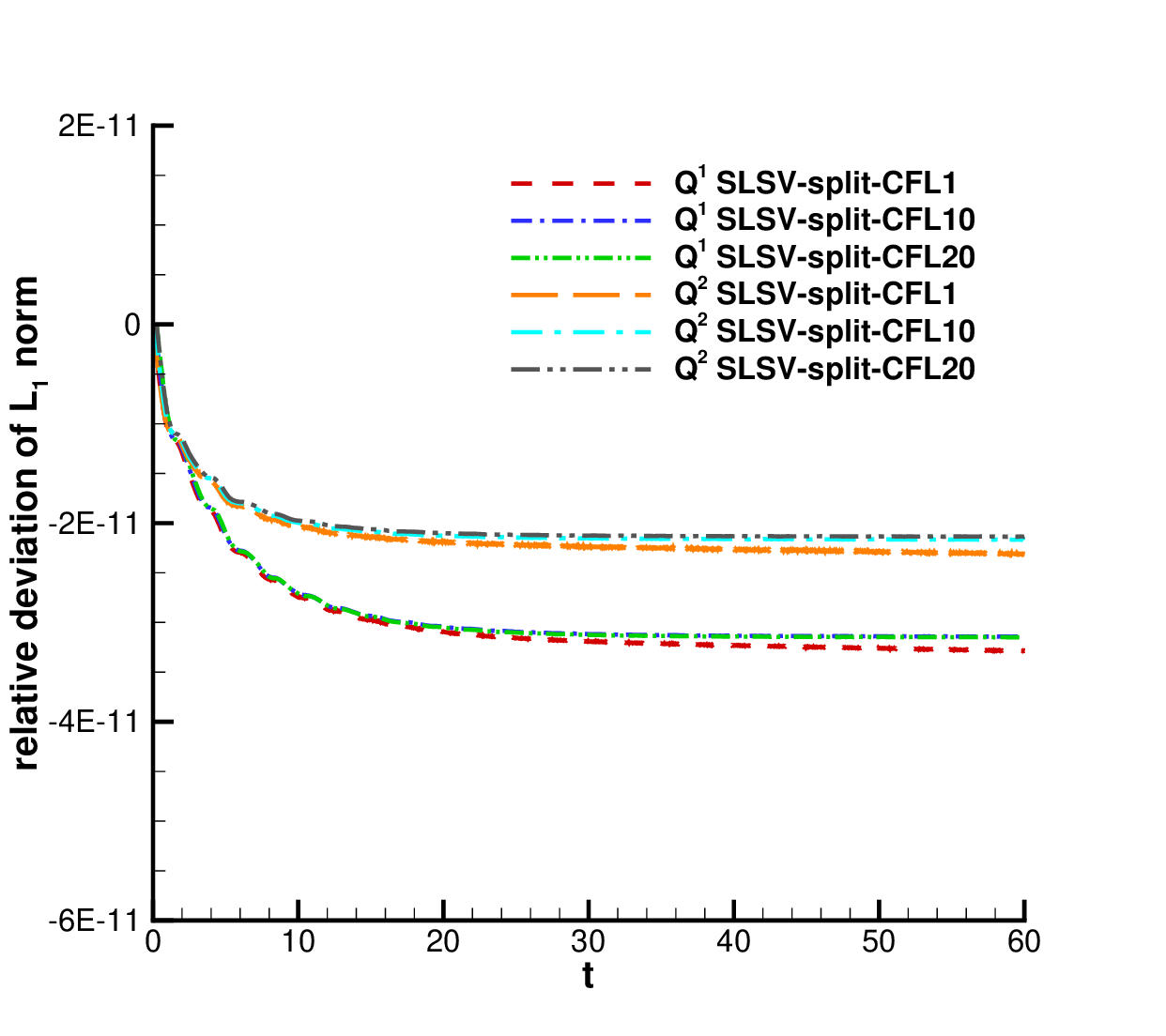}
            \includegraphics[width=0.37\textwidth]{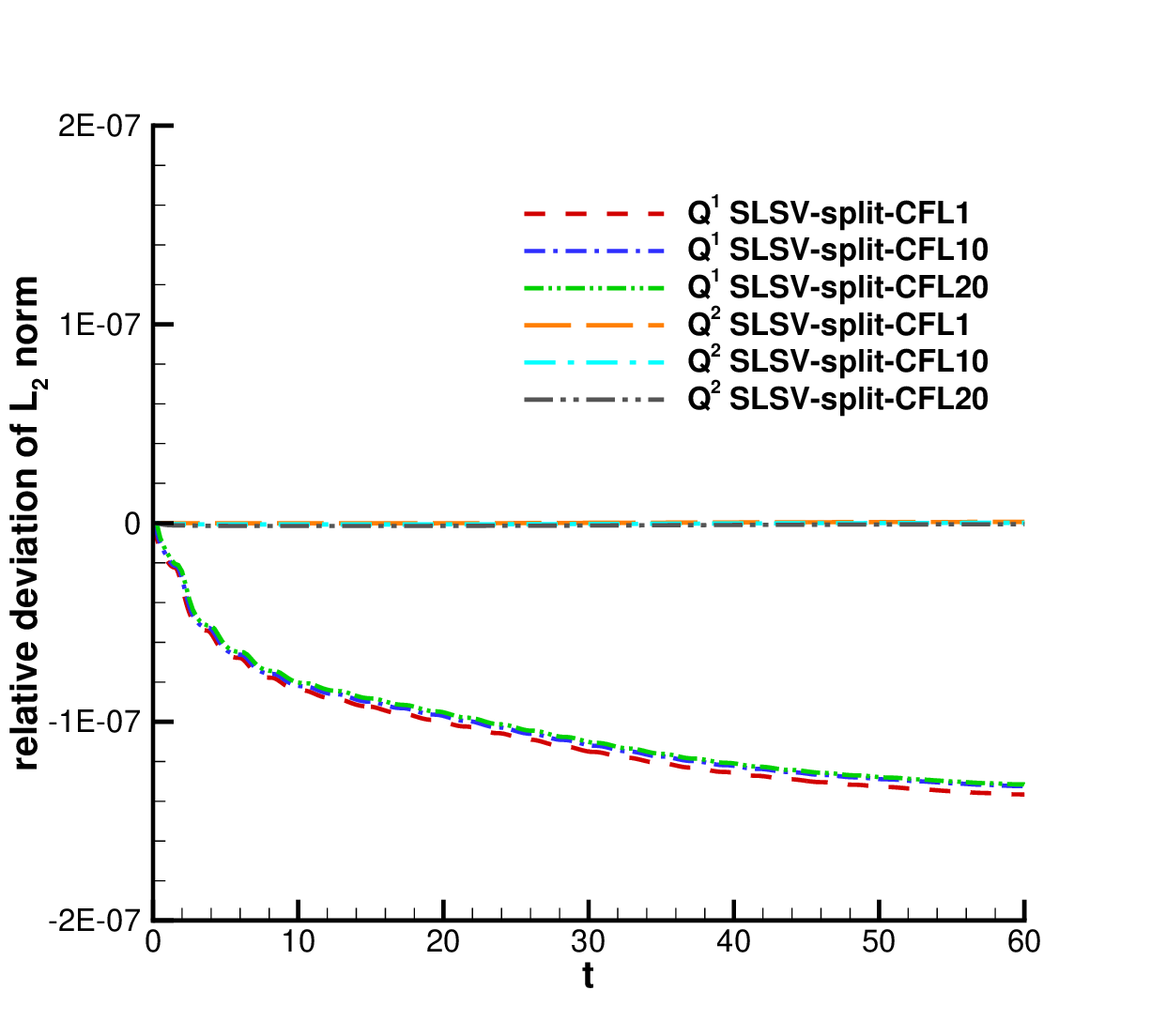}\\
            \includegraphics[width=0.37\textwidth]{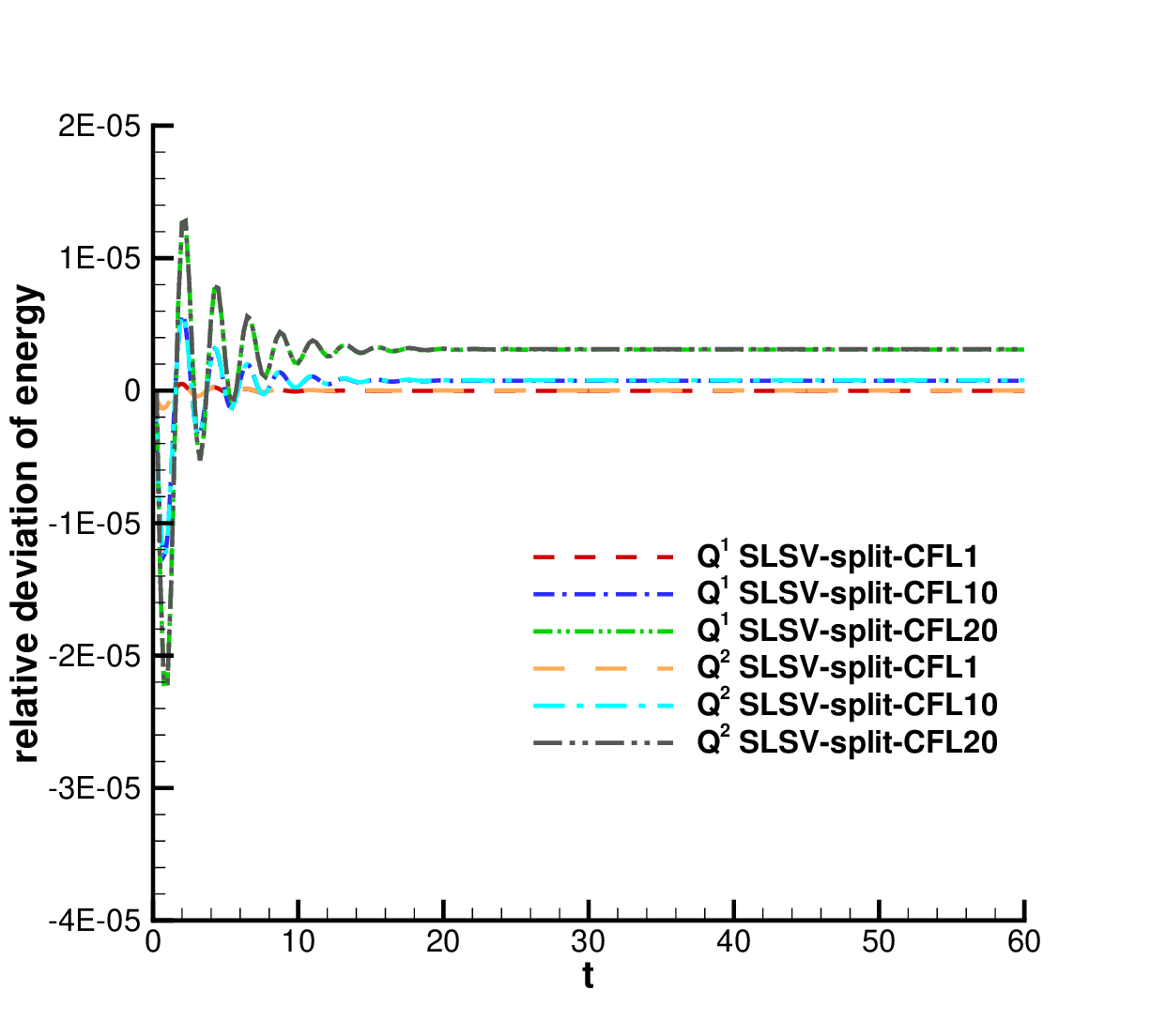}
            \includegraphics[width=0.37\textwidth]{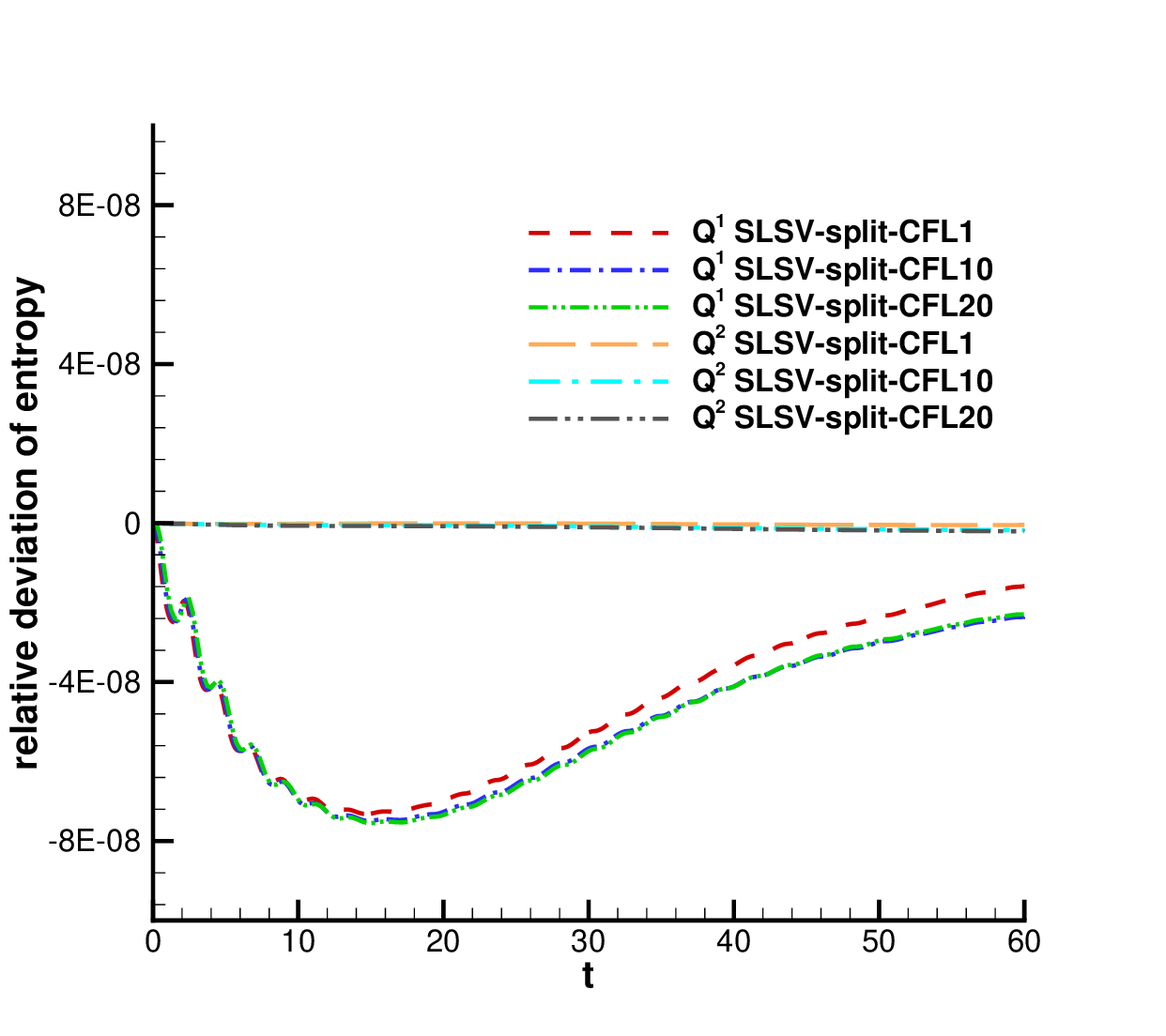}
            \caption{\label{Fig:weak_quantity}  Weak Landau damping. Time evolution of the relative deviations of $L^1$ (upper left) and $L^2$ (upper right) norms of the solution as well as the discrete kinetic energy (lower left) and entropy (lower right), using a mesh of $160 \times 160$ elements.}
        \end{figure}

Since $\alpha$ is   small, we expect to see results that closely agree with the linear theory, where the electric field decays exponentially. In Figure \ref{Fig:weak_E2} we present this decay provided by
    $$
    \ln \left(\|E(t, \cdot)\|_{L_2}\right):=\ln \sqrt{\int_{0}^{4 \pi}|E(t, x)|^2 d x},
    $$
    versus time for ${\mathcal Q}^1$ and ${\mathcal Q}^2$ splitting SLSV schemes using a mesh of $160 \times 160$ elements and different CFLs. Our computed decay rate from simulations with CFL$=1$ matches the linear decay rate, $\gamma=-0.1533$ (the solid line in the same plots). For CFL$=10$, both results match well with the theoretical value; for CFL$=20$, the third order characteristics tracing scheme exhibits superior performance, compared with the second order one, in capturing the correct damping rate in the long run. In Figure 
    \ref{Fig:weak_quantity}, we present the deviations of several quantities that are conserved by the continuous Vlasov-Poisson system from their initial values: $\|f\|_{L_1},\|f\|_{L_2}$, total energy, and entropy. In particular, we observe 
    the following numerical resluts:
    (1)The error for the $L^1$ norm (on the order of $10^{-9}$ ) is due to the truncation of the velocity domain, which can be further reduced by using a larger velocity domain in simulations.
    (2) In general, the ${\mathcal Q}^2$  SLSV method does a better job in conserving these physical norms than the ${\mathcal Q}^1$ SLSV method.

    \begin{example}(Strong Landau damping) 
    \end{example} 
    Consider strong Landau damping for the VP system. The initial condition is the same as weak one \eqref{weak landou damping initial}, but with a larger perturbation parameter $\alpha=0.5$. The computational domain is $[0,4 \pi] \times[-2 \pi, 2 \pi]$ with periodic boundary in $x$- direction. 
    
    In Table \ref{Tab: strong landau damping}, we test the spatial convergence of the proposed SLSV scheme with the third order characteristic tracing scheme. We set CFL$=0.1$ to minimize the error from time discretization. The well-known time reversibility of the VP system is used to test the order of convergence. In particular, one can integrate the VP system forward to some time $T$, and then reverse the velocity field of the solution and continue to integrate the system by the same amount of time $T$. Then, the solution should recover the initial condition with reverse velocity field, which can be used as a reference solution. We show the $L^2$ and $L^{\infty}$ errors and the corresponding orders of convergence for ${\mathcal Q}^k$ SLSV schemes with CFL$=0.1$ in Table \ref{Tab: strong landau damping}. Again the 
    $(k+1)$th order convergence is observed for ${\mathcal Q}^k$ SLSV scheme,  as expected. 
    \begin{table*}[!ht]
        \centering
        \footnotesize
        \caption{\label{Tab: strong landau damping} Strong Landau damping: errors and convergence orders of $L^2$ and $L^{\infty}$  at $T=0.5$ with $CFL=0.1$.} 
        \vspace{0.4em} \centering%  xiaowu
        \begin{tabular}{c|c|cccc}%{\textwidth}
        %\toprule[0.001pt]
        \hline
        \multirow{2}{*}{$k$}   &\multirow{2}{*}{$N^2$}   & \multicolumn{4}{c}{ $CFL=0.1$}   \\
        \cline{3-6}
        & &$L^2$ &Order&$L^{\infty}$& Order\\
        %\midrule[0.001pt]
        \hline
        ~ &$16^2$ &   3.82E-03 &- &     2.17E-02 &- \\
        ~ &$32^2$ &   1.23E-03 &     1.64 &   6.59E-03 &     1.72 \\
        1 &$64^2$ &   3.45E-04 &     1.83 &   1.77E-03 &     1.90 \\
        ~ &$128^2$ &   9.07E-05 &    1.93 &   4.52E-04 &     1.97 \\
        ~ &$256^2$ &   2.31E-05 &    1.97 &   1.14E-04 &   1.99 \\
        \hline
        ~ &$16^2$ &     6.31E-04 &- &     3.00E-03 &- \\
        ~ &$32^2$ &     9.01E-05 &     2.81 &     5.80E-04 &     2.37 \\
        2 &$64^2$ &     1.13E-05 &     3.00 &     8.13E-05 &     2.84 \\
        ~ &$128^2$ &     1.41E-06 &     3.00 &     1.07E-05 &     2.93 \\
        ~ &$256^2$ &     1.76E-07 &     3.00 &     1.37E-06 &     2.96  \\
        \hline
        \end{tabular}
        %\vspace{\baselineskip}
    \end{table*}

\begin{figure}[!h]
        \centering
        \includegraphics[width=0.37\textwidth]{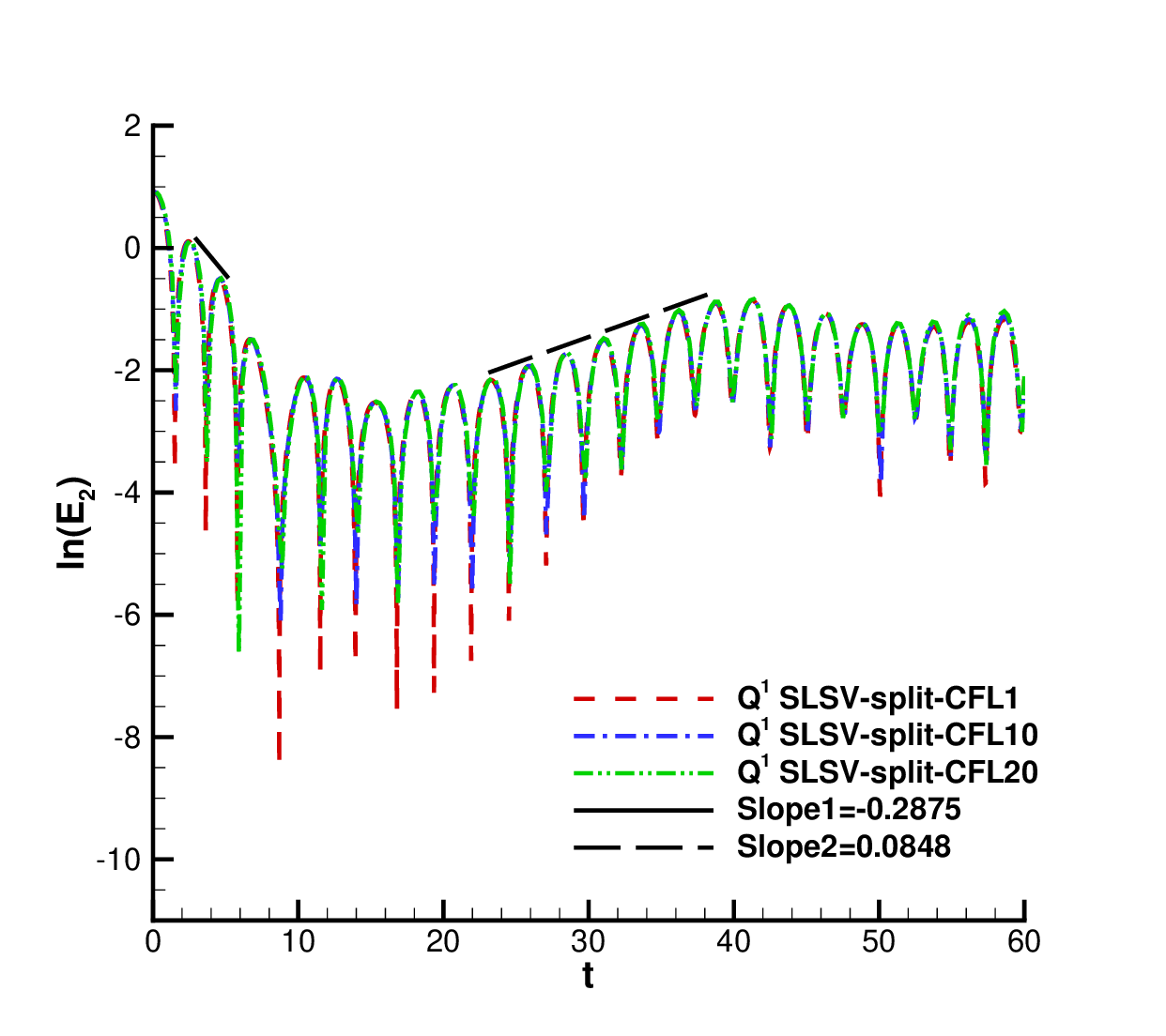}
        \includegraphics[width=0.37\textwidth]{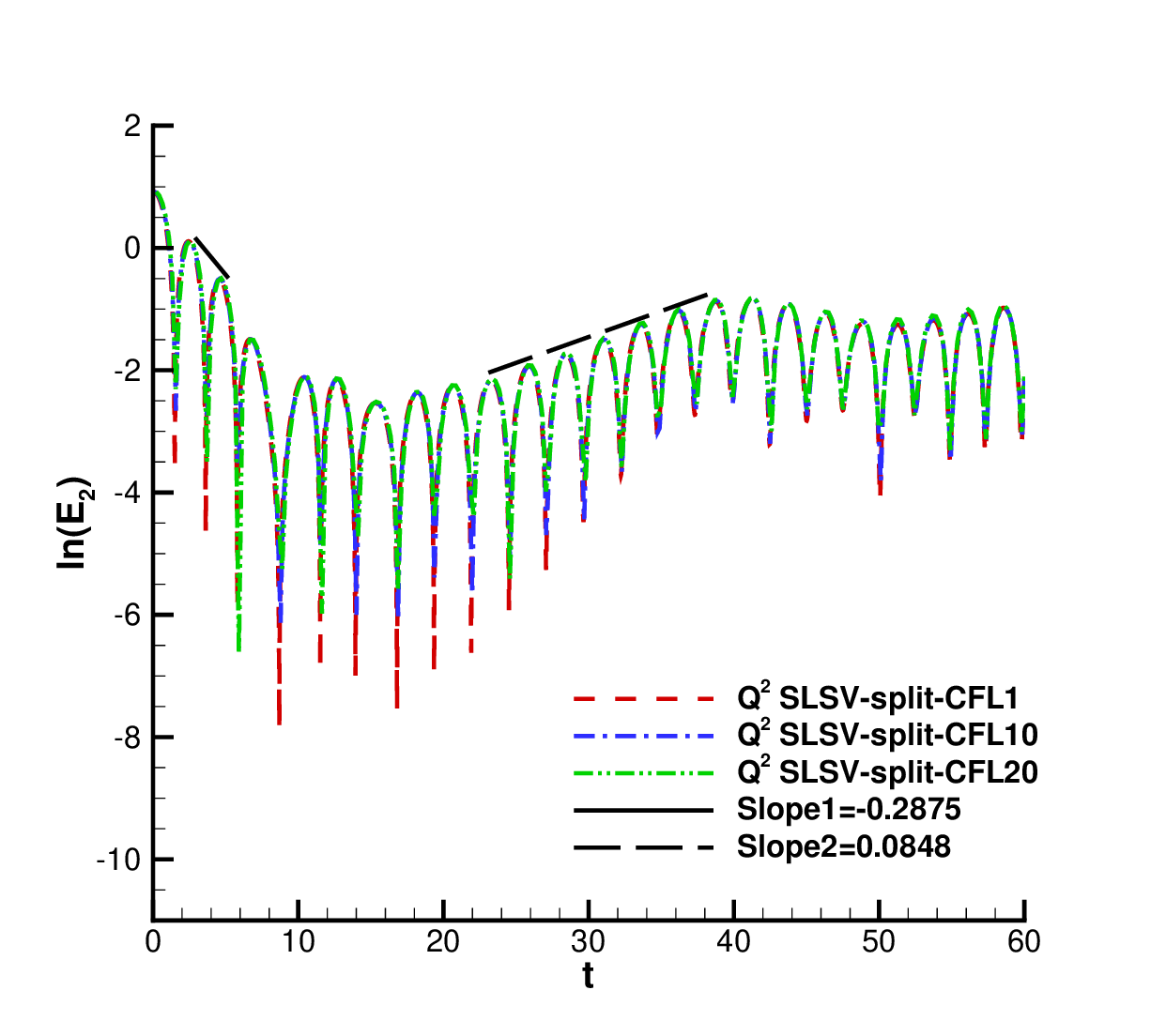}
        \caption{\label{Fig: strong E2} Strong Landau damping: The SLSV schemes are equipped with the PP limiter. Time evolution of the electric field in $L^2$, using a mesh of $160 \times 160$ elements..}
    \end{figure}

\begin{figure}[!h]
        \centering
        \includegraphics[width=0.37\textwidth]{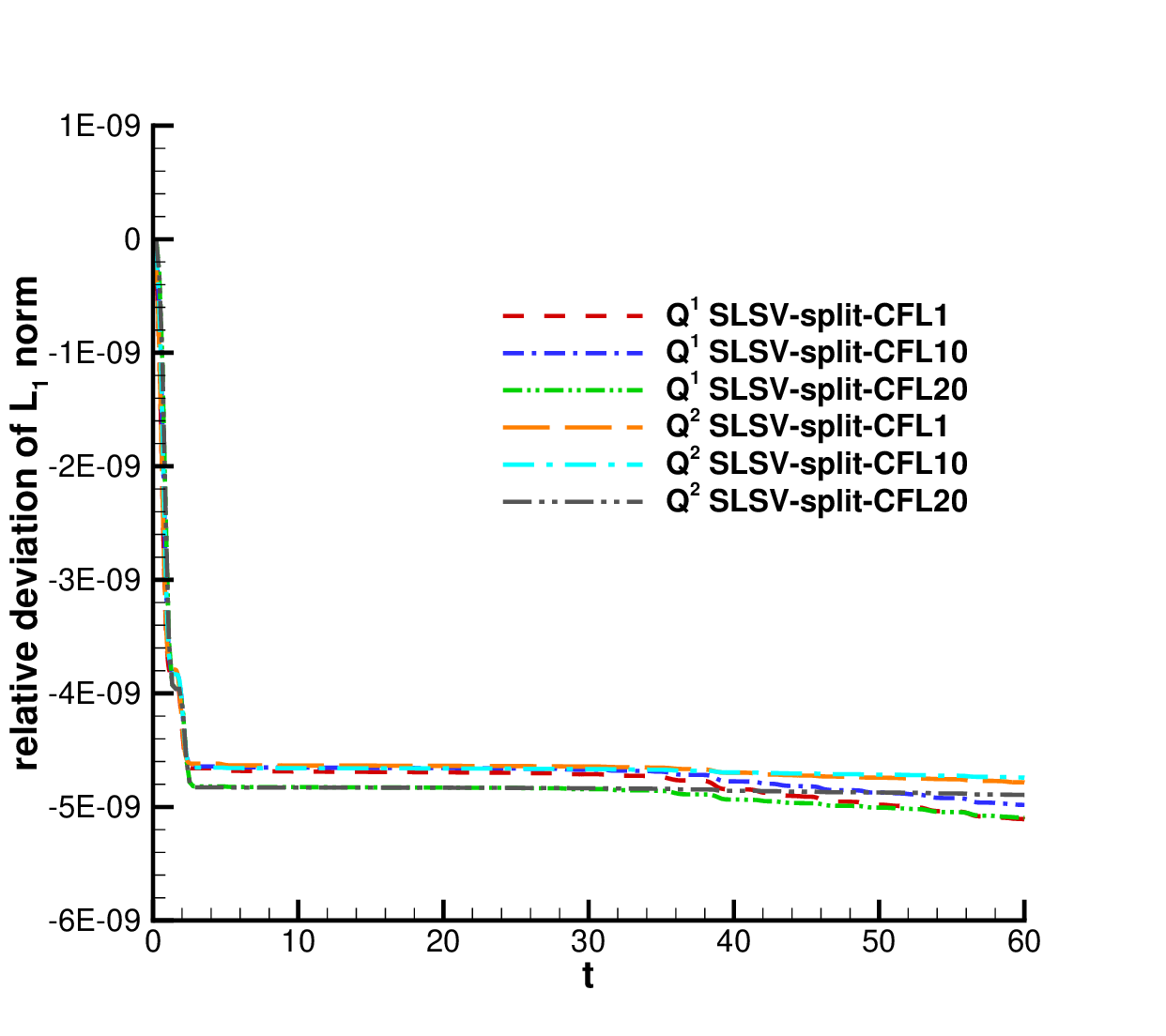}
        \includegraphics[width=0.37\textwidth]{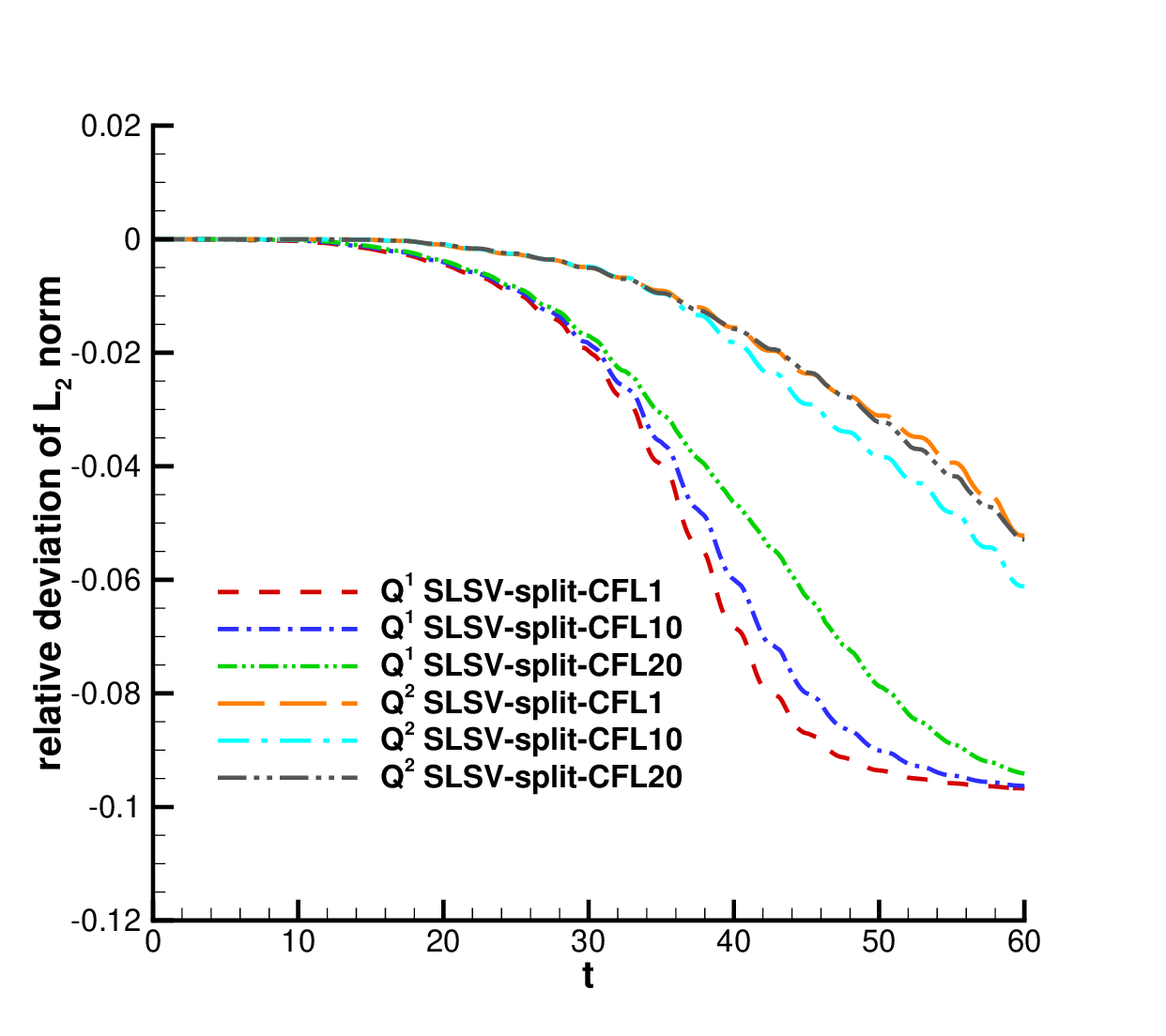}\\
        \includegraphics[width=0.37\textwidth]{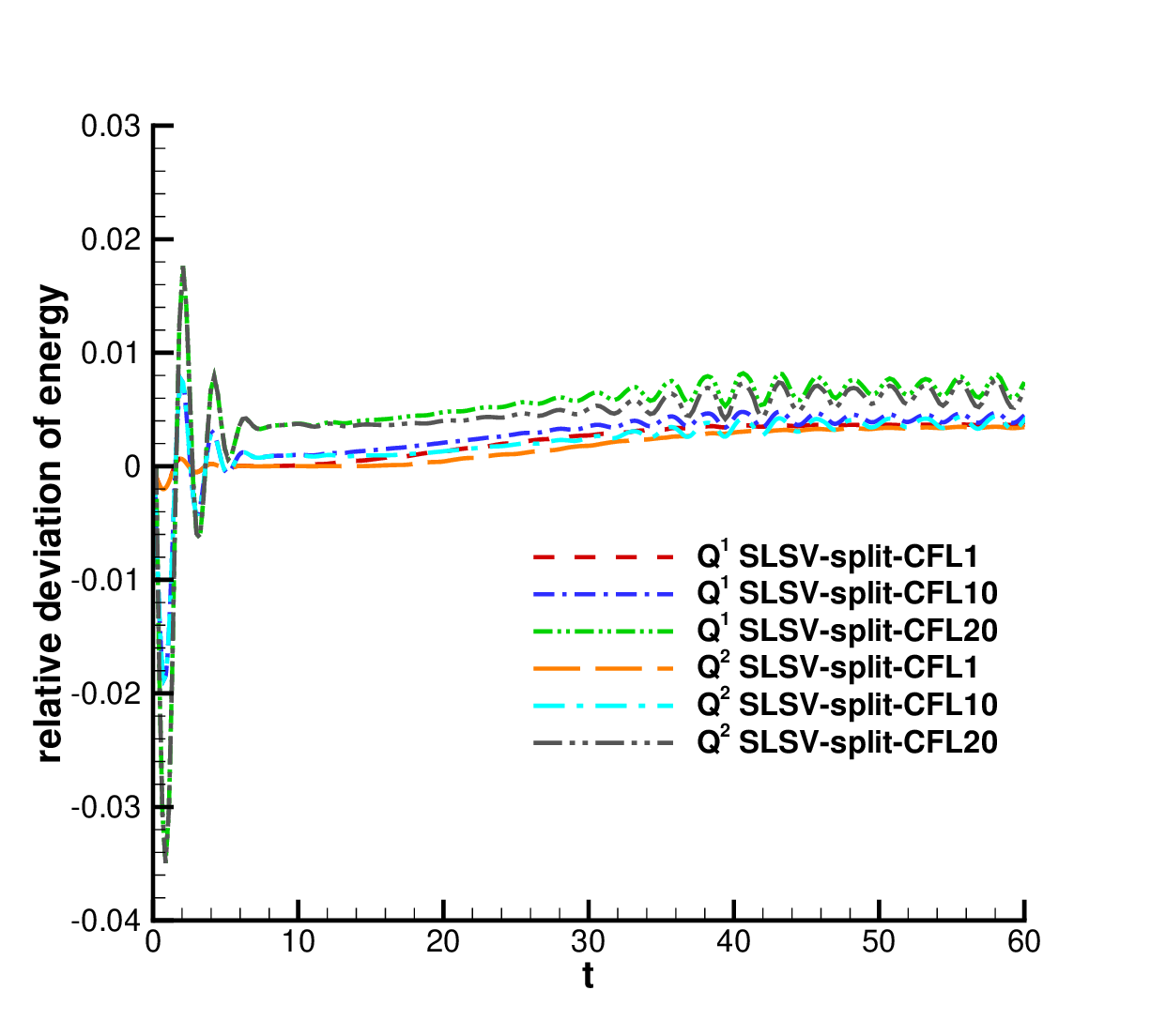}
        \includegraphics[width=0.37\textwidth]{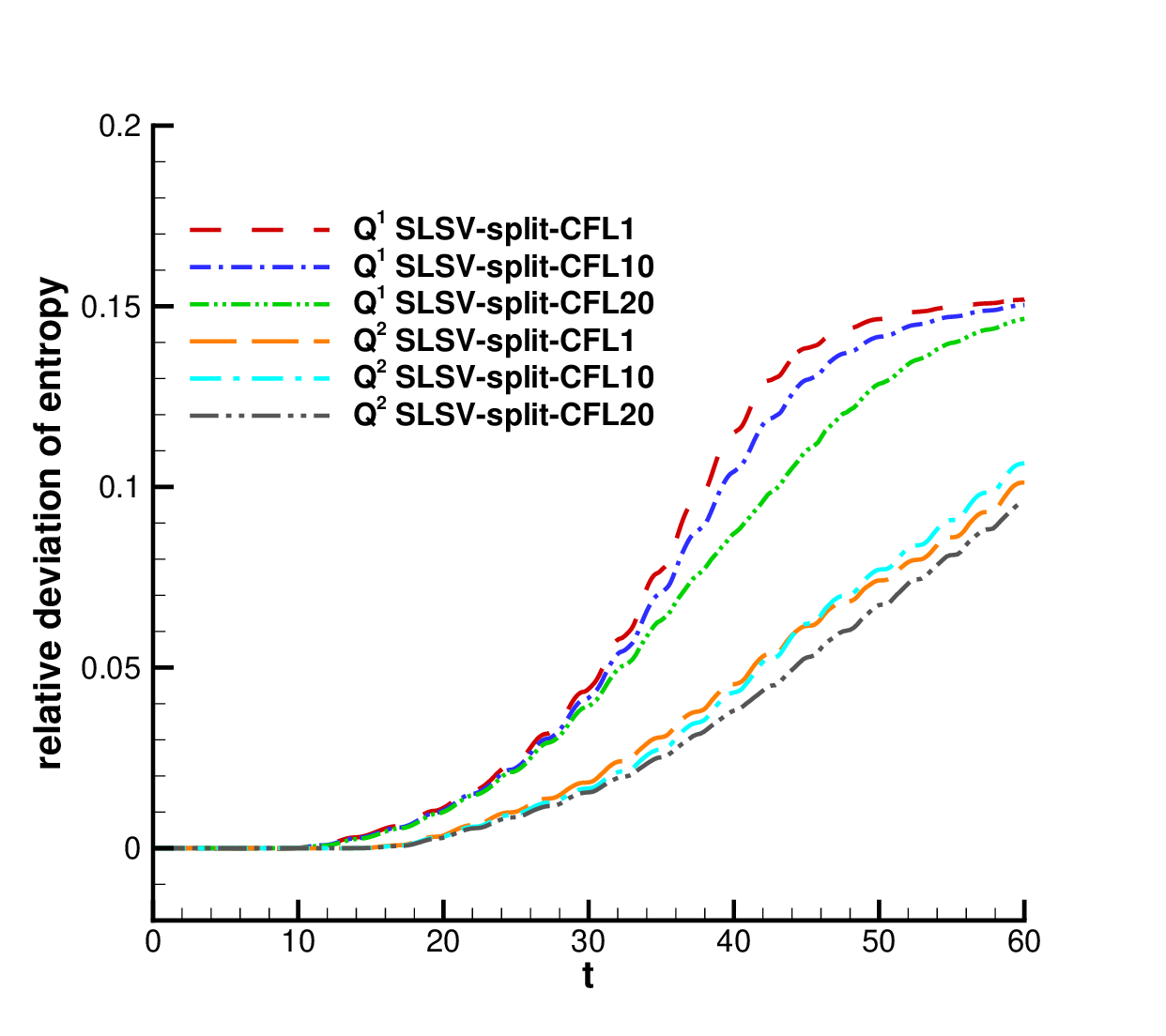}
        \caption{\label{Fig:Strong quantities}  Strong Landau damping. Time evolution of the relative deviations of $L^1$ (upper left) and $L^2$ (upper right) norms of the solution as well as the discrete kinetic energy (lower left) and entropy (lower right), using a mesh of $160 \times 160$ elements.}
    \end{figure}

\begin{figure}[!h]
        \centering
        \includegraphics[width=0.37\textwidth]{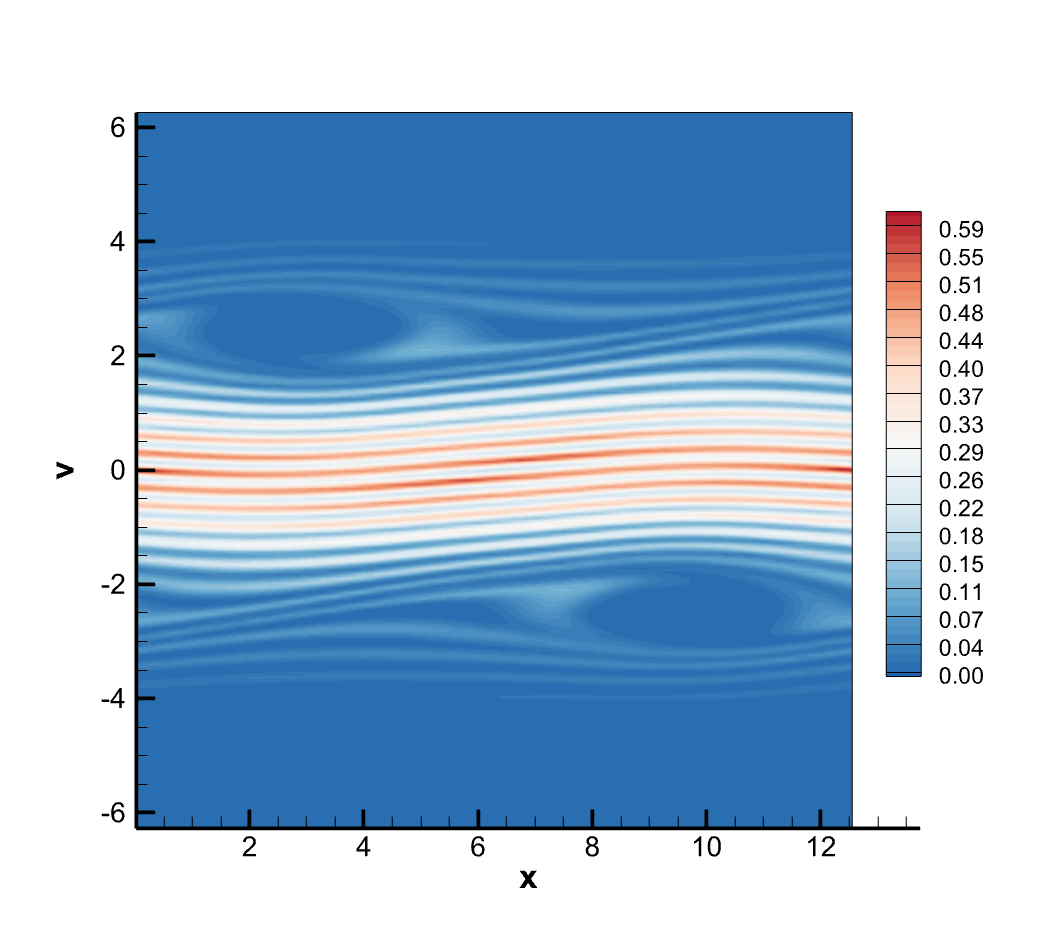}
        \includegraphics[width=0.37\textwidth]{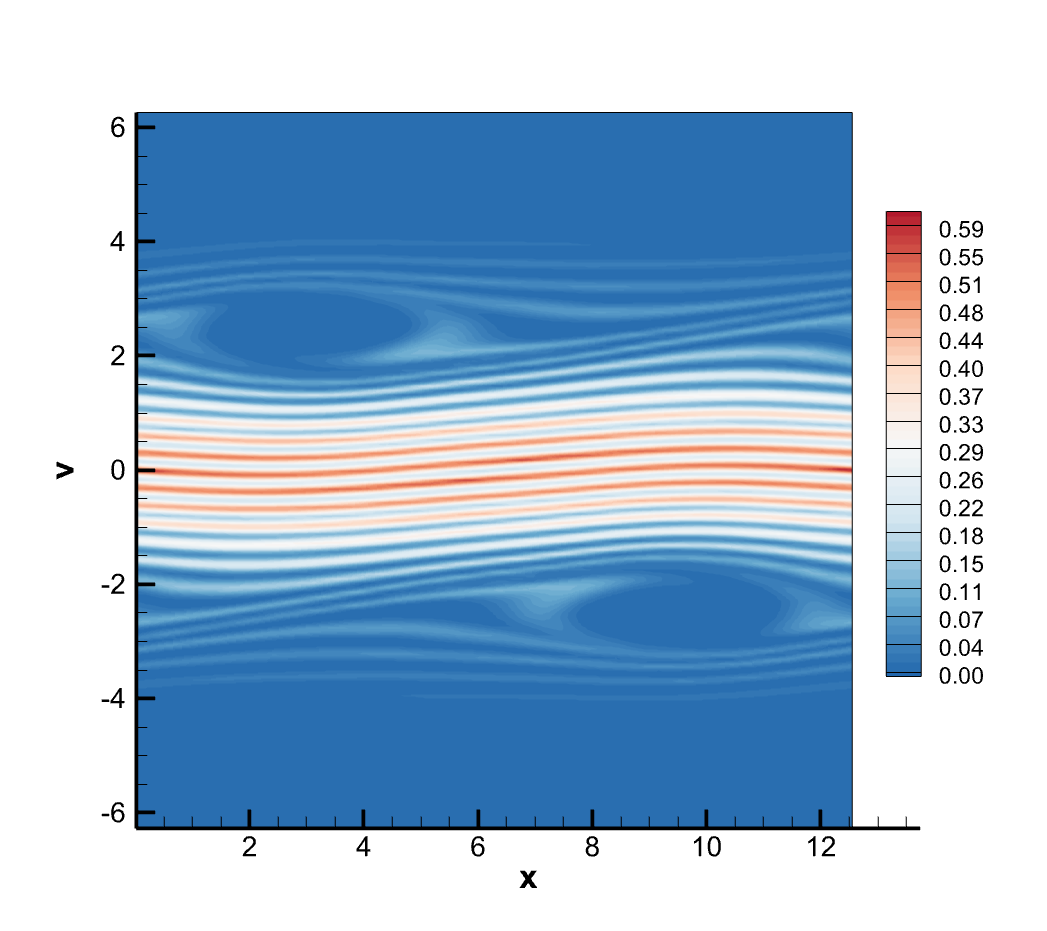}\\
        \includegraphics[width=0.37\textwidth]{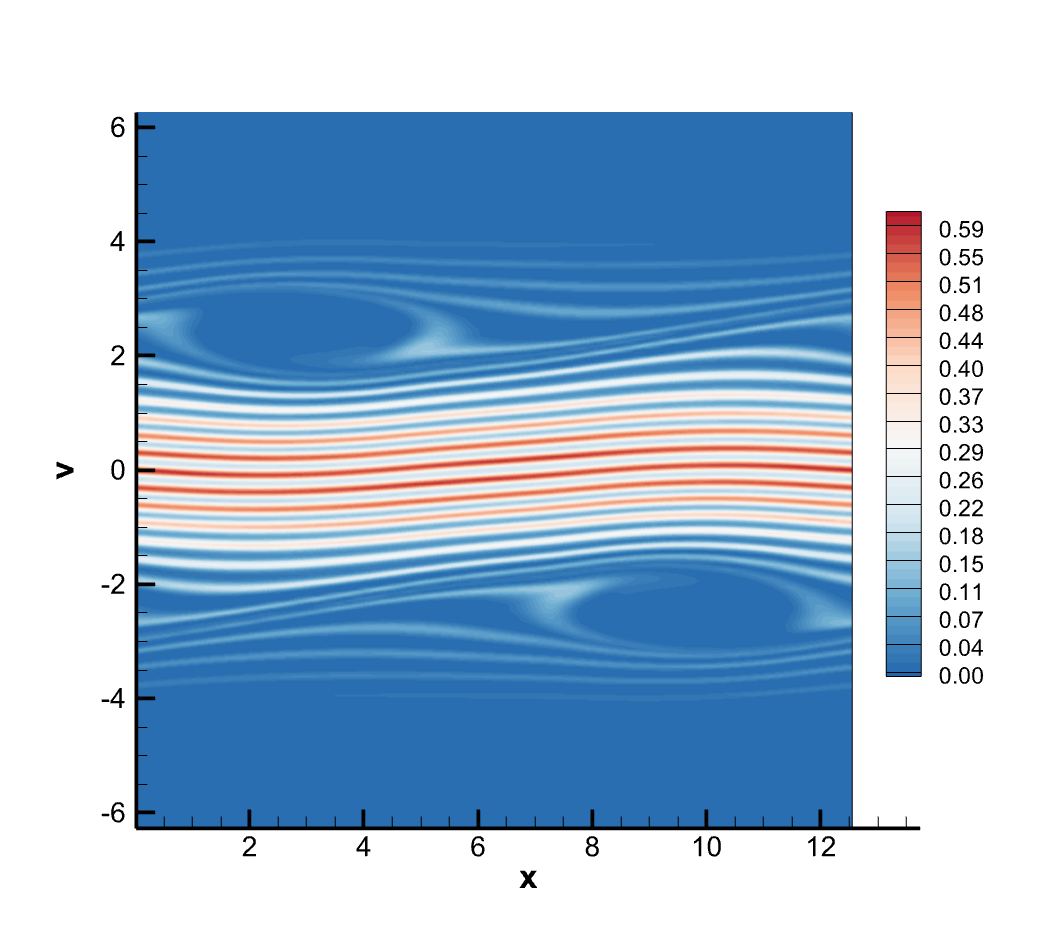}
        \includegraphics[width=0.37\textwidth]{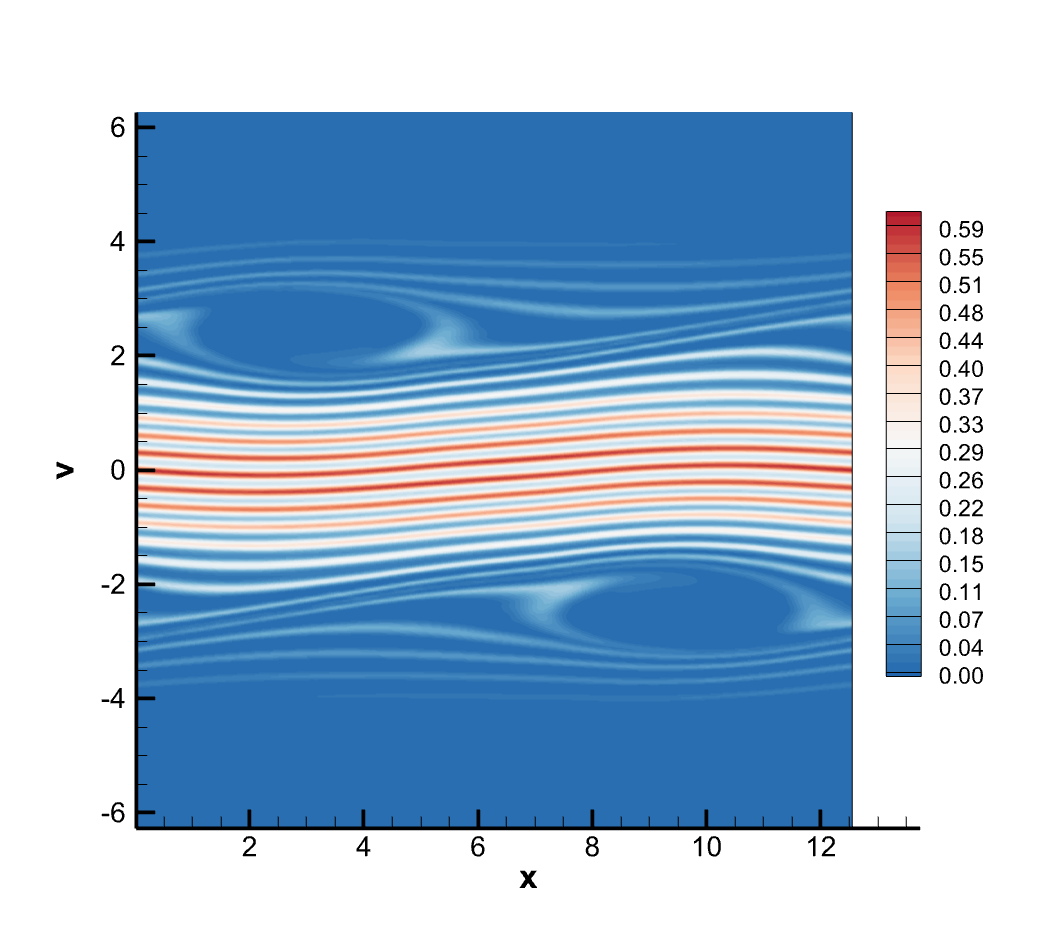}
        \caption{\label{Fig: strong field}  Strong Landau damping with the spatial mesh of $160 \times 160$ at $T=40$. Left:  $CFL = 10$. Right: $ CFL = 20$. Upper: ${\mathcal Q}^1$ SLSV. Bottom: ${\mathcal Q}^2$ SLSV.}
    \end{figure}

    We show the time evolution of the electric field in the $L^2$ norm (in semi-log scale) in Figure \ref{Fig: strong E2} and decay rates are computed by sampling the solution at data points. 
    We find the linear decay rate $\gamma_1=-0.2875$ (measured as the slope of a line originating from the local maximum of the second peak to the third peak), as well as the growth rate $\gamma_2=0.0848$ (measured as the slope of a line originating from the local maximum of the tenth peak to the sixteenth peak), and they agree with the results reported in the literature 
    \cite{VPS1,VPS2,VPS3,2011splitSLDG}.  In addition,  we plot  
    in Figure \ref{Fig:Strong quantities} time evolutions of the relative derivation of the discrete $L^1$ norm, $L^2$ norm, energy and entropy. In particular, we observe that:  (1) The error for the $L^1$ norm (on the order of $10^{-9}$ ) is due to the truncation of the velocity domain, which can be further reduced by using a larger velocity domain in simulations; (2) In general, the ${\mathcal Q}^2$ SLSV method does a better job in conserving these physical norms than the ${\mathcal Q}^1$ SLDG method; (3) Compared to the SLSV schemes with larger CFLs, the schemes with smaller CFLs are able to better preserve the energy and entropy, but perform worse in conserving the $L^2$ norm and entropy. In Figure \ref{Fig: strong field}, we present the contour plots of the solutions at $T=40$ computed by the ${\mathcal Q}^2$ SLSV scheme method with the mesh of $160 \times 160$ elements. We observe that the SLSV scheme with CFL $=10$ and  CFL $=20$ generate very consistent numerical results. Therefore, with larger CFL condition, the SLSV method is still stable and able to generate decent results: the main structures of the solution are captured.

   \begin{example}(Two stream instability I.)
    \end{example}  
    Consider two stream instability, with an unstable initial distribution function,
        \begin{equation*}\label{two stream1 initial}
            f(x, v, t=0)=\frac{2}{7 \sqrt{2 \pi}}\left(1+5 v^2\right)(1+\frac {\alpha((\cos (2 k x)+\cos (3 k x))} {1.2}+\cos (k x))) \exp \left(-\frac{v^2}{2}\right),
        \end{equation*}
    with $\alpha=0.01, k=0.5$ on the computation domain $[0,4 \pi] \times[-10,10]$ and the periodic boundary condition.

    We first present the spatial convergence of the SLSV schemes in Table \ref{Tab: Two stream1 accruate}. We  adopt the time reversibility of VP system to test the order of convergence. Slightly less than $(k+1)$th order accuracy is observed for the ${\mathcal Q}^k,k=1,2$, SLSV scheme. 
    % This example is more difficult to calculate than the previous ones.
    \begin{table*}[!ht]
        \centering
        \footnotesize
        \caption{\label{Tab: Two stream1 accruate}Two stream instability I: errors and convergence orders of $L^2$ and $L^{\infty}$   at $T=0.5$ with $CFL=0.1$.} 
        \vspace{0.4em} \centering%  xiaowu
        \begin{tabular}{c|c|cccc}%{\textwidth}
        %\toprule[0.001pt]
        \hline
        \multirow{2}{*}{$k$}   &\multirow{2}{*}{$N^2$}   & \multicolumn{4}{c}{ $CFL=0.1$}   \\
        \cline{3-6}
        & &$L^2$ &Order&$L^{\infty}$& Order\\
        %\midrule[0.001pt]
        \hline
        ~ &$16^2$ &     2.80E-03 &- &     8.77E-03 &- \\
        ~ &$32^2$ &     3.37E-04 &     3.06 &     1.88E-03 &     2.22 \\
        1 &$64^2$ &     9.74E-05 &     1.79 &     5.77E-04 &     1.71 \\
        ~ &$128^2$ &     4.63E-05 &     1.07 &     2.94E-04 &     0.98 \\
        ~ &$256^2$ &     2.98E-05 &     0.64 &     1.81E-04 &     0.70 \\
        \hline
        ~ &$16^2$ &     2.02E-04 &- &     1.01E-03 &- \\
        ~ &$32^2$ &     1.28E-04 &     0.66 &     9.27E-04 &     0.12 \\
        2 &$64^2$ &     3.34E-05 &     1.94 &     2.50E-04 &     1.89 \\
        ~ &$128^2$ &     7.01E-06 &     2.25 &     5.30E-05 &     2.23 \\
        ~ &$256^2$ &     1.28E-06 &     2.46 &     9.44E-06 &     2.49 \\
        \hline
        \end{tabular}
        %\vspace{\baselineskip}
    \end{table*}

\begin{figure}[!h]
            \centering
            \includegraphics[width=0.37\textwidth]{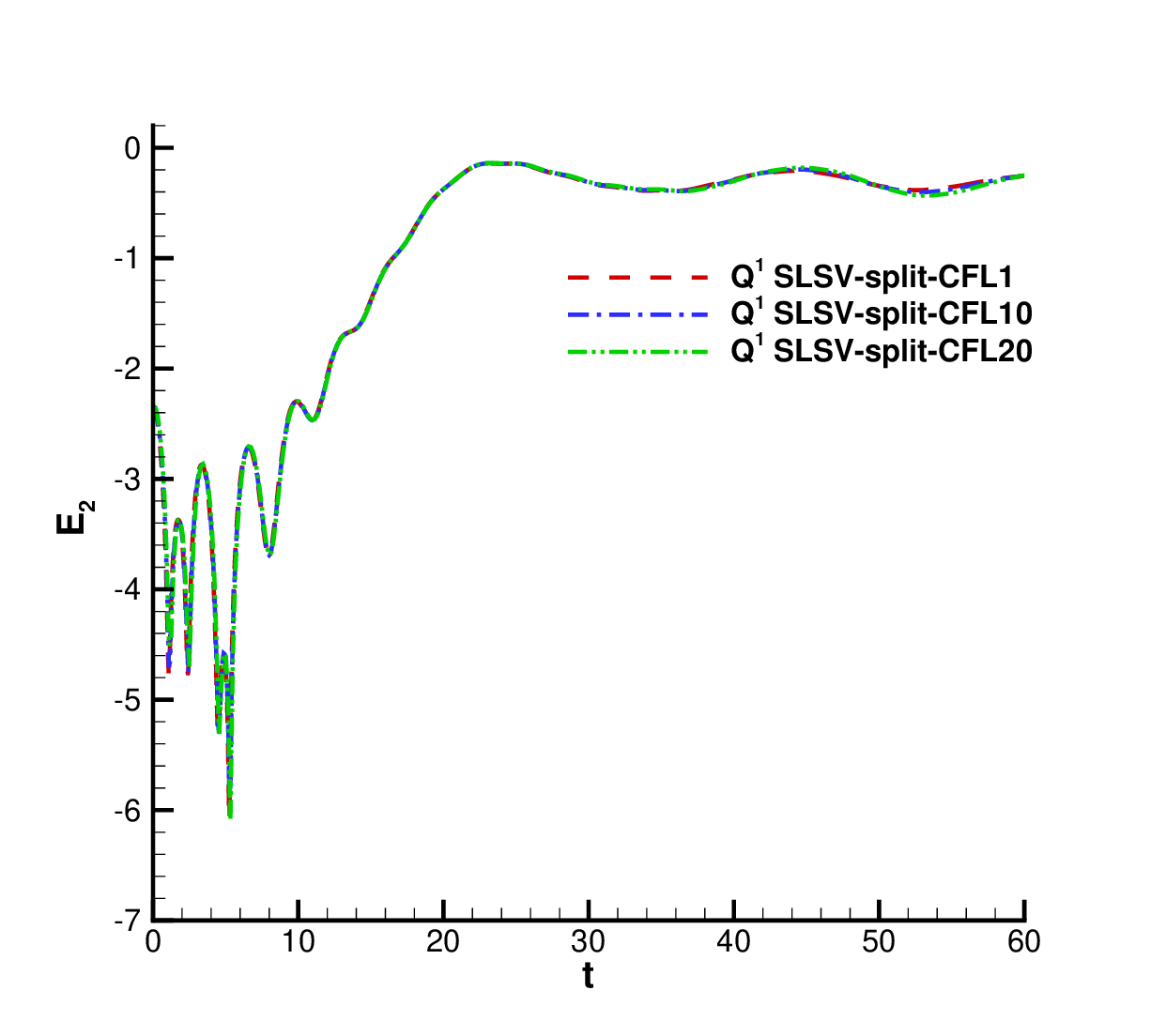}
            \includegraphics[width=0.37\textwidth]{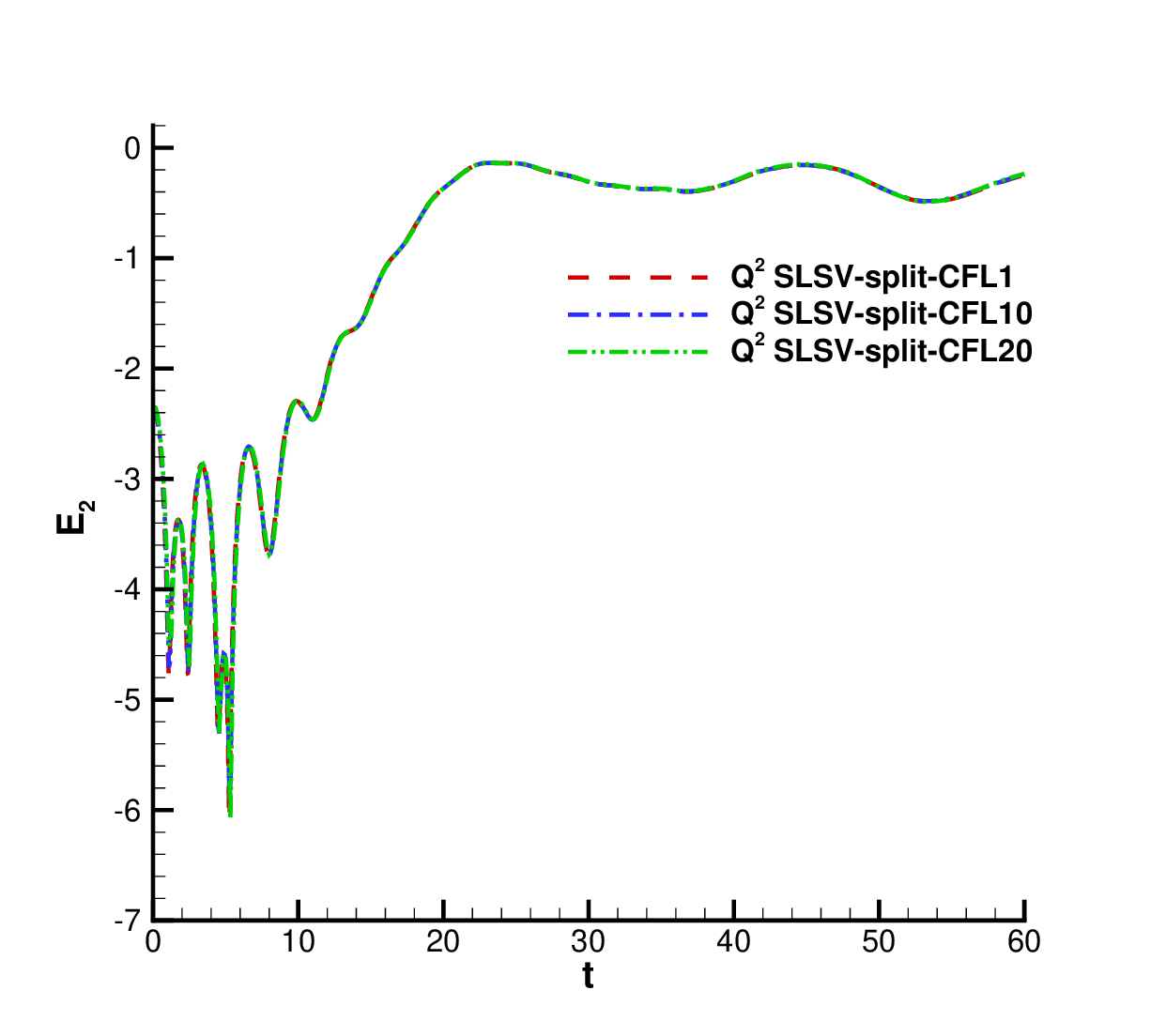}
            \caption{\label{Fig: Two stream1 E2} Two stream I: The SLSV schemes are equipped with the PP limiter. Time evolution of the electric field in $L^2$, using a mesh of $160 \times 160$ elements.}
        \end{figure}

\begin{figure}[!h]
            \centering
            \includegraphics[width=0.37\textwidth]{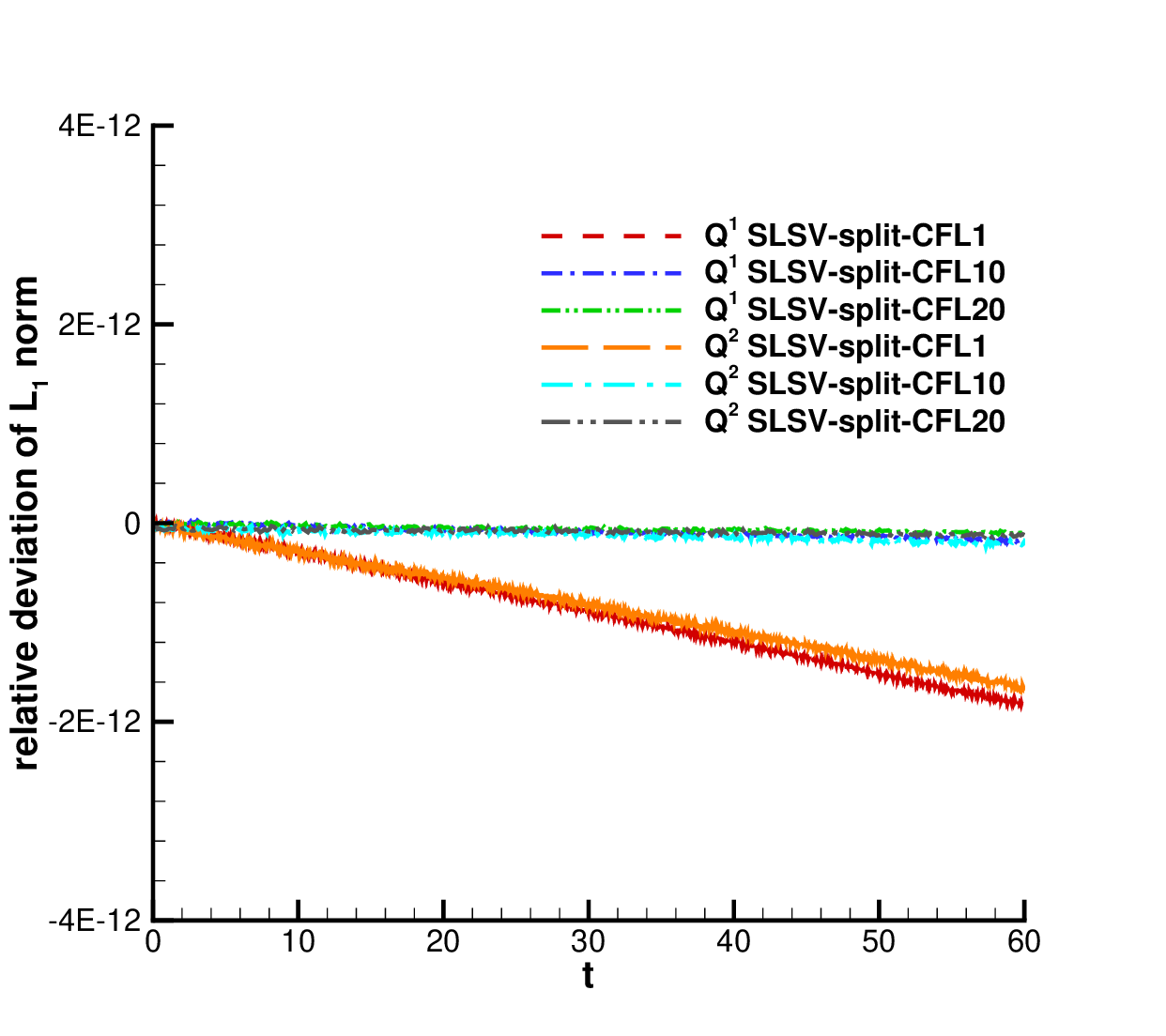}
            \includegraphics[width=0.37\textwidth]{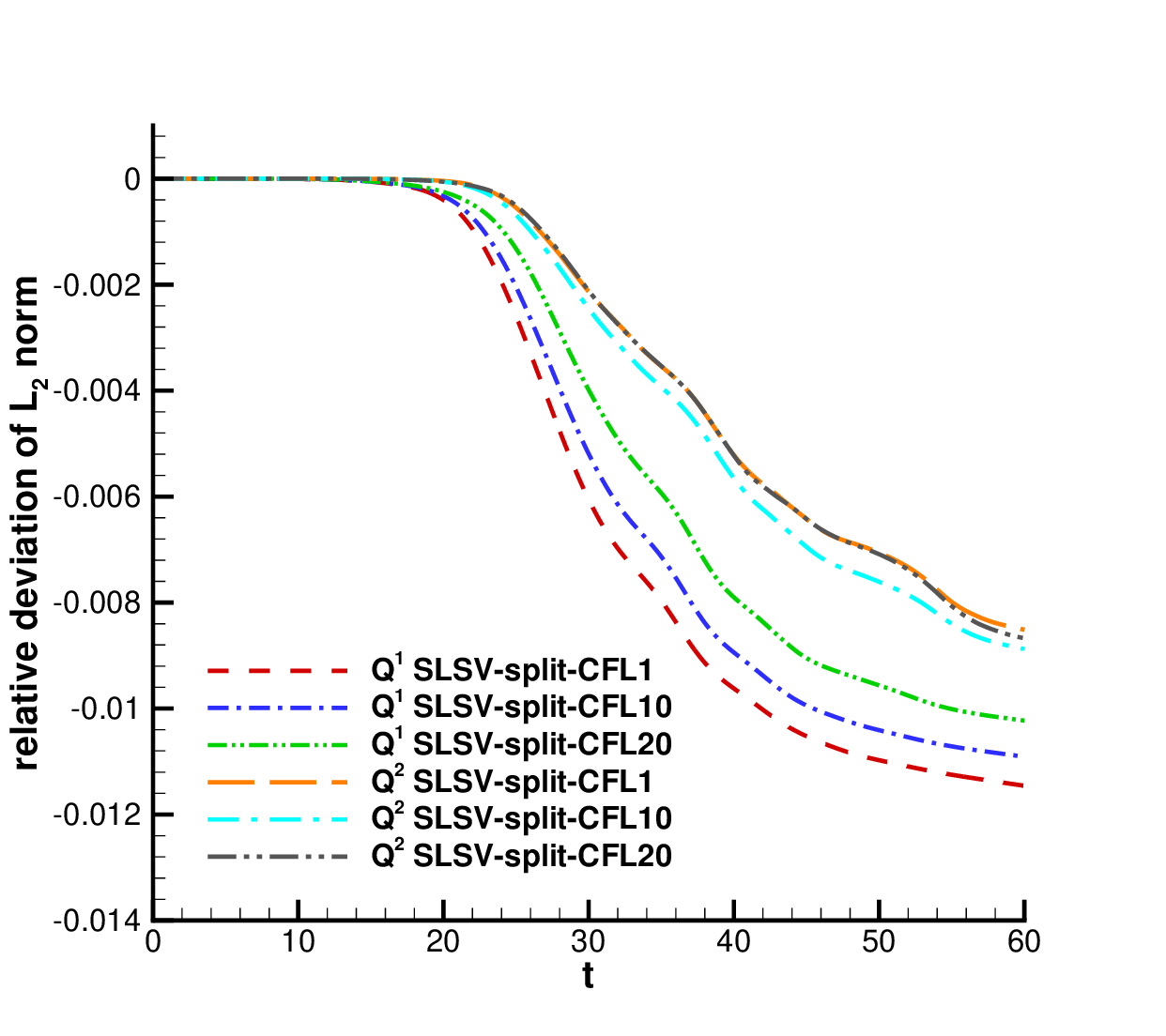}\\
            \includegraphics[width=0.37\textwidth]{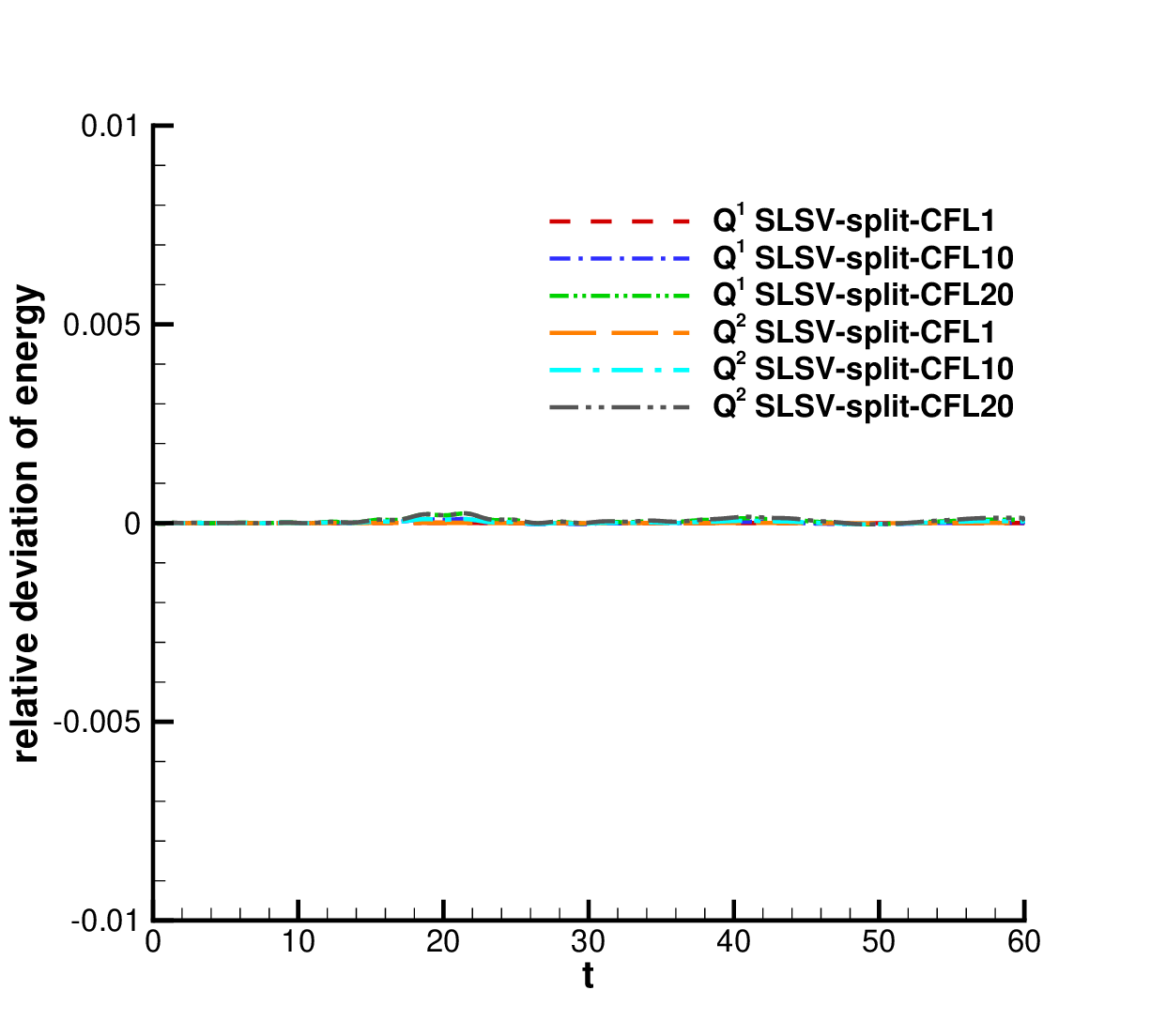}
            \includegraphics[width=0.37\textwidth]{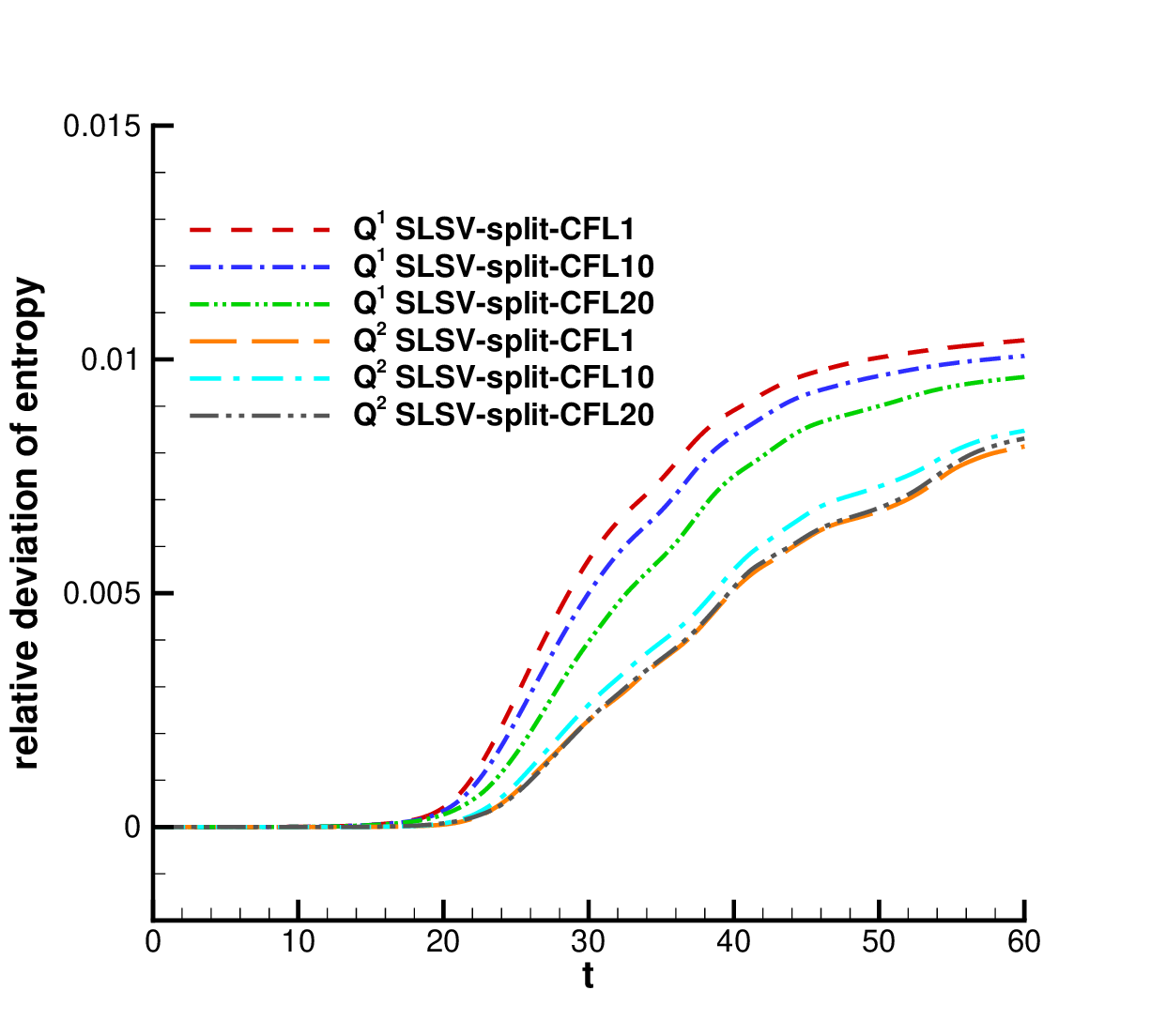}
            \caption{\label{Fig: Two stream1 quantities}  Two stream I. Time evolution of the relative deviations of $L^1$ (upper left) and $L^2$ (upper right) norms of the solution as well as the discrete kinetic energy (lower left) and entropy (lower right), using a mesh of $160 \times 160$ elements.}
        \end{figure}

\begin{figure}[!h]
            \centering
            \includegraphics[width=0.37\textwidth]{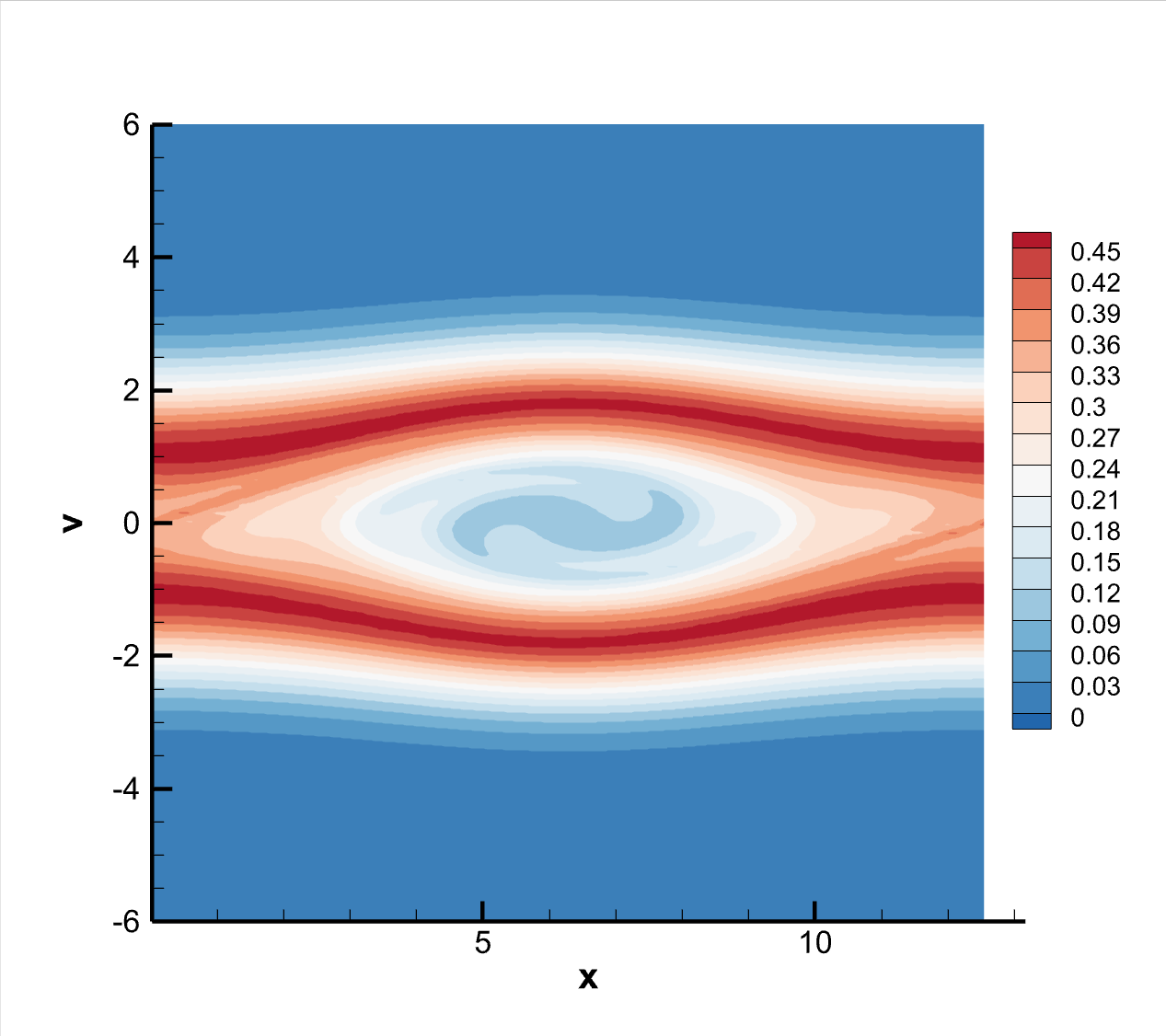}
            \includegraphics[width=0.37\textwidth]{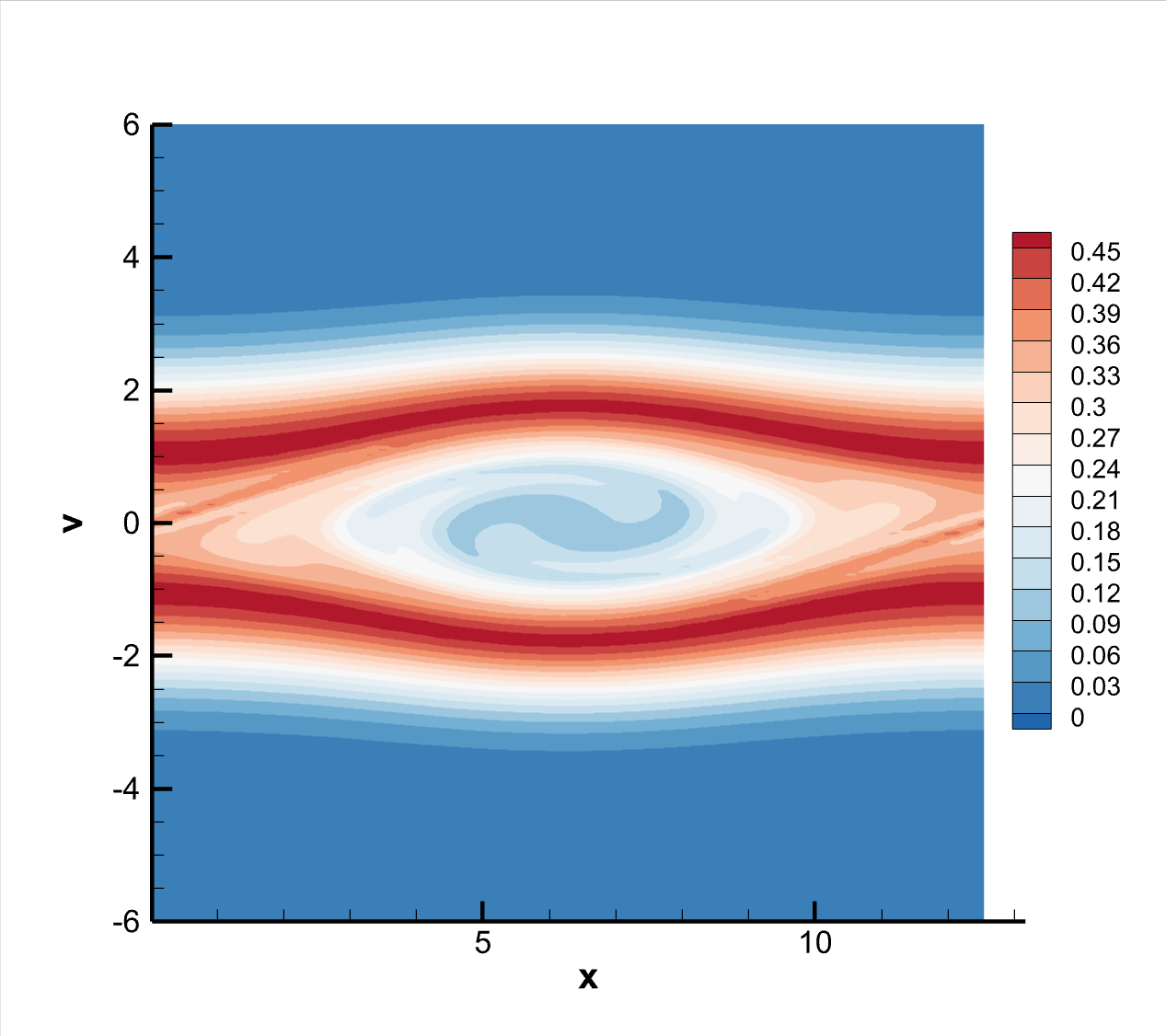}\\
            \includegraphics[width=0.37\textwidth]{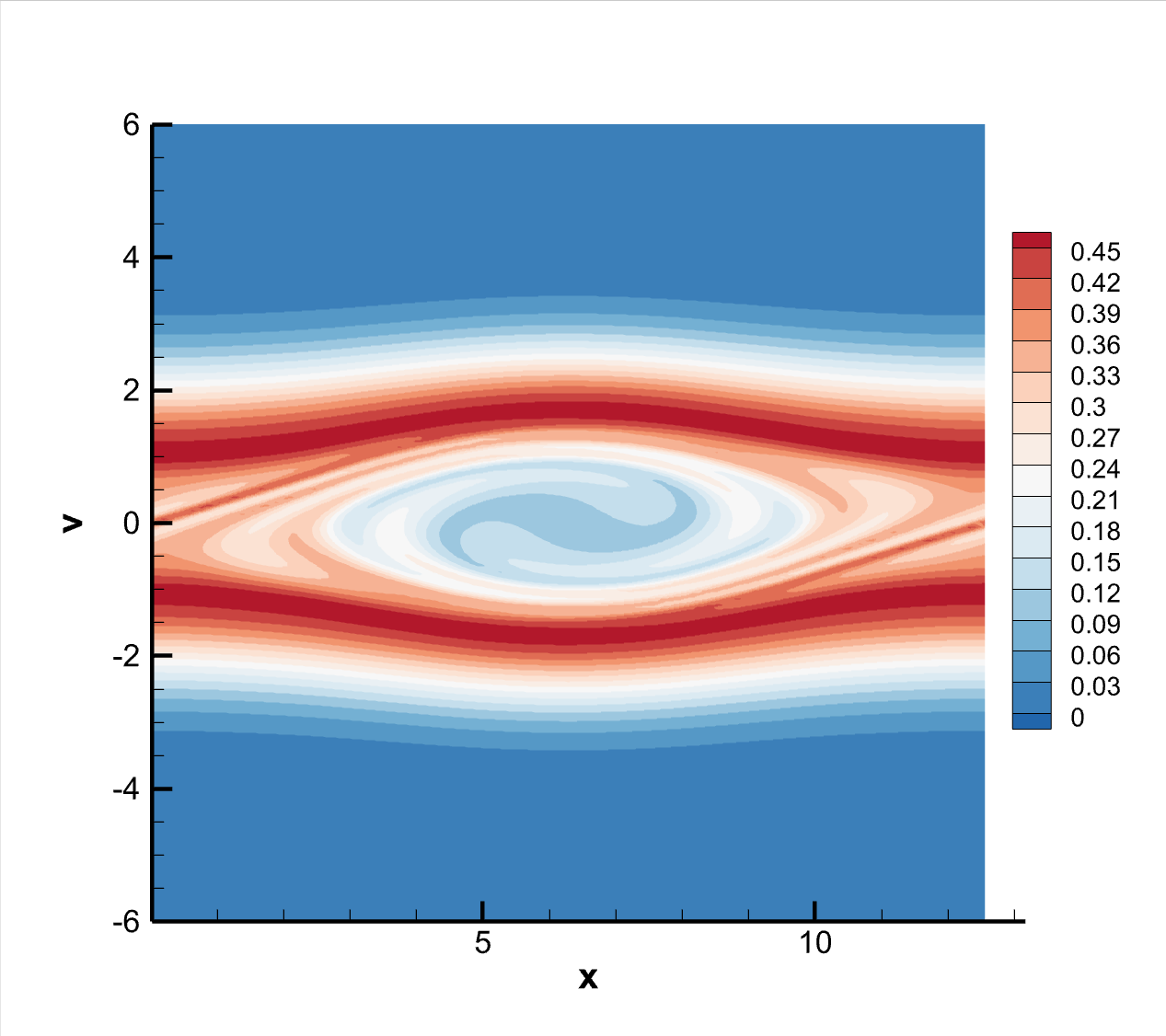}
            \includegraphics[width=0.37\textwidth]{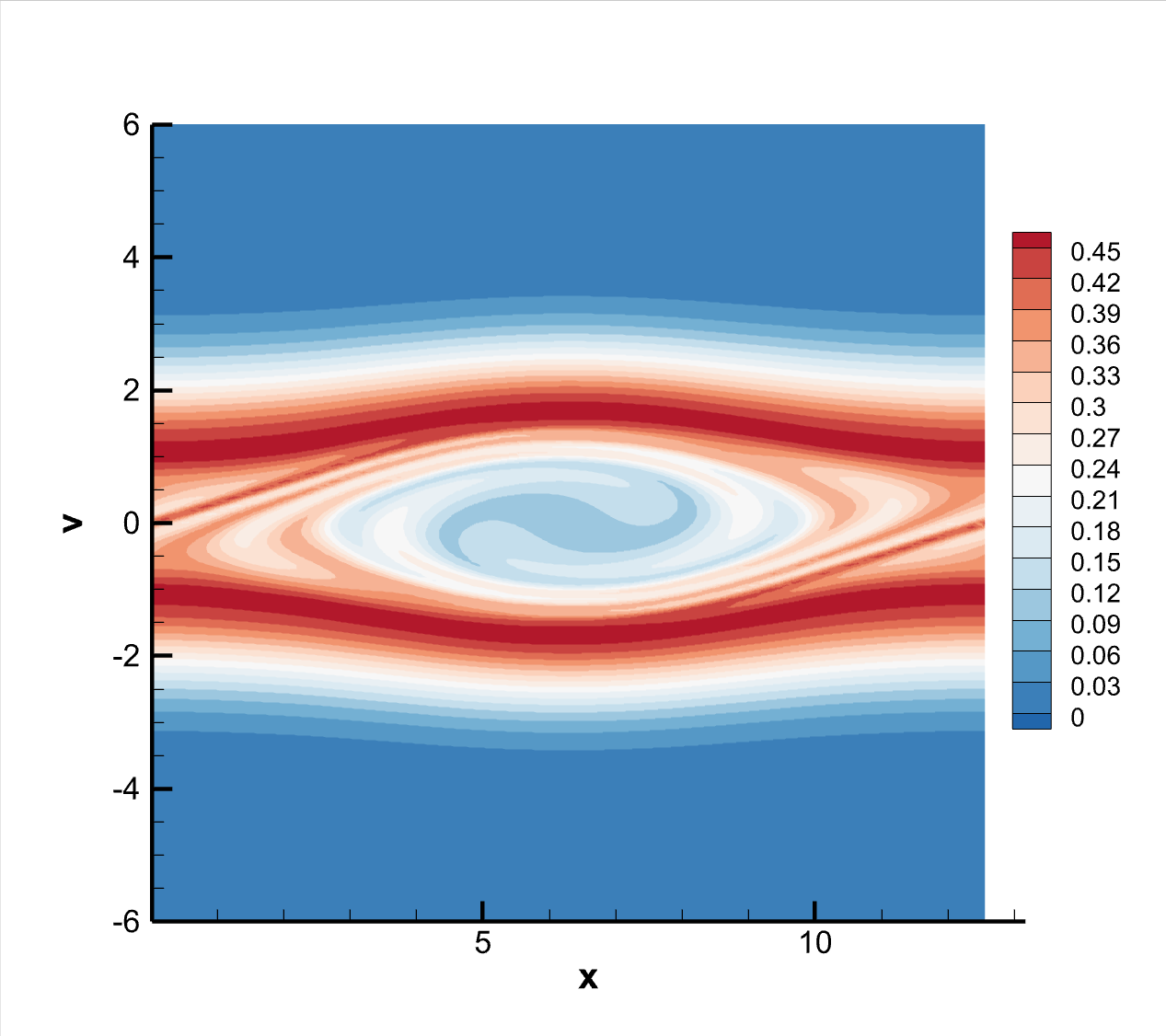}
            \caption{\label{Fig: Two stream1 field} Two stream I with the spatial mesh of $160 \times 160  at T=53$. Left: $ CFL = 10$. Right:  $CFL = 20$. Upper: ${\mathcal Q}^1$ SLSV. Bottom: ${\mathcal Q}^2$ SLSV.}
        \end{figure}

    We show in Figure \ref{Fig: Two stream1 E2}
    the time evolution of the electric field in the $L^2$ norm (in semi-log scale) for the ${\mathcal Q}^1$ SLSV and ${\mathcal Q}^2$ SLSV schemes using a mesh of $160 \times 160$ elements and different CFLs, and in Figure \ref{Fig: Two stream1 quantities} the relative derivation of the discrete $L^1$ norm, $L^2$ norm, energy and entropy, from which 
    we observe that all methods are able to conserve the $L^1$ norm up to the truncation error from the velocity domain. 
    The ability of SLSV methods to conserve these physical norms is satisfactory. In Figure \ref{Fig: Two stream1 field}, we plot the numerical solutions of phase space profiles at $T=53$.

%******************************************************
    \begin{example}(Two stream instability II.)
    \end{example} 
    Consider the symmetric two stream instability \cite{twosteamII2,twosteamII1}, with the perturbed equilibrium as the initial condition
        \begin{equation*}
            f(x, v, t=0)=\frac{1}{2 v_t \sqrt{2 \pi}}\left[\exp \left(-\frac{(v-u)^2}{2 v_{t h}^2}\right)+\exp \left(-\frac{(v+u)^2}{2 v_{t h}^2}\right)\right](1+0.05 \cos (k x)),
        \end{equation*}
    where $u=0.99, k=\frac{2}{13}$, and $v_t=0.3$. We let $v_{\max }=2 \pi$. 

    We first present the spatial convergence of the SLSV scheme in Table \ref{Tab: two stream2 accruate}. The order of convergence is examined using the same approach as for the strong Landau damping, i.e., the time reversibility of the VP system. 
    \begin{table*}[!ht]
        \centering
        \footnotesize
        \caption{\label{Tab: two stream2 accruate}Two stream instability II: errors and convergence orders of $L^2$ and $L^{\infty}$ for SLSV method, at $T=0.5$, with $CFL=0.1$.} 
        \vspace{0.4em} \centering%  xiaowu
        \begin{tabular}{c|c|cccc}%{\textwidth}
        %\toprule[0.001pt]
        \hline
        \multirow{2}{*}{$k$}   &\multirow{2}{*}{$N^2$}   & \multicolumn{4}{c}{ $CFL=0.1$}   \\
        \cline{3-6}
        & &$L^2$ &Order&$L^{\infty}$& Order\\
        %\midrule[0.001pt]
        \hline
        ~ &$16^2$ &     3.28E-02 &- &     1.67E-01 &- \\
        ~ &$32^2$ &     1.07E-02 &     1.61 &     6.40E-02 &     1.38 \\
        1 &$64^2$ &     3.73E-03 &     1.53 &     2.37E-02 &     1.44 \\
        ~ &$128^2$ &     1.11E-03 &     1.75 &     6.11E-03 &     1.95 \\
        ~ &$256^2$ &     3.00E-04 &     1.89 &     1.46E-03 &     2.06 \\
        \hline
        ~ &$16^2$ &     1.38E-02 &- &     9.95E-02 &- \\
        ~ &$32^2$ &     2.23E-03 &     2.63 &     1.42E-02 &     2.81 \\
        2 &$64^2$ &     5.22E-04 &     2.09 &     2.97E-03 &     2.26 \\
        ~ &$128^2$ &     6.21E-05 &     3.07 &     4.27E-04 &     2.80 \\
        ~ &$256^2$ &     7.82E-06 &     2.99 &     5.47E-05 &     2.96 \\
        \hline
        \end{tabular}
        %\vspace{\baselineskip}
    \end{table*}

\begin{figure}[!h]
        \centering
        \includegraphics[width=0.37\textwidth]{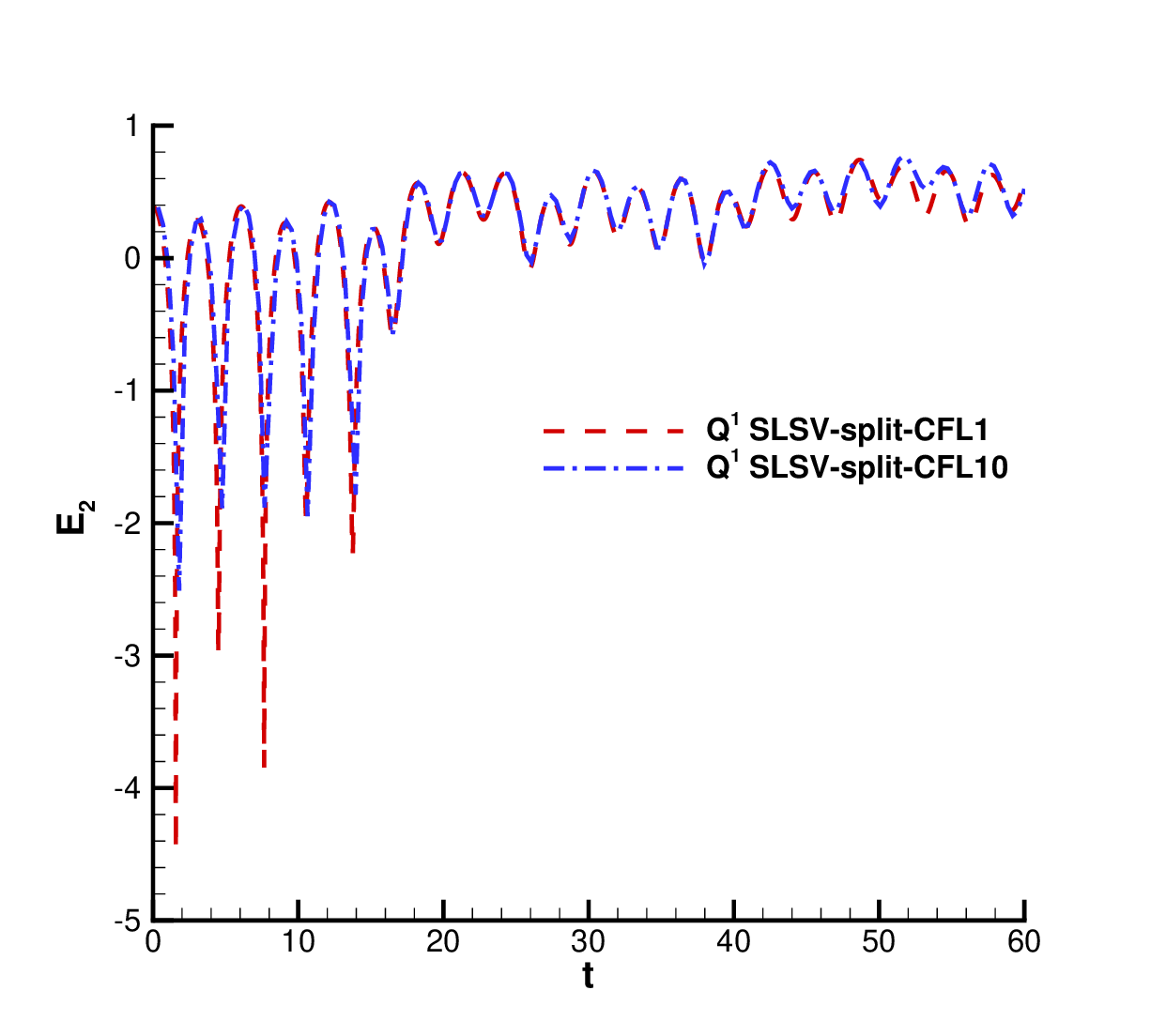}
        \includegraphics[width=0.37\textwidth]{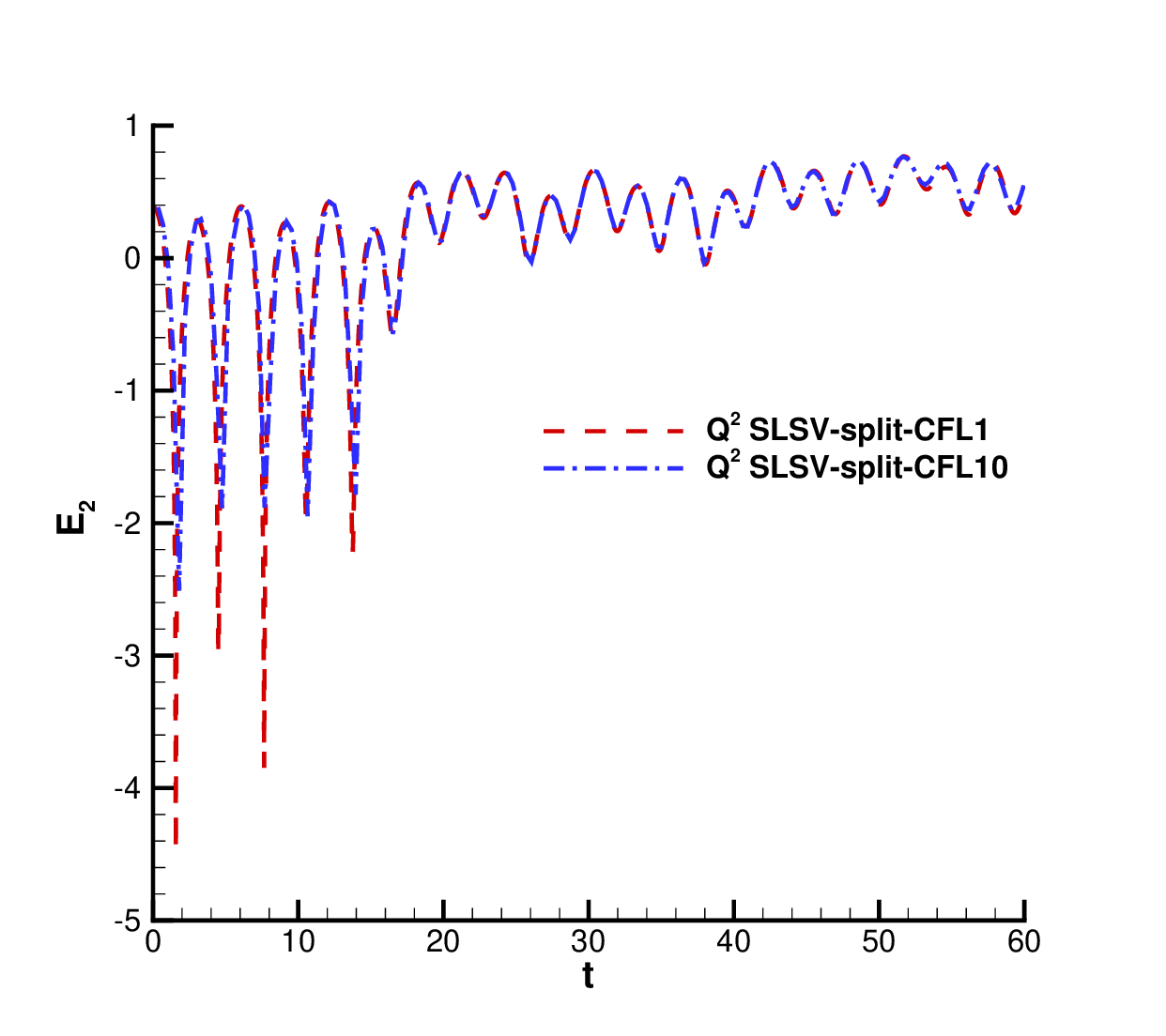}
        \caption{\label{Fig: two stream2 E2} Two stream II: The SLSV schemes are equipped with the PP limiter. Time evolution of the electric field in $L^2$, using a mesh of $160 \times 160$ elements.}
    \end{figure}

\begin{figure}[!h]
        \centering
        \includegraphics[width=0.37\textwidth]{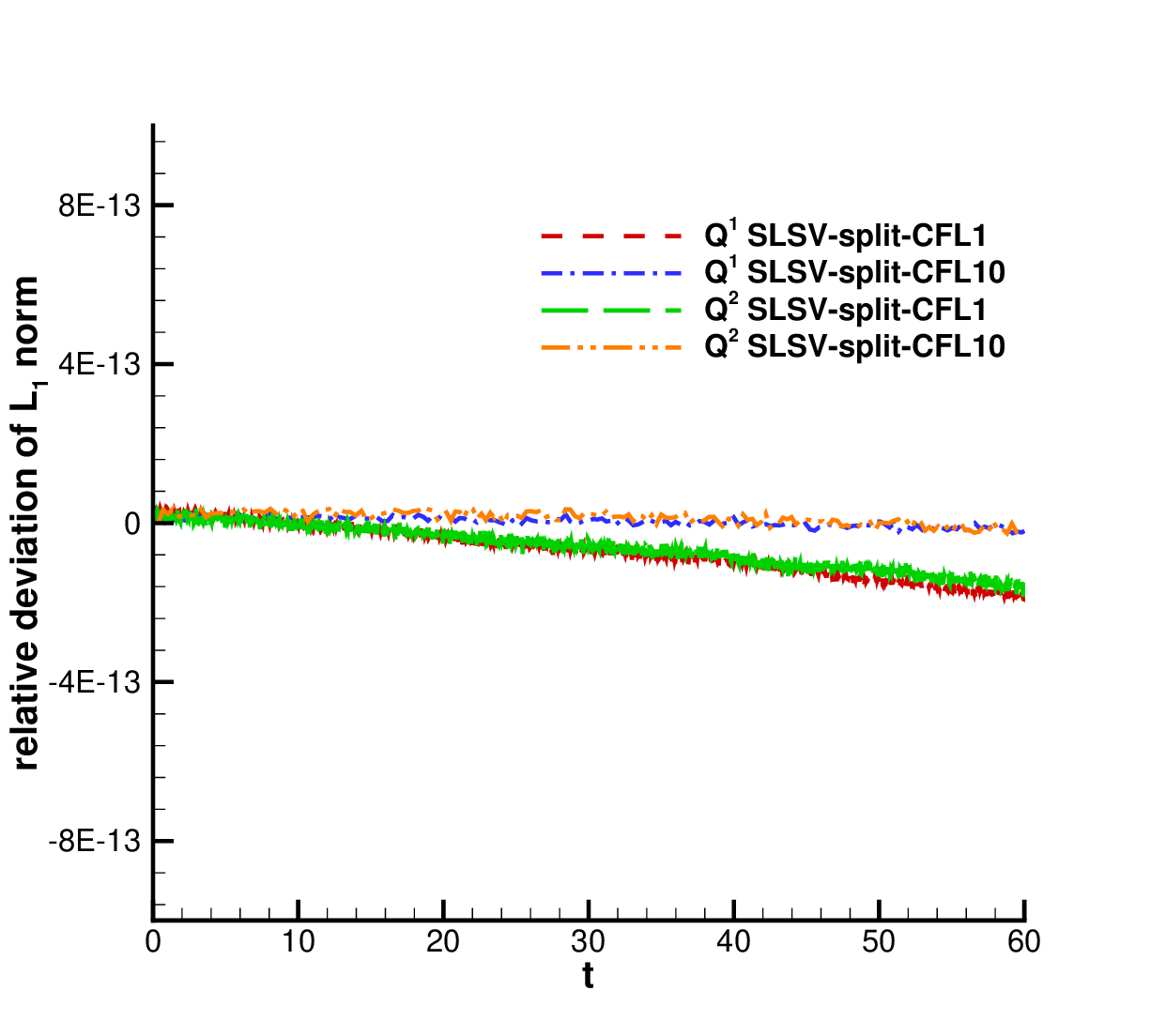}
        \includegraphics[width=0.37\textwidth]{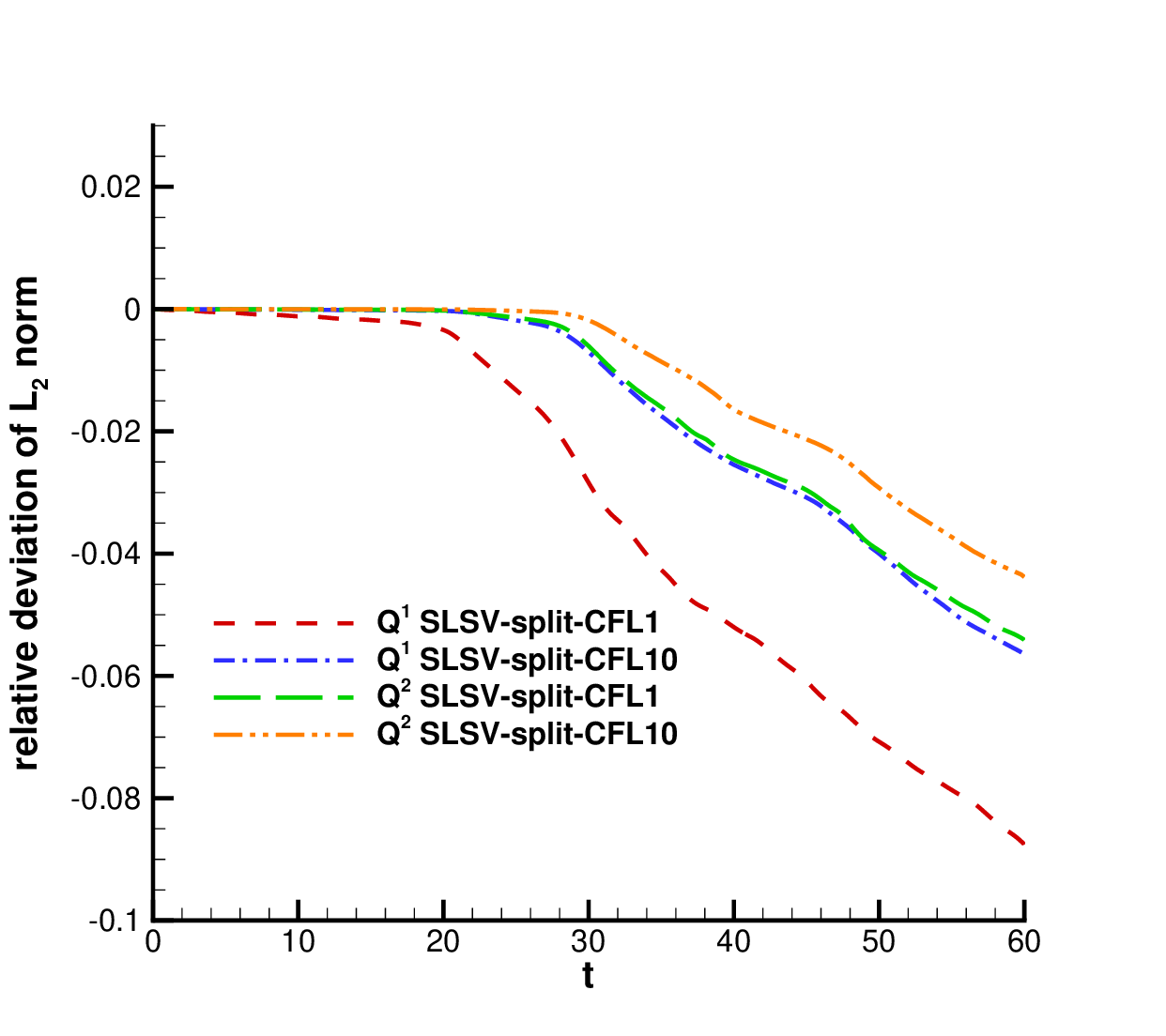}\\
        \includegraphics[width=0.37\textwidth]{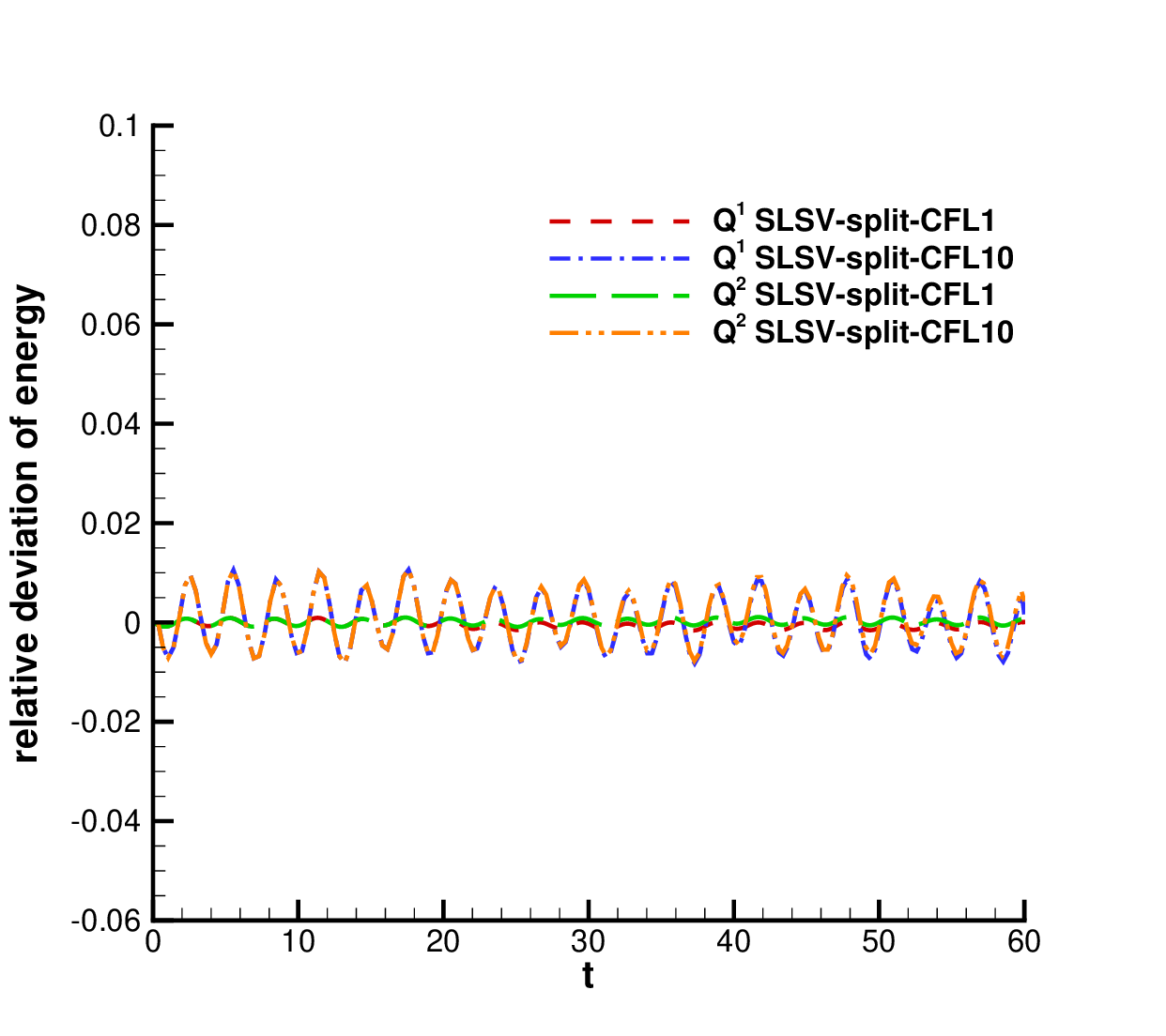}
        \includegraphics[width=0.37\textwidth]{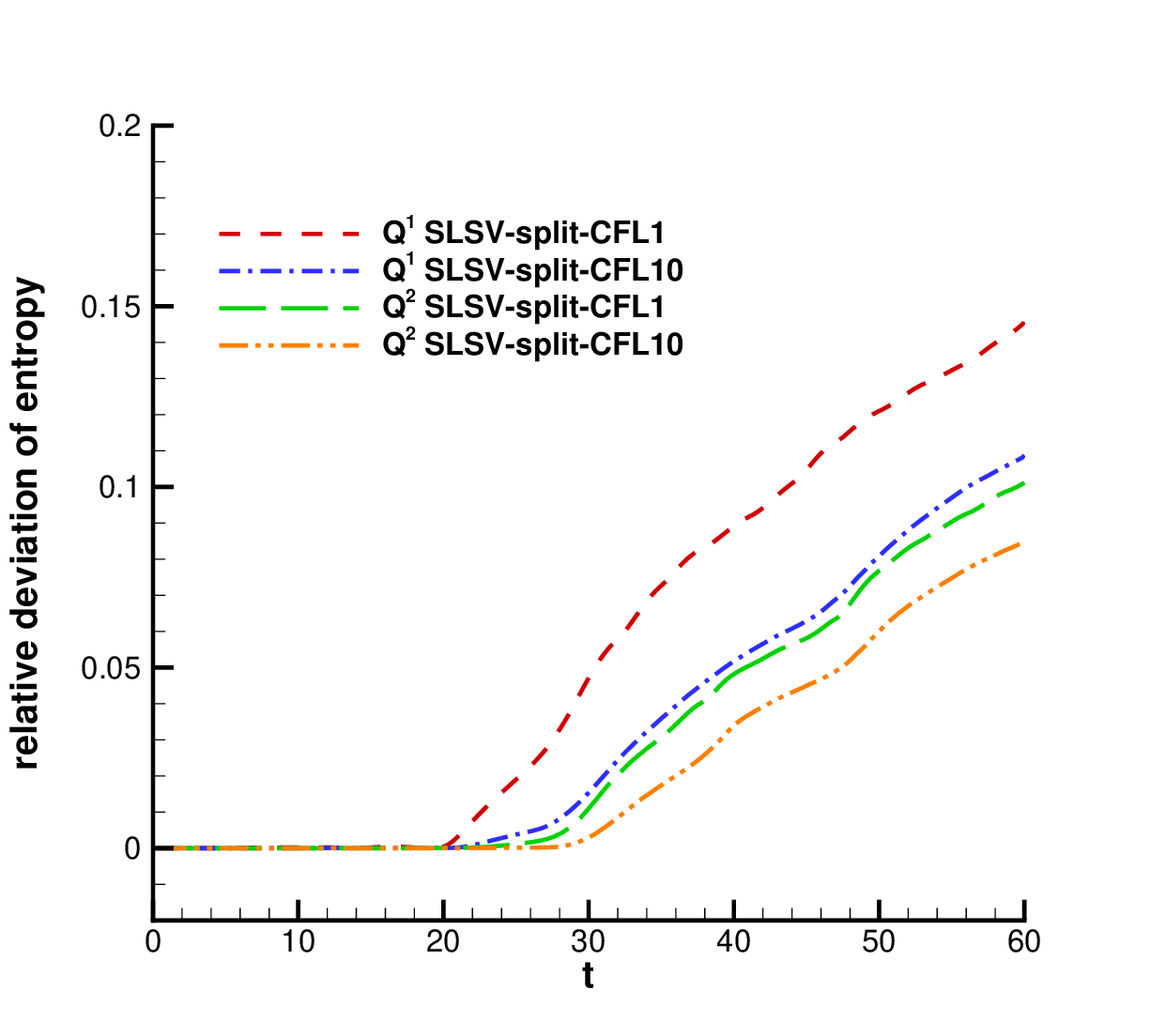}
        \caption{\label{Fig: two stream2 quantities}  Two stream II. Time evolution of the relative deviations of $L^1$ (upper left) and $L^2$ (upper right) norms of the solution as well as the discrete kinetic energy (lower left) and entropy (lower right), using a mesh of $160 \times 160$ elements.}
    \end{figure}

\begin{figure}[!h]
        \centering
        \includegraphics[width=0.37\textwidth]{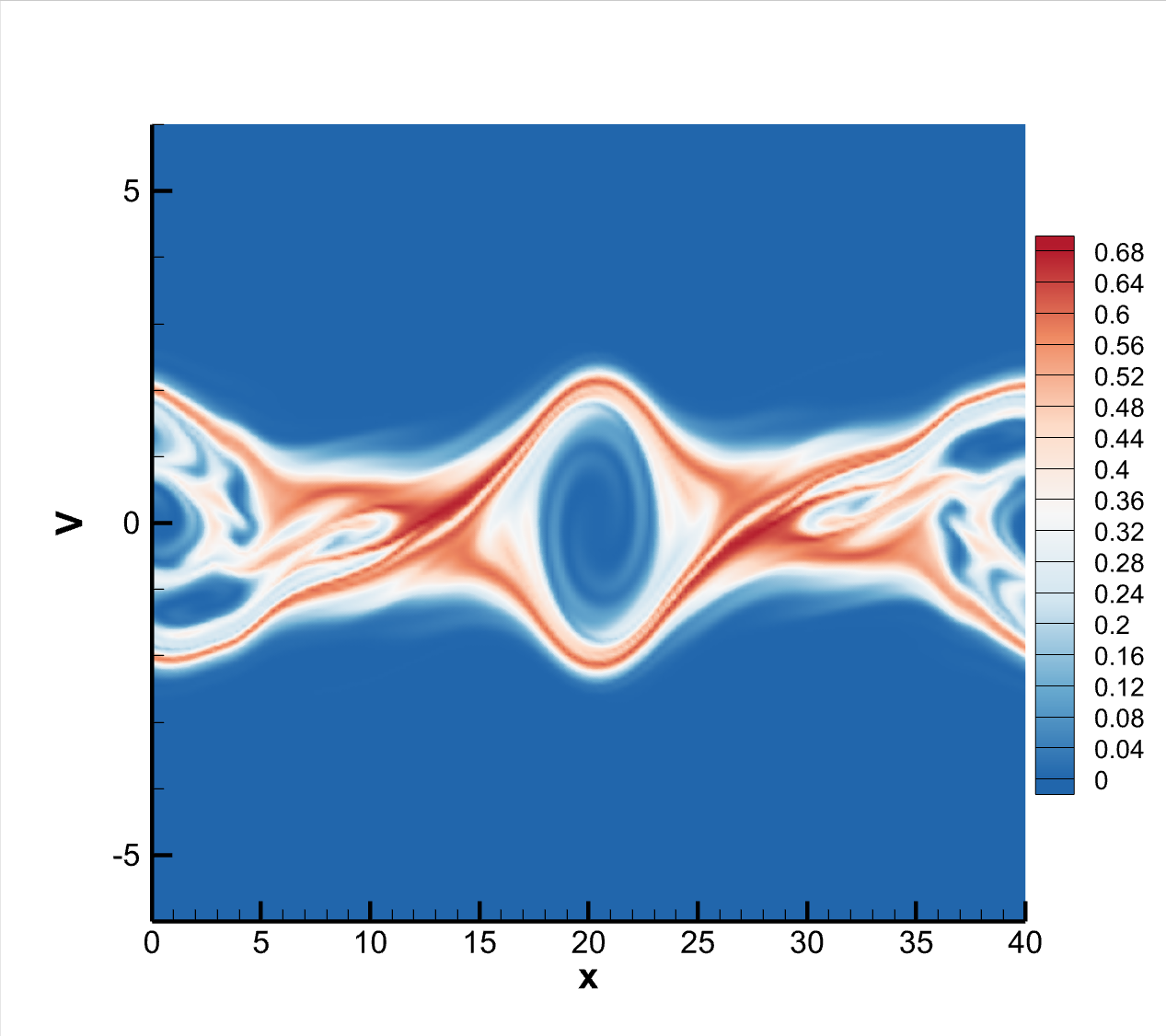}
        \includegraphics[width=0.37\textwidth]{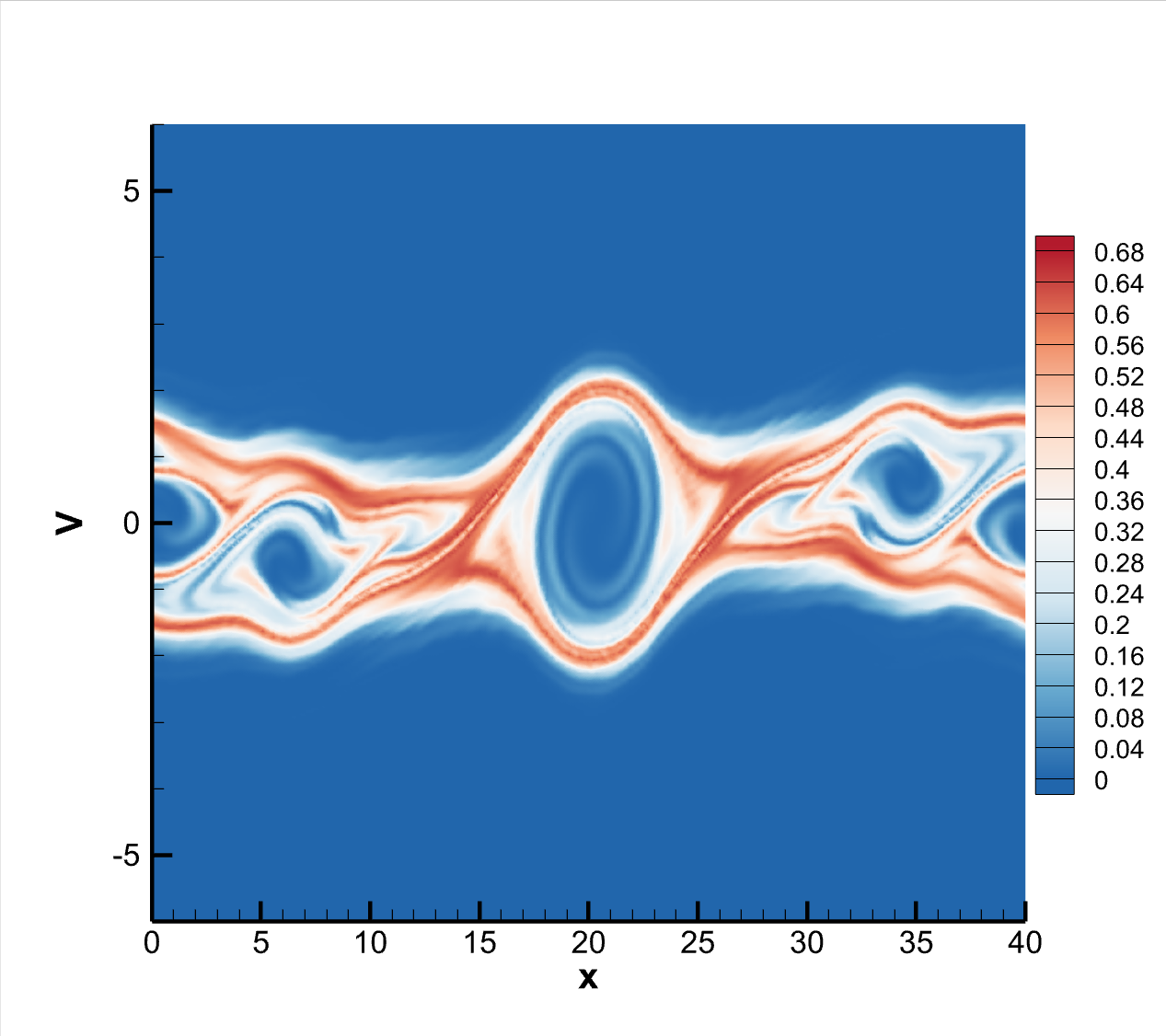}\\
        \includegraphics[width=0.37\textwidth]{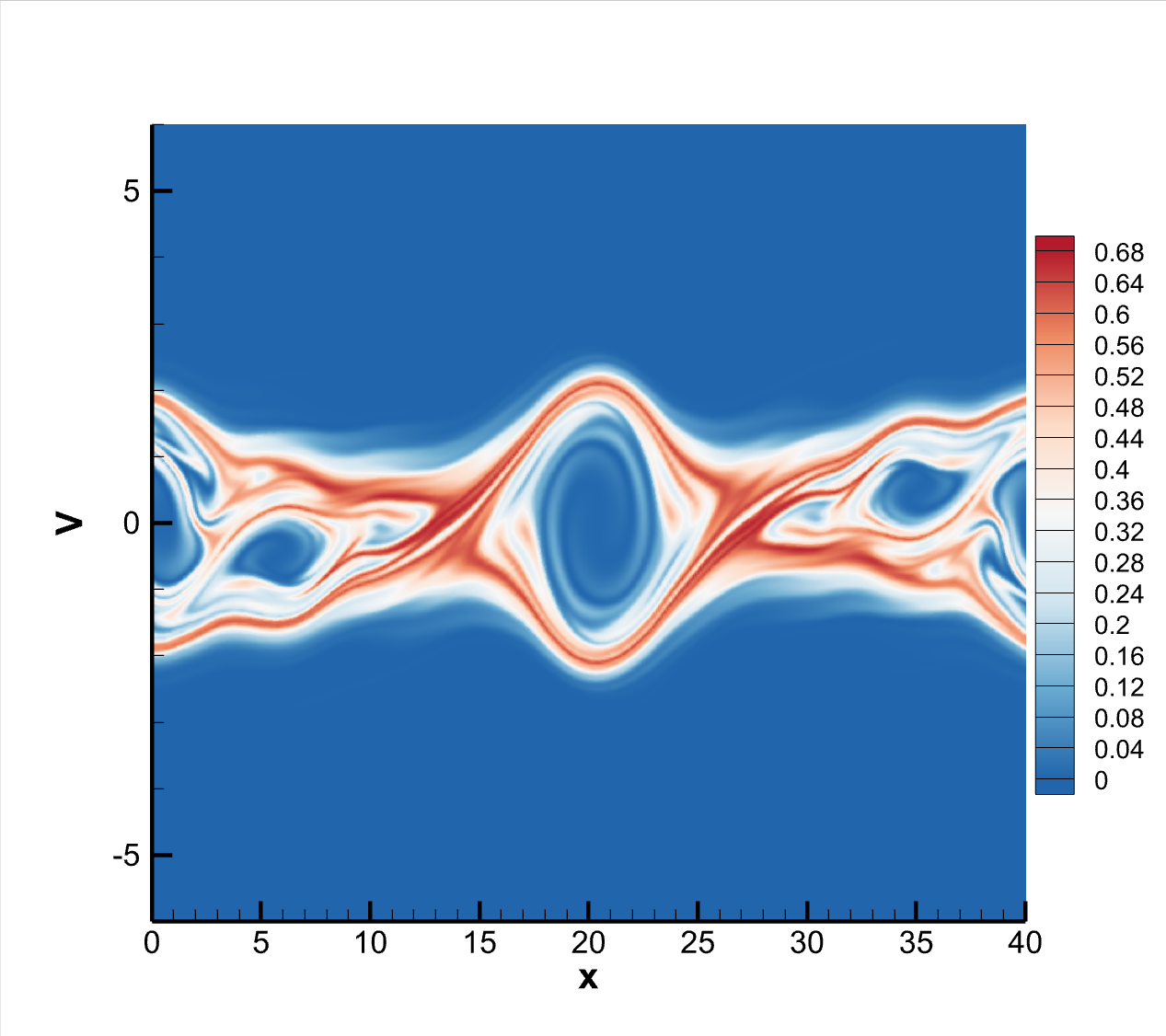}
        \includegraphics[width=0.37\textwidth]{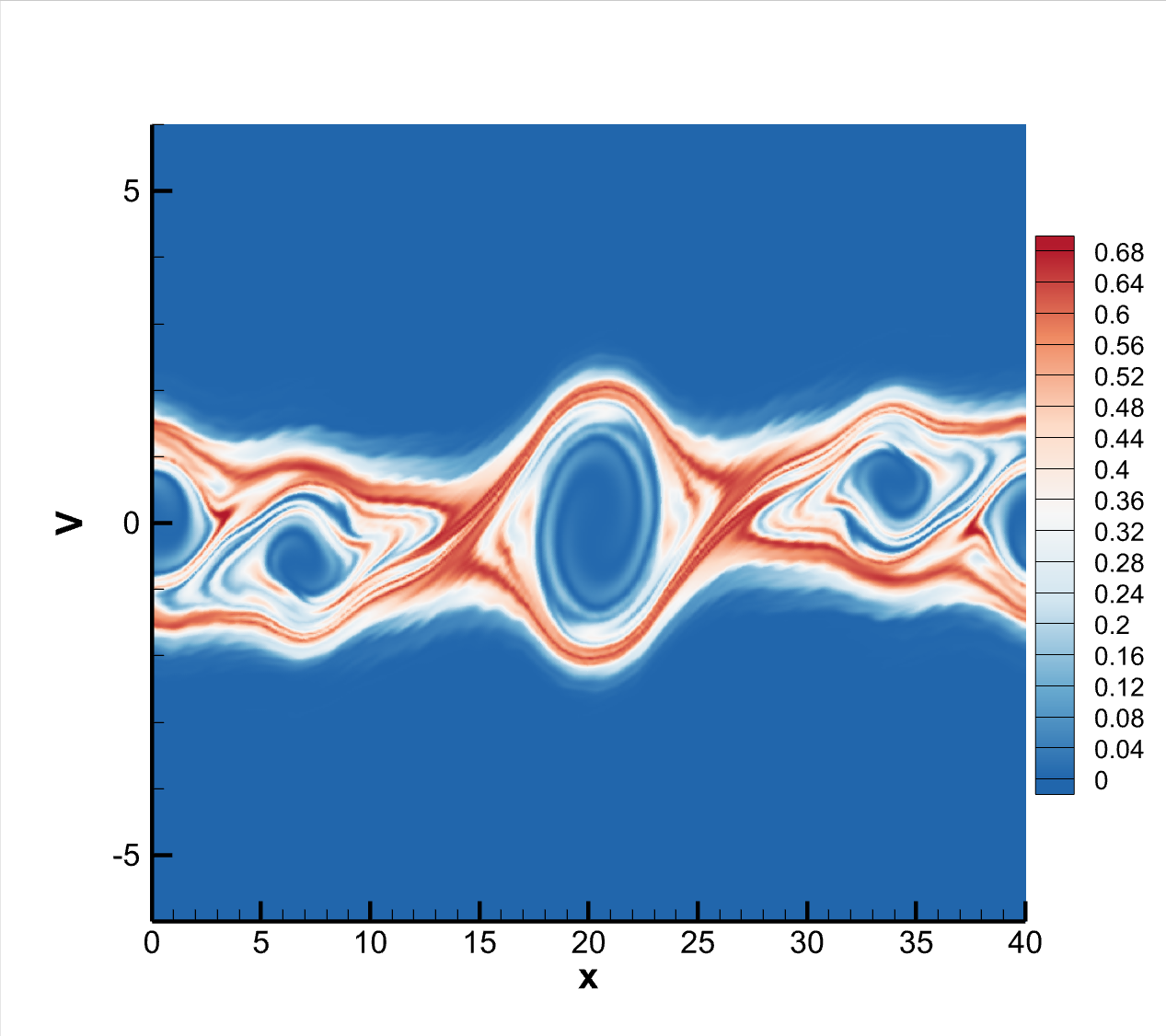}
        \caption{\label{Fig: two stream2 field} Two stream II with the spatial mesh of $160 \times 160$  at $T=70$. (a): ${\mathcal Q}^1$ SLSV with $CFL=10$. (b): ${\mathcal Q}^1$ SLSV with $CFL=20$. (c): ${\mathcal Q}^2$ SLSV with $CFL=10$.(d): ${\mathcal Q}^2$ SLSV with $CFL=20$.}
    \end{figure}

   We plot the time evolution of the electric field in the $L^2$ and $L^{\infty}$ norms (in semi-log scale) in Figure \ref{Fig: two stream2 E2}, which is benchmarked against the results reported in the literature. 
   Time evolution of the relative derivation of the discrete $L^1$ norm, $L^2$ norm, energy and entropy is presented 
   in Figure \ref{Fig: two stream2 quantities}. Figure \ref{Fig: two stream2 field} shows the numerical solutions of phase space profiles at $T=70$, computed by the ${\mathcal Q}^k,k=1,2,$ SLSV method with CFL$=10,20$. All the results show that the SLSV method is accurate for the capture of detail structures.

\section{Conclusion}
In this paper,  a novel high-order semi-Lagrangian spectral volume method based on 
operator splitting framework is proposed  and studied for nonlinear Vlasov–Poisson 
system. 
The proposed SLSV method integrates  both advantages of the semi-Lagrangian and spectral volume approaches
 and has the following desired properties:
 \begin{itemize}
 \item It is a high-order method with arbitrary high-order spatial accuracy. 
 \item It preserves mass and positivity, and  is unconditionally energy stable. 
 \item It is flux-free and feasible for parallel computing, 
 as the high-dimensional original problem is decomposed into one-dimensional advection equations,
 which can be solved in parallel.
\item It has the ability to handle large time steps,  and thus 
predicts more accurate approximation  over long time intervals and more efficient for long time simulation. 
 \end{itemize}
The SLSV-based Vlasov-Poisson solver is applied to classical benchmark problems, including Landau damping and two-stream instabilities, to evaluate its performance in VP simulations. Extensive numerical experiments have been conducted to validate the efficiency, accuracy, and robustness of the proposed scheme, confirming its effectiveness in addressing complex nonlinear phenomena.

\bigskip 

{\bf Funding} This research is supported by National Natural Science Foundation of China under grants No. 12271049, 12201052, and
the Guangdong Provincial Key Laboratory of Inter-disciplinary Research and Application for Data Science of project code 2022B1212010006, BNBU Research Grant with No. of UICR0700035-22 at Beijing Normal-Hong Kong Baptist University,Zhuhai, PR China,
Guangdong basic and applied basic research foundation [2025A1515012182], National Key Laboratory for Computational Physics [6142A05230201].

{\bf Data Availability} Data sharing not applicable to this article as no datasets were generated or analyzed during
the current study.

\section {Declarations} \ \ \newline

{\bf Conflict of interest} The author declares no conflict of interest.

% \bibliographystyle{siamplain}
% \bibliography{references}

\end{document}